\documentclass[10pt,a4paper]{article}

\usepackage{lineno,hyperref}
\modulolinenumbers[5]

\usepackage[english]{babel}
\usepackage{graphicx}
\usepackage{marginnote}

\usepackage{tikz}
\usetikzlibrary{patterns}
\usepackage{bbm}
\usetikzlibrary{arrows}
\usepackage{mathtools}
 \usepackage{multirow}
\usepackage{amssymb}
\usepackage{amsfonts}
\usepackage{subfig}

\usepackage{algorithm}     
\usepackage{algpseudocode} 

\newtheorem{theorem2}{Theorem}[section]
\newtheorem{remark}[theorem2]{Remark}

\newcommand{\email}[1]{\hspace*{\stretch{1}}\emph{\texttt{#1}}}

\makeatletter
\def\blfootnote{\xdef\@thefnmark{$\star$}\@footnotetext}
\makeatother
\newenvironment{Authors}%
  {\begin{center}\begin{bfseries}}%
  {\end{bfseries}\end{center}}
\newenvironment{Addresses}%
  {\begin{flushleft}\begin{itshape}}%
  {\end{itshape}\end{flushleft}}

 \usepackage{fancyhdr}  
  
\fancypagestyle{plain}{
	\fancyhead{}
	\fancyhead[C]{\hfill submitted to Journal of Computational Physics, October 2017}
}

\begin{document}

\thispagestyle{plain}

\title{A Reduced Basis Technique for Long-Time Unsteady Turbulent Flows}
 \date{}
 
 \maketitle
\vspace{-50pt}

 \begin{Authors}
{Lambert Fick}$^{1}$,
{Yvon Maday}$^{2,3}$,
{Anthony T Patera}$^{4}$,
Tommaso Taddei$^{2}$
\end{Authors}

\begin{Addresses}
$^1$   Texas A{\&}M University, 
Department of Nuclear Engineering, USA
 \email{lambert.fick@tamu.edu} \\ 
$^2$   Sorbonne Universit{\'e}s,
Laboratoire Jacques-Louis Lions, France       \email{taddei@ljll.math.upmc.fr,
maday@ann.jussieu.fr} \\ 
$^{3}$ Brown University, Division of Applied Mathematics, USA
       \email{yvon{\_}jean{\_}maday@brown.edu} \\
$^{4}$ MIT, Department of Mechanical Engineering, USA\\
       \email{patera@mit.edu} \\             
\end{Addresses}

\begin{abstract}
We present a reduced basis technique for long-time integration of parametrized incompressible turbulent flows.
The new contributions are threefold. 
First, we propose a constrained Galerkin formulation that corrects the standard Galerkin statement by incorporating prior information about the long-time attractor.
For explicit and semi-implicit time discretizations, our statement reads as  a constrained quadratic programming problem where the objective function is the Euclidean norm of the error in the reduced Galerkin (algebraic) formulation, while the constraints correspond to  bounds for the maximum and minimum value of the coefficients of the $N$-term expansion.
Second, we propose an \emph{a posteriori} error indicator, which corresponds to the dual norm of the residual associated with the time-averaged momentum equation.
We demonstrate that the error indicator
is highly-correlated with the error in mean flow prediction, and 
 can be efficiently computed through an  offline/online strategy.
 Third, we propose a Greedy algorithm for the construction of an approximation space/procedure valid over a range of parameters; the Greedy is informed by the \emph{a posteriori} error indicator developed in this paper. 
We illustrate our approach and we demonstrate its effectiveness 
by studying the dependence of a two-dimensional turbulent lid-driven cavity flow on the Reynolds number.
\end{abstract}


\section{Introduction}
\label{sec:goal_project}
For turbulent flows, estimation of the entire solution trajectory through a low-dimensional Reduced Order Model (ROM) is infeasible due to the slow decay of the Kolmogorov $N$-width, and due to the sensitivity of the dynamical system to perturbations. Nevertheless, it might still be possible to estimate various  moments of the solution associated to a Direct Numerical Simulation (DNS).
The goal of this work is to develop a Reduced Basis (RB) technique for long-time integration of turbulent flows.
Our equations of interest are the unsteady incompressible Navier-Stokes  equations for high-Reynolds number flows with no-slip boundary conditions:
\begin{equation}
\label{eq:unsteady_Navier_Stokes}
\left\{
\begin{array}{ll}
\partial_t u + (u \cdot \nabla) u - \frac{1}{{\rm Re}} \Delta u + \nabla p = f & {\rm in} \, \Omega \times \mathbb{R}_+, \\[3mm]
\nabla \cdot  u  = 0 & {\rm in} \, \Omega \times \mathbb{R}_+, \\[3mm]
u=g  & {\rm on} \, \partial \Omega \times \mathbb{R}_+, \\[3mm]
u=u_0  &   {\rm on} \,   \Omega \times \{ 0 \}, \\[3mm]
\end{array}
\right.
\end{equation}
where $\Omega \subset \mathbb{R}^d$, 
and $f,g, u_0$ are suitable fields.
We denote by $\mu \in \mathcal{P} \subset \mathbb{R}^P$ the set of parameters associated with the equations.

We consider two separate problems: the solution reproduction problem, and the parametric problem.
In the \emph{solution reproduction problem}, given the velocity DNS data 
$\{ u(t_{\rm s}^k; \mu)\}_{k=1}^K$, 
and possibly the pressure DNS data
$\{  p(t_{\rm s}^k;\mu)    \}_{k=1}^K$, 
at the sampling times 
$\{ t_{\rm s}^k \}_{k=1}^K \subset \mathbb{R}_+$,
 we wish to construct a ROM that approximates 
--- in a sense that will be  defined soon ---  
 the original DNS data for the same value of the parameter $\mu$. In the \emph{parametric problem}, we wish to construct a ROM that approximates 
 the DNS data for all values of $\mu$ in a prescribed parameter range $\mathcal{P} \subset \mathbb{R}^P$.
 For the parametric problem, we wish to control the offline costs associated with the construction of the reduced space:
this implicates  a Greedy (rather than POD)  strategy  in parameter.
 Although the solution reproduction problem might be of limited interest in practice, it represents the first step towards the development of a ROM for the  parametric problem.
 
 Following \cite{balajewicz2012stabilization}, we quantify the accuracy of the ROM by computing the error in the long-time average $\langle u \rangle(x) := \lim_{T \to \infty} \frac{1}{T} \int_0^T u(x,t) \, dt$, 
 and the error in the turbulent kinetic energy
${\rm TKE}(t) = \frac{1}{2} \int_{\Omega} \| u(x,t) -  \langle u \rangle(x)   \|_2^2 \, dx$ where 
$\| \cdot \|_2$ is the Euclidean norm.
We remark that, at present, there is no universally-accepted notion of ROM accuracy for turbulent flows.
In \cite{wang2012proper}, the authors evaluate different ROMs based on five different criteria: the kinetic energy spectrum, the mean velocity components, the Reynolds stresses, the root mean square values of the velocity fluctuations, and the time evaluations of the POD coefficients. 
In \cite{cordier2013identification}, the authors consider just the time evaluations and the power spectra of selected POD coefficients.
 From an engineering perspective, the definition of accuracy is entirely determined by the particular quantity of interest we wish to predict: for this reason, we envision that for several applications accurate estimates of the long-time averages and possibly of the turbulent kinetic energy might suffice.

The most popular approach for the solution reproduction problem is the so-called POD-Galerkin  method
\cite{deane1991low,ma2000dynamics,ma2002low,galletti2004low}:
first,
we generate a reduced space $\mathcal{Z}^{\rm u} = \{ \zeta_n \}_{n=1}^{N}$ for the velocity field by applying the Proper Orthogonal Decomposition (POD, \cite{berkooz1993proper,kunisch2002galerkin,kahlbacher2007galerkin}) in the $L^2$ inner product; 
then, we estimate the velocity field for each time-step 
$t_{\rm g}^1,\ldots,t_{\rm g}^J$
as $\hat{u}^j(\cdot)= \sum_{n=1}^{N}  a_n^j \,  \zeta_n (\cdot)$ where the coefficients $ \mathbf{a}^j = [a_1^j,\ldots,a_{N}^j ]$ are computed by projecting the momentum equation onto the space $\mathcal{Z}^{\rm u}$. 
Since all  DNS data for the velocity are divergence-free, it is straighforward to verify that POD modes $\zeta_1,\ldots,\zeta_{N}$ are divergence-free, and so is the ROM solution.

As observed by several authors, ROMs based on $L^2$ POD-Galerkin are prone to instabilities
\cite{carlberg2015galerkin,rempfer2000low}. This can be explained through a physical argument.
In the limit of high-Reynolds numbers, large-scale flow features are broken down into smaller and smaller scales until the scales are fine enough that viscous forces can dissipate their energy  (\cite{moin1998direct}).
This implies that small-scale modes have significant influence on the dynamics. POD modes  based on the $L^2$ inner product are biased toward large, high-energy scales: since large scales are not endowed with the natural energy dissipation tendency  of the smaller lower-energy viscous scales, this leads to instabilities and/or large errors in the estimate of the turbulent kinetic energy.

To address the issue of stability, several strategies have been proposed: (i) including dissipation via a closure model, (ii) modifying the POD basis by including functions that resolve a range of scales, (iii) employing a minimum residual formulation, (iv) employing stabilizing inner products,  (v) calibration methods, and (vi) generating the reduced space through Dynamic Mode Decomposition.
Below, we briefly describe each strategy, and we provide some references;
we  remark that most of the works presented below are restricted to either laminar flows or short-time integration;
therefore, they  do not directly address the problem of interest (the long-time integration of fully-turbulent flows).
We also recall that other topics are treated in the literature:
in particular, Noack et al. \cite{noack2005need} proposed to incorporate pressure in the ROM for
cases with other than no-slip boundary conditions. For the problem considered in this paper, the discussion of  
\cite{noack2005need} is not relevant, and is here omitted.
\begin{itemize}
\item[(i)]
Starting with the pioneering work in \cite{aubry1988dynamics}, several authors have proposed to include dissipation through the vehicle of a closure model.
A first class of models is designed and motivated by analogy with Large Eddy Simulation (see \cite{sagaut2006large} for an introduction to LES): in this respect, 
Couplet et al.  (\cite{couplet2003intermodal}) 
and Noack et al. (\cite[section 4.2]{noack2003hierarchy}) 
 provide theoretical and numerical evidence that the energy transfer among $L^2$-based POD modes is similar to the energy transfer among Fourier modes, 
and for this reason LES ideas based on the energy cascade concept might be promising for POD-ROMs.  
We remark that in \cite{couplet2003intermodal,noack2003hierarchy} the POD space is built 
for the fluctuating field $u_{\rm f} = u - \langle u \rangle$ for a fixed value of the parameters (solution reproduction problem) based on the $L^2$ inner product:
for this reason, it appears difficult to rigorously apply these ideas in the parametric setting
in which we must combine modes associated with different parameters.
Another class of closure models is based on the extension of stabilization techniques originally introduced in the Finite Element or Spectral framework:
two notable examples are the Spectral Vanishing Viscosity Model (SVVM, \cite{sirisup2004spectral}, see also  \cite{maday2016online}) 
originally presented by Tadmor in \cite{tadmor1989convergence} for spectral discretization of nonlinear conservation laws for controlling high-wave number oscillations, and the 
SUPG stabilization discussed in  \cite{giere2015supg}. 
We refer to \cite{wang2012proper} for a numerical comparison of four closure models for incompressible Navier-Stokes equations:
the mixing-length model, the Smagorinsky model, the variational multiscale model, and the dynamic subgrid-scale model.
We further refer to \cite{xie2017approximate} for another POD closure model based on approximate deconvolution, and we  refer to \cite{san2013proper} for a numerical comparison of several closure models for the Burgers' equation.
Finally, we mention the  nonlinear Galerkin method proposed by Marion and Temam in \cite{marion1990nonlinear}, and applied to the simulation of turbulent flows in \cite{debussche1995nonlinear}.
As the above-mentioned variational multiscale method, this approach 
corrects the standard Galerkin model by exploiting 
the separation between large-scale and small-scale modes. To our knowledge, the nonlinear Galerkin method has never been applied in the model order reduction framework.
\item[(ii)]
Another approach is based on including in the POD basis functions that resolve a range of scales. Bergmann et al. in \cite{bergmann2009enablers} 
(see also \cite{bergmann2004optimisation})
proposed to 
augment the original POD basis with a second POD performed on the 
 residuals of the momentum equation (and of the mass equation in case pressure is modelled by the ROM); on the other hand, Balajewicz and Dowell proposed a Greedy technique to  include
in the basis random linear combinations of  low-energy POD modes — associated with the $L^2$ inner product.
\item[(iii)]
Minimum residual formulation was first introduced in the reduced basis framework in \cite{maday2001blackbox} for linear noncoercive problems, and then extended to fluid problems in  \cite{carlberg2011efficient,carlberg2013gnat,carlberg2015galerkin,tallet2015minimum}. Given the reduced space $\mathcal{Z}^{\rm u}$ for velocity (and possibly the reduced space $\mathcal{Z}^{\rm p}$ for pressure), after having discretized the equation in time, the latter approach computes the solution in  $\mathcal{Z}^{\rm u}$ (or $\mathcal{Z}^{\rm u} \times \mathcal{Z}^{\rm p}$) that minimizes a suitable dual residual at each time-step. 
We remark that for problems with quadratic nonlinearities
minimum residual ROMs require $\mathcal{O}(N^4)$ storage  
and the online cost for each time-step is  $\mathcal{O}(N^4)$ for semi-implicit/explicit time-discretizations
--- as opposed to $\mathcal{O}(N^3)$ for standard POD-Galerkin ROMs.
For this reason, hyper-reduction techniques are employed to reduce  the online cost and the memory constraints
(\cite{carlberg2011efficient,carlberg2013gnat,carlberg2015galerkin}).
\item[(iv)]
Iollo et al. in \cite{iollo2000stability} proposed to employ the $H^1$ inner product — rather than the more standard $L^2$ inner product — to generate the POD modes. This choice is motivated by dynamic considerations: since small-scale modes have relatively large $H^1$ norm compared to their $L^2$ norm, and recalling that small scales are responsible for energy dissipation, the use of the $H^1$ inner product leads to a more dissipative reduced order model. We remark that several other authors proposed to not employ the standard $L^2$ inner product (\cite{aubry1992spatio,barone2009stable,rowley2004model}); however, their choices were not motivated by long-time stability considerations.
\item[(v)]
If we denote by 
 $\dot{  \mathbf{a}  } = \mathcal{F}( \mathbf{a} )$
 the ROM for the coefficients of the POD expansion,
in \cite{couplet2005calibrated}, Couplet et al. proposed to calibrate the coefficients of $\mathcal{F}$ 
based on DNS data, under the assumption that $\mathcal{F}$ is a polynomial of degree $2$ in $\mathbf{a}$.
We observe that the ROM $\mathcal{F}$ depends on the particular POD basis selected; for large values of $N$ (dimension of the POD space), the calibration procedure might require a substantial number of DNS snapshots.
\item[(vi)]
Dynamic Mode Decomposition (DMD)
was first proposed by Schmid in \cite{schmid2010dynamic};
as shown by Rowley et al. in \cite{rowley2009spectral} DMD can be interpreted as an algorithm for finding the Koopman modes associated with the nonlinear discrete dynamical system obtained from the discretization of the Navier-Stokes equations. Despite several authors have proven the effectiveness  of DMD for the extraction of physically-relevant time scales and their associated spatial structures (\cite{rowley2009spectral,schmid2010dynamic,schmid2011applications}), the work by Alla and Kutz \cite{alla2016nonlinear}
represents one of the few examples of application of DMD within  the Galerkin framework.
\end{itemize}

Despite these advances, the solution reproduction problem remains an open issue, particularly for turbulent flows.
By performing a detailed analysis of the performance of the POD-Galerkin approach, we empirically demonstrate that in the case of turbulent flows POD-Galerkin  ROMs might exhibit other spurious effects such as false stable steady flows.
This demonstrates the need for a more fundamental correction to the POD-Galerkin formulation.
We remark that a similar issue has been observed in \cite{sirisup2004spectral} by Sirisup and Karniadakis for long-time integration of a POD-Galerkin ROM for a laminar flow past a cylinder,
and  
--- in a different context --- 
by Curry et al. in \cite{curry1984order} for highly-truncated spectral approximations to turbulent flows.

To our knowledge, there are very few works that systematically address  the parametric problem.
Ma and Karniadakis  (\cite{ma2002low}), Galletti et al. (\cite{galletti2004low}),
and Stabile et al. (\cite{stabile2017advances})
developed a reduced order model based on POD-Galerkin for the flow past a cylinder for a wide range of Reynolds numbers in the laminar regime. In these papers, the authors use DNS data for pre-selected Reynolds numbers to generate reduced spaces for velocity (\cite{ma2002low,galletti2004low}),
and for velocity and pressure (\cite{stabile2017advances}). 
The choice of the parameters for which the DNS data are computed is performed \emph{a priori}.
Non-adaptive explorations of the parameter space  typically require a large number of offline evaluations of the Full Order Model (FOM); for this reason, they might not be practical in our context.

The goal of this work is to develop a Model Order Reduction (MOR) procedure for the parametrized incompressible Navier-Stokes equations. The three key pieces of our MOR technique are (i) a reduced formulation for the computation of the reduced-order solution, 
(ii) an \emph{a posteriori} indicator for the error in the prediction of the mean flow,
and (iii) a $H^1$-POD-$h$Greedy  strategy for the construction of the reduced space
informed by the above-mentioned
 \emph{a posteriori} indicator.
\begin{itemize}
\item[(i)]
Our reduced formulation is based on a constrained Galerkin formulation. 
The approach is designed to correct the standard Galerkin formulation, especially for moderate values of $N$.
For explicit and semi-implicit time discretizations the formulation reads as a quadratic programming problem where the objective function corresponds to the Euclidean norm of the error in the reduced Galerkin (algebraic) formulation, while the constraints correspond to bounds for the maximum and minimum value of the coefficients 
$\{ a_n^j \}_{n=1}^N \subset  \mathbb{R}$ of the expansion.
We discuss  an actionable  procedure to estimate the lower and upper bounds associated with each coefficient of the reduced expansion based on DNS data.
\item[(ii)]
Our error indicator corresponds to the dual norm of the residual associated with the time-averaged momentum equation.
Time-averaging is here motivated by the chaotic behavior in time of the velocity field.
We verify that the error indicator can be efficiently computed through an offline/online strategy;
furthermore, we numerically demonstrate that the indicator is highly-correlated with the error in the mean flow prediction:
therefore, it is  well-suited to drive the Greedy procedure for the generation of the ROM. 
\item[(iii)]
As in the seminal work by Haasdonk and Ohlberger \cite{haasdonk2008reduced},
our  POD- $h$Greedy algorithm  
combines POD in time with Greedy in parameter.
The procedure is a simplified version of the  $h$-type Greedy proposed in
\cite{eftang2011hp}.
Given $\mu^1 \in \mathcal{P}$, 
we generate the DNS data for $\mu^1$, we apply POD
--- based on the $H^1$ inner product --- 
 to generate the reduced space 
$\mathcal{Z}_1^{\rm u}$, we
build  the POD-ROM, and we 
evaluate the error indicator  $\Delta_1^{\rm u}(\mu)$
for all $\mu \in \mathcal{P}_{\rm train} \subset \mathcal{P}$.
Then, we select $\mu^{2}$ that maximizes the error estimate $\Delta_1^{\rm u}$ over the training set $\mathcal{P}_{\rm train}$. During the second iteration, 
we perform the same steps as before for $\mu^{2}$ (generation of DNS data, POD, construction of the ROM, estimate of the error).
Finally, we select $\mu^3$ that maximizes  
$\Delta_{1,2}^{\rm u} (\mu):=\min \{  \Delta_1^{\rm u}(\mu),  \Delta_2^{\rm u}(\mu)  \}$
over  $\mathcal{P}_{\rm train}$. We then proceed  to generate $\mu^4,\ldots,\mu^L$.
At the end of the offline stage, the procedure produces $L$ different ROMs; during the online stage, 
given a new value of $\mu \in \mathcal{P}$, we first evaluate the ROMs associated with the 
$n_{\rm cand}$ nearest anchor points, and then we select the ROM that minimizes the error indicator.
\end{itemize}
We observe that in this work we restrict ourselves to linear approximation spaces 
$\mathcal{Z}^{\rm u}$
that do not depend on time: this greatly simplifies the implementation,  and 
reduces the memory constraints for long-time integration. 
We refer to 
\cite{taddei2015reduced,cagniart2016model,ohlberger2013nonlinear}
for MOR strategies based on nonlinear approximation spaces for unsteady problems.
On the other hand, we refer to \cite{urban2012new,yano2014space} for space-time approximations of linear and nonlinear parabolic problems.

The idea of employing a constrained formulation is new in the MOR framework. 
We observe that a constrained formulation has been recently proposed in 
\cite{argaud2016stabilization} in the context of steady-state data assimilation:
as in our work, the constraints in \cite{argaud2016stabilization} provide further information about the solution manifold; however, while in our work the constraints are designed to compensate for the effect of the unmodelled dynamics, in \cite{argaud2016stabilization} the constraints are designed to limit the effect of experimental noise. As opposed to calibration techniques and also stabilized ROMs, the hyper-parameters of the ROM (the lower and upper bounds for the coefficients of the expansion) are here tuned directly through sparse DNS data, for an arbitrary reduced space $\mathcal{Z}^{\rm u}$, without having to evaluate the ROM for several tentative candidates. This feature of the approach greatly simplifies the implementation of the method, and in practice reduces the offline costs.

The time-averaged error indicator is also new. In \cite{haasdonk2008reduced}, the authors employ  a residual estimator that measures the error in the entire trajectory: for turbulent flows, this metric is not appropriate due to the chaotic nature of the dynamical system. This explains the importance of our new error indicator for the problem at hand.

The POD-Greedy algorithm was first proposed  in \cite{haasdonk2008reduced}, and then analyzed in \cite{haasdonk2013convergence}. 
The algorithm in \cite{haasdonk2008reduced} combines data from different parameters to generate  a single reduced space for the entire parameter space $\mathcal{P}$.
On the other hand, in our approach we  build a reduced space for each of them.
Recalling the definitions of \cite{eftang2011hp}, the 
algorithm of \cite{haasdonk2008reduced} corresponds to  a POD-$p$Greedy, 
while our approach corresponds to a POD-$h$Greedy.
For turbulent flows, we empirically show in Appendix \ref{sec:p_refinement} 
 that combining modes associated with different values of the parameters might lead to poor performance. On the other hand, $h$-refinement leads to more accurate and stable ROMs. 

The paper is organized as follows.
In section \ref{sec:lid_driven_cavity_problem}, we introduce the model problem considered in this work.
In section \ref{sec:reproduction_problem}, we consider the solution reproduction problem.
First, we consider  the  POD-Galerkin approach: we introduce the formulation, and   we assess the numerical performance.
Then, we present our  constrained POD-Galerkin approach: as in the previous case, we discuss the formulation, and then we numerically assess the performance.
In section \ref{sec:parametric_problem}, we consider the parametric problem:
first, we present the POD-$h$Greedy approach; 
second, we discuss how to adapt the constrained Galerkin formulation to the parametric setting;
third, we propose the time-averaged error indicator; 
and fourth, we present the numerical assessment.
 In section \ref{sec:conclusions}, we  offer some concluding remarks, and we discuss potential extensions of the current approach.
A number of appendices provide further analysis and numerical investigations: 
in Appendix \ref{sec:lid_driven_appendix} we provide an analysis of the model problem considered;
in Appendix \ref{sec:cv_POD} we discuss the selection of the sampling times $\{ t_{\rm s}^k \}_{k=1}^K$;
in Appendix \ref{sec:ROM_stability} we propose a suitable definition of stability for ROMs;
in Appendix \ref{sec:constrained_stability} we investigate the  robustness of the  constrained formulation proposed in this paper;
 in Appendix \ref{sec:p_refinement} we illustrate the problem of $p$-refinement for the parametric case;
  and in Appendix \ref{sec:error_indicator_offline_online} we describe the offline/online strategy employed to compute the error indicator.

\section{A lid-driven cavity problem}
\label{sec:lid_driven_cavity_problem}
We consider the following unsteady lid-driven cavity problem:  
\begin{subequations}
\label{eq:lid_driven_cavity}
\begin{equation}
\left\{
\begin{array}{ll}
\partial_t u + (u \cdot \nabla) u - \nu({\rm Re}) \Delta u + \nabla p = 0 & {\rm in} \, \Omega \times \mathbb{R}_+ ,\\[3mm]
\nabla \cdot  u  = 0 & {\rm in} \, \Omega \times \mathbb{R}_+, \\[3mm]
u=g(x)  & {\rm on} \, \Gamma_{\rm top} \times \mathbb{R}_+, \\[3mm]
u=0 &   {\rm on} \,  \partial \Omega \setminus \Gamma_{\rm top} \times \mathbb{R}_+, \\[3mm]
u=0  &   {\rm on} \,   \Omega \times \{ 0 \}, \\[3mm]
\end{array}
\right.
\end{equation}
where 
the velocity $u:\Omega \times \mathbb{R}_+ \to \mathbb{R}^2$ is a two-dimensional vector field,
the pressure $p: \Omega \times \mathbb{R}_+ \to \mathbb{R}$ is a scalar field,
$\nu({\rm Re}) = \frac{1}{{\rm Re}}$,
$\Omega=(-1,1)^2$, 
$\Gamma_{\rm top}  = \{ x \in \bar{\Omega}: x_2=1 \}$,  the Dirichlet datum is given by
\begin{equation}
g(x) =
\left[
\begin{array}{c}
(1-x_1^2)^2 \\
0 \\
\end{array}
\right],
\end{equation}
and the Laplacian $\Delta$ should be interpreted as component-wise.
We remark that in \eqref{eq:lid_driven_cavity} time is non-dimensionalized by the convective scaling 
(i.e., dimensional boxside half-length divided by dimensional maximum lid velocity).
\end{subequations}
The problem corresponds to a isothermal, incompressible, two-dimensional flow inside a square cavity driven by a prescribed lid velocity.
The problem is a well-known prototypical example used to validate numerical schemes and reduced order models (\cite{balajewicz2012stabilization,balajewicz2012new,cazemier1998proper,shankar2000fluid,terragni2011local,lorenzi2016pod});
unlike in the more standard lid-driven cavity problem with $g(x) = [1,0]$,
here we regularize the singularity near the upper corners of the cavity.

In this paper, we study the dependence of the flow on the Reynolds number, that is $\mu = {\rm Re}$.
It is well-known (\cite{shankar2000fluid}) that the flow
exhibits a long-time unsteady but stationary solution   for ${\rm Re} > {\rm Re_c}$  (\cite{shankar2000fluid});
here stationarity implies that all statistics are invariant under a shift in time (\cite{pope2001turbulent}).
Since we are interested in long-time unsteady flows, we here consider 
${\rm Re} \in \mathcal{P} = 
[ {\rm Re}_{\rm LB}, {\rm Re}_{\rm UB}   ]
= [15000,25000]$:
for all values of ${\rm Re}$ in $\mathcal{P}$ the flow is asymptotically statistically stationary.
Balajewicz and Dowell considered  the same problem --- 
for a single value of $\rm Re$ --- 
in
\cite{balajewicz2012stabilization}; we remark that they define the viscosity as 
$\nu({\rm Re}) = \frac{2}{{\rm Re}}$, and they consider the case ${\rm Re} = 30000$.

 In view of the development of the ROM for \eqref{eq:lid_driven_cavity} it is convenient to consider the lifted equations.
If we denote by $R_g$ the two-dimensional vector field defined as the solution to the following Stokes problem:
\begin{subequations}
\label{eq:navier_stokes_unsteady_lifted}
\begin{equation}
\label{eq:stokes_lift}
\left\{
\begin{array}{ll}
-\Delta R_g
+
\nabla \lambda
= 0
&
{\rm in}
\,
\Omega,
\\[3mm]
\nabla \cdot R_g = 0 
&
{\rm in}
\,
\Omega,
\\[3mm]
 R_g=g
&
{\rm on}
\,
\Gamma_{\rm top},
\\[3mm]
 R_g=0
&
{\rm on}
\,
\partial \Omega \setminus
\Gamma_{\rm top},
\\
\end{array}
\right.
\end{equation}
we can define the lifted velocity solution $\mathring{u}= u - R_g$ as the solution to:
\begin{equation}
\left\{
\begin{array}{ll}
\partial_t \mathring{u}
+
\left( (\mathring{u} + R_g) \cdot \nabla \right)  
(\mathring{u} + R_g)
- \frac{1}{\rm Re} \Delta \left(
\mathring{u}
+
R_g\right)
+ \nabla p = 
0
&
{\rm in}
\,
\Omega \times \mathbb{R}_+,
\\[3mm]
\nabla \cdot \mathring{u} = 0
&
{\rm in}
\,
\Omega  
\times \mathbb{R}_+,
\\[3mm]
\mathring{u} =0
&
{\rm on}
\,
\partial \Omega \times \mathbb{R}_+,
\\[3mm]
\mathring{u}(t=0)= -R_g
&
 {\rm on} \,   \Omega \times \{ 0 \}.
\\[3mm]
\end{array}
\right.
\end{equation}
\end{subequations}
Then, if we introduce the spaces $V:= [H_0^1(\Omega)]^2$, and 
$Q= \{ q \in L^2(\Omega): \int_{\Omega} q \, dx=0 \}$, 
we can define the weak form associated with \eqref{eq:navier_stokes_unsteady_lifted}:
find $(\mathring{u},p) \in \mathcal{V} \times \mathcal{Q}$ such that
for a.e. $t>0$ 
\begin{subequations}
\label{eq:weak_form_cavity_problem}
\begin{equation}
\left\{
\begin{array}{l}
\langle \partial_t \mathring{u}(t), v \rangle_{\star}
+
\frac{1}{\rm Re}
(\mathring{u}(t) + R_g, v )_V
+
c(\mathring{u}(t) + R_g,  \mathring{u}(t)+ R_g, v)
+
b(v, p(t))
=
0
\\[2mm]
\hfill
\forall \, v \in V ,
\\[3mm]
b(\mathring{u}(t), q) = 0
\hfill
\forall \, q \in Q, \\
\end{array}
\right.
\end{equation}
where
$ \mathcal{V} = 
\{ 
v \in L_{\rm loc}^2(\mathbb{R}_+; V):
\, \partial_t v \in L_{\rm loc}^2(\mathbb{R}_+; V')
 \}$,  $ \mathcal{Q} =  L_{\rm loc}^2(\mathbb{R}_+; Q)$,
 $\langle \cdot, \cdot \rangle_{\star}$ denotes the pairing between $V'$ and $V$
which (for our smoothness assumptions and numerical approximations)   can be evaluated in terms of the pivot space $L^2$,
 $(w,v)_V = \int_{\Omega} \, \nabla \, w  : \nabla  v \, dx$
 is the inner product associated with $V$, 
 and
\begin{equation}
c(w,u,v) = \int_{\Omega} \, (w \cdot \nabla) u \cdot v \, dx,
\quad
b(v,q) = -\int_{\Omega} \, (\nabla \cdot v) q \, dx.
\end{equation}
\end{subequations}

We resort to a $\mathbb{Q}_M-\mathbb{Q}_{M-2}$ spectral element (\cite{patera1984spectral})
discretization in space, and to an explicit three-step Adams-Bashforth (AB3)/ implicit two-step Adams-Moulton (AM2)  discretization in time.
DNS simulations are performed using the open-source software \texttt{nek5000} (\cite{nek5000-web-page}).
 We refer to the spectral element literature (see, e.g.,  \cite{bernardi1997spectral,canuto2012spectral,karniadakis2013spectral,deville2002high})
for further details about the spectral element method and its implementation for fluid dynamics problems.
More in detail, we consider a $16$ by $16$ structured quadrilateral mesh,  we consider $M=8$, and we resort to an
equispaced time grid $\{ t_{\rm g}^j= j \Delta t \}_{j=0}^J$, with $\Delta t= 5 \cdot 10^{-3}$.
We estimate the  long-time averaged velocity field  as\footnote{In the current implementation, 
$\langle u \rangle_{\rm g}$ is computed inside the time integration loop of the Full Order Model.}:
\begin{equation}
\label{eq:mean_flow_def}
\langle u \rangle_{\rm g}
=
\frac{\Delta t}{T-T_0}
\sum_{j=J_0 + 1}^J u^j,
\end{equation}
where $ T_0=500$, $ T=t_{\rm g}^{J}$, and 
$J_0$ is such that $t_{\rm g}^{J_0} = T_0=500$.   Consequently, we estimate the instantaneous turbulent kinetic energy as
\begin{equation}
\label{eq:TKEj}
{\rm TKE}^j:=
\frac{1}{2} \int_{\Omega}
\,
\| u^j  - \langle u \rangle_{\rm g} \|_2^2 \, dx.
\end{equation}
In Appendix \ref{sec:lid_driven_appendix}, we provide a detailed analysis of the solution to the lid-driven cavity problem  \eqref{eq:lid_driven_cavity}.

In order to generate (and, later,   assess) the ROM, we collect data at the sampling times
$\{ t_{\rm s}^k= T_0 + 
\Delta t_{\rm s}
k  \}_{k=1}^K$ with $\Delta t_{\rm s} = 1$.
We observe that $\{ t_{\rm s}^k   \}_{k=1}^K \subset \{ t_{\rm g}^j  \}_{j=J_0}^J$, and $K \ll J$: this is dictated by memory constraints.
We further observe that we do not collect data in the transient region: this is motivated by the fact that we are here ultimately interested in the long-time dynamics.
In the remainder of the paper, we use the subscript  ``s" to indicate the sampling times, and the subscript ``g" to indicate the time discretization. 
Furthermore, we use the symbol  $\langle \cdot \rangle_{\rm s}$ to indicate time averages performed based on the sampling times, and the symbol 
$\langle \cdot \rangle_{\rm g}$ to indicate time averages performed based on the time grid $\{ t_{\rm g}^j  \}_{j=J_0}^J$.
In Appendix \ref{sec:cv_POD}, we comment on the choice of $\Delta t_{\rm s}$ and $K$.

\section{The solution reproduction problem}
\label{sec:reproduction_problem}
In this section, we propose a MOR procedure for the solution reproduction problem.
As explained in the introduction, the solution reproduction problem is of limited practical interest; however, it represents a key intermediate step towards the development of a MOR procedure for the parametric problem.
Algorithm \ref{reproduction_problem} outlines the general offline/online paradigm for the solution reproduction problem. We recall that the offline stage  is expected to be expensive  and is performed once, while the online stage
should be inexpensive and
is performed many times
---
this distinction is of little relevance here, but will be crucial in section \ref{sec:parametric_problem} for the parametric problem.

\begin{algorithm}[H]                      
\caption{Offline/online paradigm for the solution reproduction problem}     
\label{reproduction_problem}                           

\begin{flushleft}
\textbf{Task:}
find an estimate of $\mathring{u}=\mathring{u}(x,t)$ of the form
$\hat{u}(x,t) = \sum_{n=1}^N a_n(t) \, \zeta_n(x)$.
\end{flushleft}

\textbf{Offline stage}
\medskip

\begin{algorithmic}[1]

\State
Generate the DNS data $\{ \mathring{u}^k := \mathring{u}(t_{\rm s}^k) \}_{k=1}^K \subset V$.
\vspace{4pt}

\State
Generate the reduced space $\mathcal{Z}^{\rm u} = {\rm span} \{  \zeta_n \}_{n=1}^N$.
\vspace{4pt}

\State
Formulate the Reduced Order Model.
\vspace{4pt}
\end{algorithmic}

\medskip

\textbf{Online stage}
\medskip

\begin{algorithmic}[1]

\State
Estimate the coefficients $\{ a_n^j= a_n(t_{\rm g}^j)  \}_{n=1}^N$
for $j=0,1,\ldots,J$.
\vspace{4pt}

\State
Compute the QOIs 
(e.g., mean flow, TKE,...)
\end{algorithmic}
\end{algorithm}

As anticipated in section \ref{sec:lid_driven_cavity_problem}, 
we here  generate a ROM for the lifted velocity field 
$\mathring{u} =  u - R_g$, where $R_g$ is the solution to the Stokes problem
\eqref{eq:stokes_lift}.
Reduction of the lifted equations is preferable from the MOR perspective  since 
it  greatly simplifies the imposition of essential  (Dirichlet) inhomogenous boundary conditions. 
We observe that in the Fluid Mechanics literature many authors consider $R_g = \langle u \rangle_{\rm g}$; however, the  latter choice of the lift cannot be  extended to the parametric case.

This section is organized as follows.
In section \ref{sec:galerkin_ROM}, we present the POD-Galerkin ROM.
We first introduce the formulation, we review Proper Orthogonal Decomposition for the generation of the reduced space,
and then we present numerical results that highlight the limitations of the approach. 
In section \ref{sec:cgalerkin_ROM}, we present the constrained POD-Galerkin ROM proposed in this paper.
As for POD-Galerkin, we first present and motivate the mathematical statement, and then we present a number of numerical results to motivate the approach.

\subsection{The POD-Galerkin ROM}
\label{sec:galerkin_ROM}

\subsubsection{The Galerkin formulation}
\label{sec:galerkin_ROM_formulation}

Given the reduced space
$\mathcal{Z}^{\rm u} = {\rm span} \{  \zeta_n \}_{n=1}^N \subset V_{\rm div} = \{ v \in V: \nabla \cdot v = 0 \}$, 
we seek $\hat{u} \in \mathcal{V}_N:= H_{\rm loc}^1(\mathbb{R}_+; \mathcal{Z}^{\rm u})$ such that
\begin{equation}
\label{eq:galerkin_rom_non_discretized}
\left\{
\begin{array}{l}
\displaystyle{
\frac{d}{dt}
(  \hat{u}(t), v )_{L^2(\Omega)}
+
\frac{1}{\rm Re}
(\hat{u}(t) + R_g, v)_V
+
c(\hat{u}(t) + R_g, \hat{u}(t) + R_g, v)
=
0
}
\\
\hfill
\forall \, v \in \mathcal{Z}^{\rm u},
\\[3mm]
\hat{u}(0)
=
-\Pi_{\mathcal{Z}^{\rm u}}^{L^2} R_g ,
\\
\end{array}
\right.
\end{equation}
where $\Pi_{\mathcal{Z}^{\rm u}}^{L^2} : [L^2(\Omega)]^2 \to \mathcal{Z}^{\rm u}$ is the $L^2(\Omega)$-projection operator on $\mathcal{Z}^{\rm u}$,
and
$(\cdot, \cdot)_{L^2(\Omega)}$ is the $L^2(\Omega)$ inner product.
If we employ a semi-implicit time discretization, we obtain: 
 \begin{equation}
\label{eq:galerkin_rom_discretized}
\begin{array}{l}
\displaystyle{
\left(
\frac{\hat{u}^{j+1}-\hat{u}^j}{\Delta t},
v
\right)_{L^2(\Omega)}
+
\frac{1}{\rm Re}
(\hat{u}^{j+1} + R_g, v    )_V
+
c(\hat{u}^{j} + R_g, \hat{u}^{j+1} + R_g, v)
 =
0
}
\\
\hfill
\forall \, v \in \mathcal{Z}^{\rm u}, 
 \; \;
 j=0,1,\ldots,
 \\
\end{array}
\end{equation}
where $\Delta t = t_{\rm g}^{j+1} - t_{\rm g}^{j}$.
We remark that the time scheme is not the same used by \texttt{nek5000}.
The Galerkin formulation \eqref{eq:galerkin_rom_discretized} leads to the following  algebraic system  for the coefficients
$\{ \mathbf{a}^j \}_{j=0}^J$ of the $N$-term expansion: 
\begin{subequations}
\label{eq:galerkin_algebraic}
\begin{equation}
\mathbb{A}(\mathbf{a}^j; {\rm Re})
\, 
\mathbf{a}^{j+1}
=
\mathbf{F}(\mathbf{a}^j; {\rm Re}),
\qquad
 j=0,1,\ldots,
\end{equation}
where 
$
\mathbb{A}(\mathbf{a}^j; {\rm Re}) := 
\mathbb{A}_1
+
\frac{1}{\rm Re}
\mathbb{A}_2
+
\mathbb{C}(\mathbf{a}^j)$, 
$\mathbf{F}(\mathbf{a}^j; {\rm Re}):=   
\mathbb{E}    \mathbf{a}^j
-
\frac{1}{\rm Re}
\mathbf{G}$,  with 
\begin{equation}
\begin{array}{l}
\displaystyle{
(\mathbb{A}_1)_{m,n}
=
\frac{1}{\Delta t}
( \zeta_n,\zeta_m  )_{L^2(\Omega)}
+
c(R_g,\zeta_n,\zeta_m),
\quad
(\mathbb{A}_2)_{m,n}
=
( \zeta_n, \zeta_m)_V,
}
\\[3mm]
(\mathbb{C}(\mathbf{w}) )_{m,n}
=
\sum_{i=1}^N \, w_i \, c(\zeta_i, \zeta_n,\zeta_m),
\\
\end{array}
\end{equation}
and
\begin{equation}
\mathbf{G}_m 
=
( R_g, \zeta_m)_V,
\quad
(\mathbb{E} )_{m,n}
=
\frac{1}{\Delta t}
( \zeta_n,\zeta_m  )_{L^2(\Omega)} 
-
c(\zeta_n, R_g, \zeta_m).
\end{equation}
\end{subequations}

We observe that the Galerkin model for the velocity field does not contain the pressure field. 
This follows from 
(i)
the fact that the ROM is derived from the weak form of the equations,
(ii) the particular boundary conditions prescribed, and
(iii) the absence of parameters in the form $b(\cdot, \cdot)$.  
We have indeed that for certain choices of the boundary conditions the ROM should be obtained from the strong form of the Navier-Stokes equations:
in this respect, we recall that in \cite{noack2005need}  a Galerkin ROM is derived from the strong form for a laminar flow problem with 
convective boundary condition (\cite{sohankar1998low}) at the outflow. 
In the parametric case it is possible to derive a ROM that does not contain the pressure field 
if the form $b(\cdot, \cdot)$ in \eqref{eq:weak_form_cavity_problem} is parameter-independent; otherwise, it is not possible in general to generate a space $\mathcal{Z}^{\rm u}$ such that 
$b(z, \cdot) \equiv 0$ for all $z \in \mathcal{Z}^{\rm u}$ and for all values of the parameters.  
Since the bilinear form $b(\cdot, \cdot)$ in \eqref{eq:navier_stokes_unsteady_lifted} does not depend on the Reynolds number,  we will be able in section \ref{sec:parametric_problem} to generate a ROM for the velocity only. 
We remark that the case of parametrized $b$ form corresponds to the case of geometric parametrizations, which is of particular interest for applications.
A potential strategy to handle this issue is to resort to the Piola's transform
(see \cite{lovgren2006reduced}).
We  refer to a future work for a detailed discussion of this case.
We also refer to the Reduced Basis literature (\cite{rozza2007stability,rozza2013reduced,lassila2014model,ballarin2017reduced}) for a thorough discussion about fluid problems in parametrized domains for low-to-moderate Reynolds number flows.
 
The algebraic formulation \eqref{eq:galerkin_algebraic} is the starting point for the development of the offline/online decomposition. The matrices $\mathbb{A}_1$, $\mathbb{A}_2$, $\mathbb{E}$, the third-order tensor $\mathbb{C}$ and the vector $\mathbf{G}$ can be pre-computed during the offline stage. Therefore, during the online stage, the method only requires 
$\mathcal{O}(N^3)$  storage, and the online cost is $\mathcal{O}(N^3 J)$. Provided that  $N$ is much smaller than the spatial mesh-size $\mathcal{N}$, the Galerkin ROM is significantly less expensive and less memory-demanding than the Full Order Model.
Other choices of the time discretization lead to similar reduced systems that allow the same offline/online decomposition.
 
\subsubsection{Construction of the reduced space: Proper Orthogonal Decomposition}
\label{sec:galerkin_POD}

We employ Proper Orthogonal Decomposition 
(POD, \cite{berkooz1993proper,kunisch2002galerkin,sirovich1987turbulence}) to generate the reduced space $\mathcal{Z}^{\rm u}$.
Below we briefly review the numerical strategy --- known as \emph{method of snapshots} (\cite{sirovich1987turbulence})
--- employed for the computation of the POD modes.
We refer to \cite{volkwein2011model} for a review of the theoretical results concerning the optimality properties of POD.

Given the snapshot set $\{ \mathring{u}^k \}_{k=1}^K$, we assemble the Gramian $\mathbb{U}  \in \mathbb{R}^{K,K}$
$\mathbb{U}_{k,k'} = ( \mathring{u}^k, \mathring{u}^{k'})_{\star}$
where   $( \cdot, \cdot )_{\star}$ is a suitable inner product that will be introduced soon;
then, we compute the first $N$ eigenmodes of the symmetric matrix $\mathbb{U}$:
\begin{subequations}
\label{eq:method_snapshots}
\begin{equation}
\mathbb{U} \boldsymbol{\zeta}_n = \lambda_n  \, \boldsymbol{\zeta}_n,
\qquad
\lambda_1 \geq \ldots \geq \lambda_K \geq 0;
\end{equation}
finally, we define the POD modes as
\begin{equation}
\zeta_n := \sum_{k=1}^K \, \left(  \boldsymbol{\zeta}_n    \right)_k \, \mathring{u}^k,
\qquad
n=1,\ldots,N.
\end{equation}
\end{subequations}
It is easy to show that $\zeta_1,\ldots,\zeta_N$ can be chosen to be orthogonal in the  $( \cdot, \cdot )_{\star}$-inner product; for stability reasons, we also orthonormalize the POD modes so that $(\zeta_n,\zeta_{n'})_{\star} = \delta_{n,n'}$, $n,n' = 1,\ldots,N$.  In
Appendix  \ref{sec:cv_POD}, we discuss the choice of the sampling times $\{  t_{\rm s}^k\}_k$, and we propose a numerical   technique to  assess  the accuracy of the POD space for the full trajectory.

In this work, we employ the $H_0^1(\Omega)$ inner product:
\begin{equation}
\label{eq:inner_product_POD}
(w,v)_{\star}
=
(w,v)_V
=
\int_{\Omega} \, \nabla w  : \nabla  v \, dx.
\end{equation}
As explained in the introduction, 
this choice is motivated by dynamic considerations. Since small-scale modes have relatively large $H^1$ norm compared to their $L^2$ norm, and recalling that small scales are responsible for energy dissipation, the use of the $H_0^1$ inner product leads to a more dissipative reduced order model (\cite{iollo2000stability}).

\subsubsection{Performance of the POD-Galerkin ROM}
\label{sec:galerkin_performance}

We assess the numerical performance of the POD-Galerkin ROM presented in this section.
We here consider the lid-driven cavity problem \eqref{eq:lid_driven_cavity} for ${\rm Re}=15000$.
We consider the time grid $\{ t_{\rm g}^j = \Delta t j \}_{j=0}^J$ with 
$\Delta t=5 \cdot 10^{-3}$ and $J=2 \cdot 10^5$ ($T=t_{\rm g}^J=10^3$), and we acquire the snapshots 
$\{ \mathring{u}^k =  \mathring{u}(t_{\rm s}^k) \}_{k=1}^K$ where 
$t_{\rm s}^k= 500 + k$ and $K=500$.
The long-time averaged velocity field
$\langle u \rangle_{\rm g}$
is estimated through \eqref{eq:mean_flow_def}.
On the other hand, we estimate the mean TKE as follows:
$$
\langle {\rm TKE} \rangle_{\rm s}
=
\frac{1}{2K}
\sum_{k=1}^K
\, 
\| u^k - \langle u \rangle_{\rm g} \|_{L^2(\Omega)}^2.
$$
Assembling and integration of the Reduced Order Model are performed in \texttt{Matlab}  \cite{MATLAB:2016}.

\begin{figure}[h!]
\centering
\subfloat[ ]  
{  \includegraphics[width=0.30\textwidth]
 {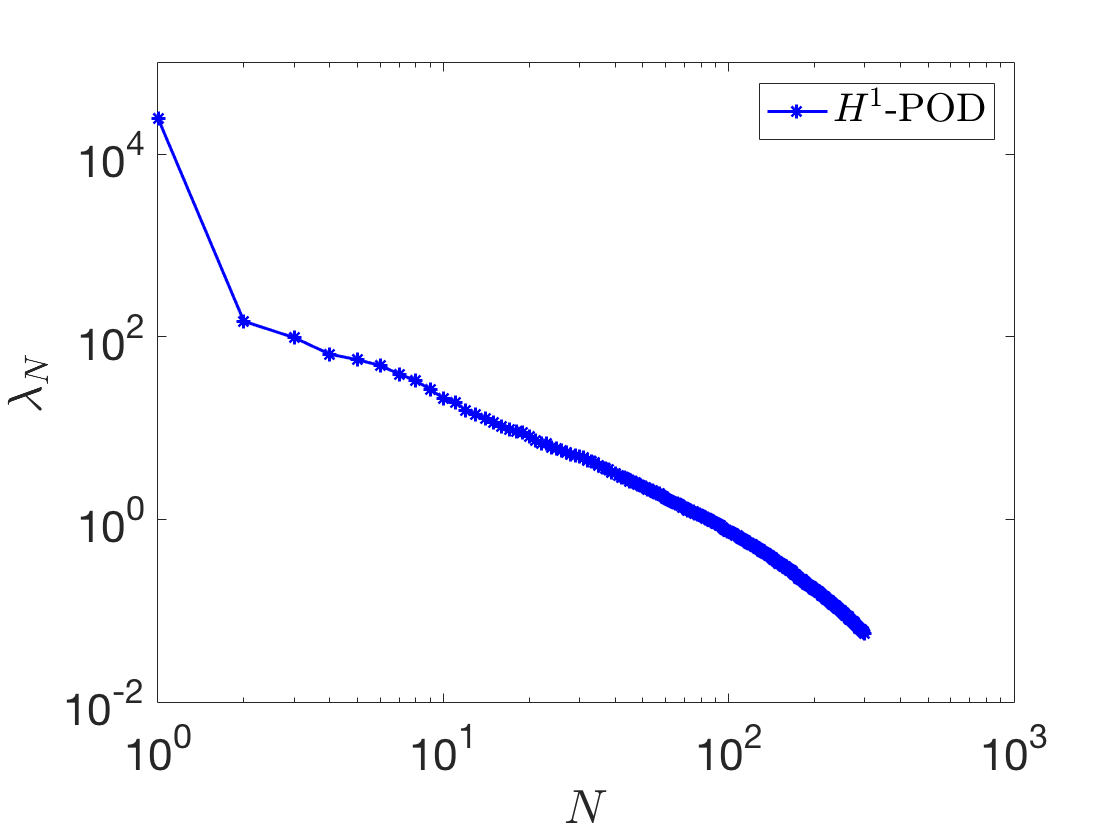}}
 ~~
 \subfloat[]
{  \includegraphics[width=0.30\textwidth]
 {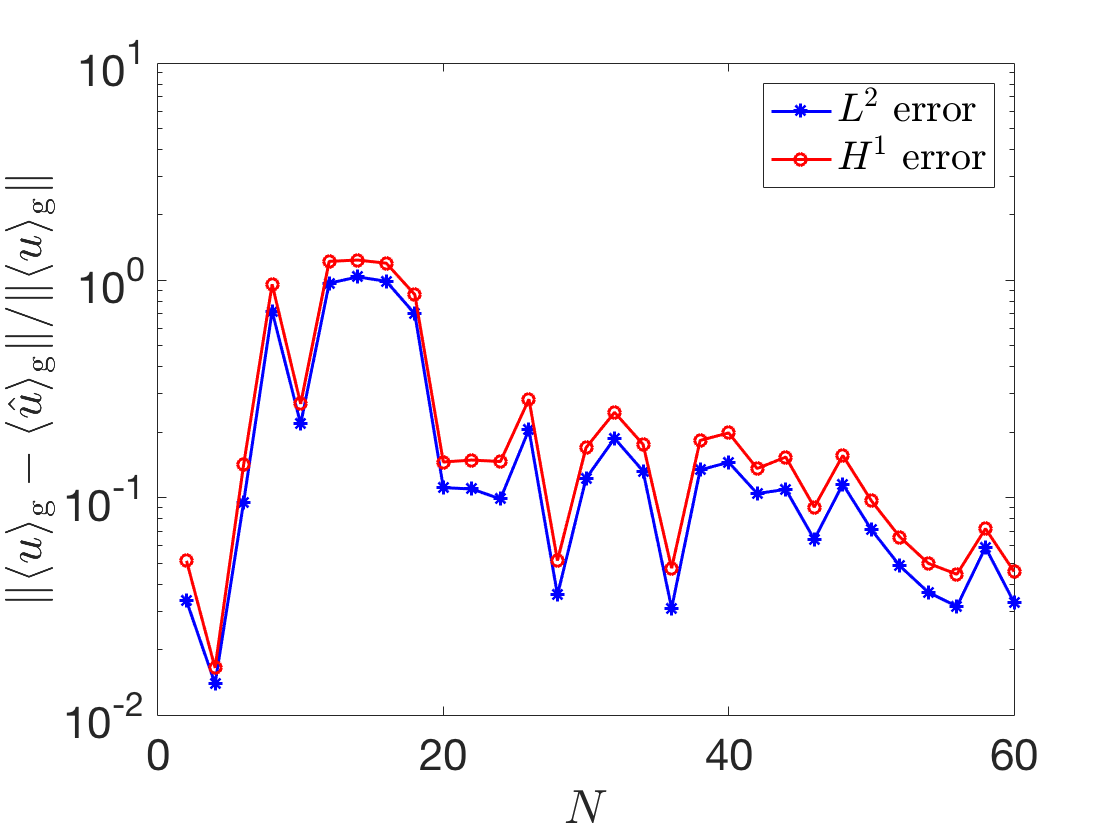}}
   ~~
 \subfloat[]
{  \includegraphics[width=0.30\textwidth]
 {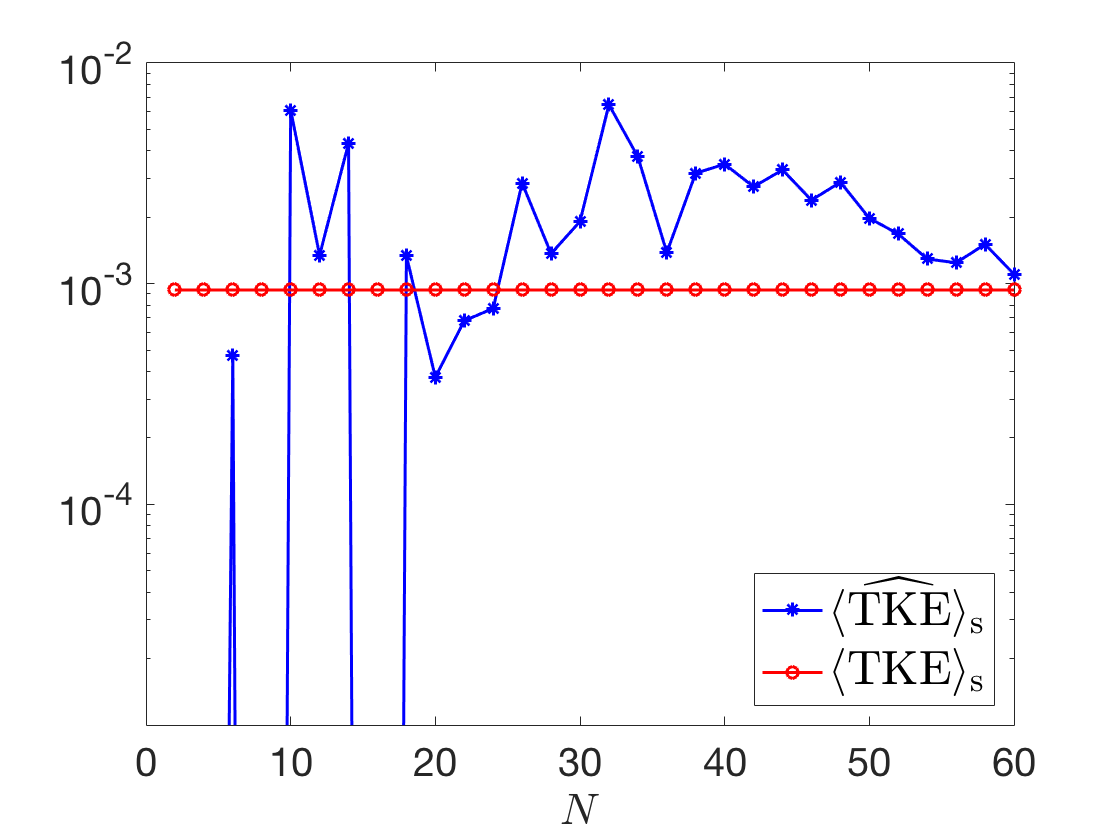}}
 
 \caption{The solution reproduction problem; POD-Galerkin.
 Figure (a): POD eigenvalues.
 Figure (b): behavior of the relative $L^2$ and $H^1$ errors in mean flow prediction with $N$.
  Figure (c): behavior of the mean TKE with $N$.
 (${\rm Re} = 15000$).
  }
 \label{fig:pod_gal_numerics}
  \end{figure}

Figure \ref{fig:pod_gal_numerics}(a) shows the behavior of the eigenvalues $\{ \lambda_n  \}_{n=1}^K$.
The first eigenmode is roughly proportional to $\langle u \rangle_{\rm g} - R_g$; 
provided that the estimate of the coefficients is accurate,
it does not contribute to the fluctuating field.  Therefore, we can identify the ratio
$$
r_N =
\frac{ \sum_{n=2}^N \lambda_n   }{  \sum_{n=2}^K \lambda_n  }
$$
as the portion of $H_0^1$ energy of the fluctuating field associated with the reduced POD space of dimension $N$.
We find that $r_N = 0.165$ for $N=2$, $r_N = 0.731$ for $N=20$, $r_N = 0.797$ for $N=30$, and 
$r_N = 0.87$ for $N=50$.
We observe that the decay with $N$ is rather slow; this suggests that accurate estimates of the entire system dynamics are out of reach for fully turbulent flows. 
Figure \ref{fig:pod_gal_numerics}(b) shows the behavior with $N$  of the relative error in the mean flow prediction:
$$
E_N^0 = \frac{ \| \langle u - \hat{u} \rangle_{\rm g}    \|_{L^2(\Omega)}   }{  \| \langle u   \rangle_{\rm g}    \|_{L^2(\Omega)}  },
\qquad
E_N^1 = \frac{ \| \langle u - \hat{u} \rangle_{\rm g}    \|_{H_0^1(\Omega)}   }{  \| \langle u   \rangle_{\rm g}    \|_{H_0^1(\Omega)}  };
$$
while
Figure \ref{fig:pod_gal_numerics}(c) shows the behavior with $N$ of the mean predicted TKE:  
$\langle \widehat{\rm TKE}  \rangle_{\rm s}$. We observe that for small values of $N$, we predict a false stable steady flow, while for moderate values of $N$ we substantially overestimate the TKE-
 Finally, for $N \gtrsim 50$ we observe a slow convergence of the Galerkin ROM to the mean values predicted by the high-fidelity model.

\begin{figure}[h!]
\centering
\subfloat[ $N=20$]  
{  \includegraphics[width=0.30\textwidth]
 {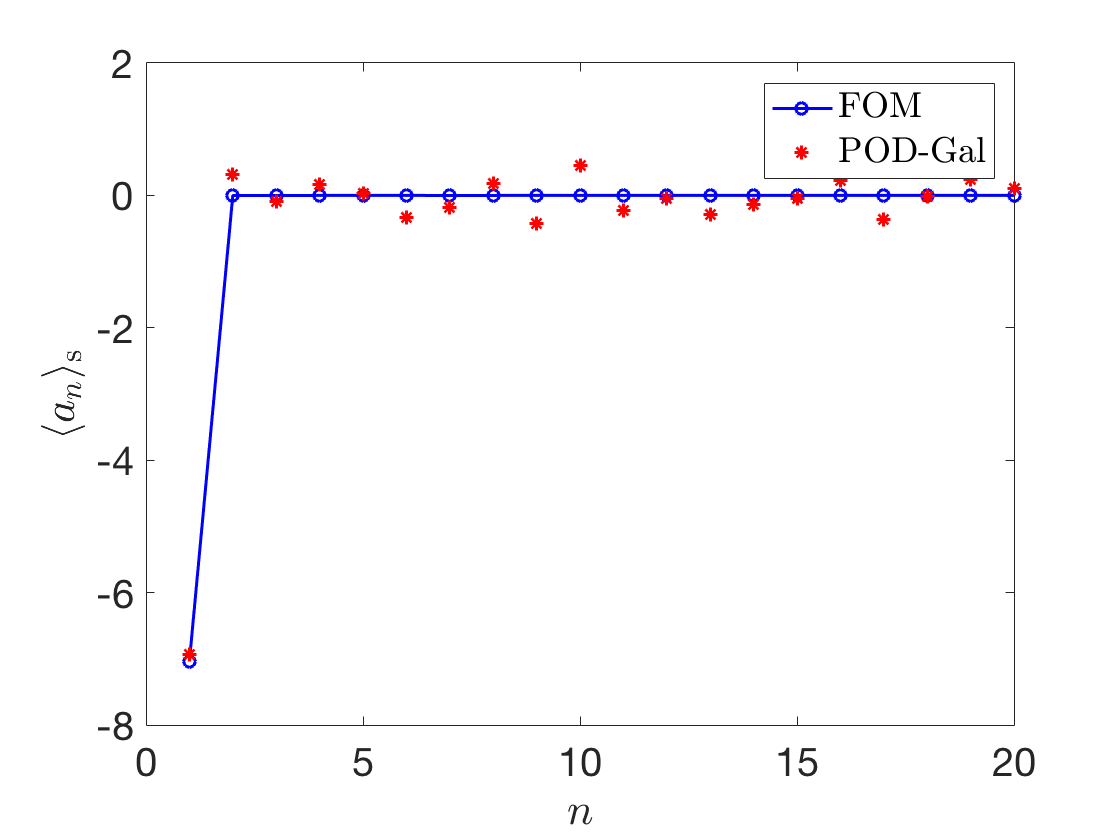}}
 ~~
  \subfloat[ $N=40$]  
{  \includegraphics[width=0.30\textwidth]
 {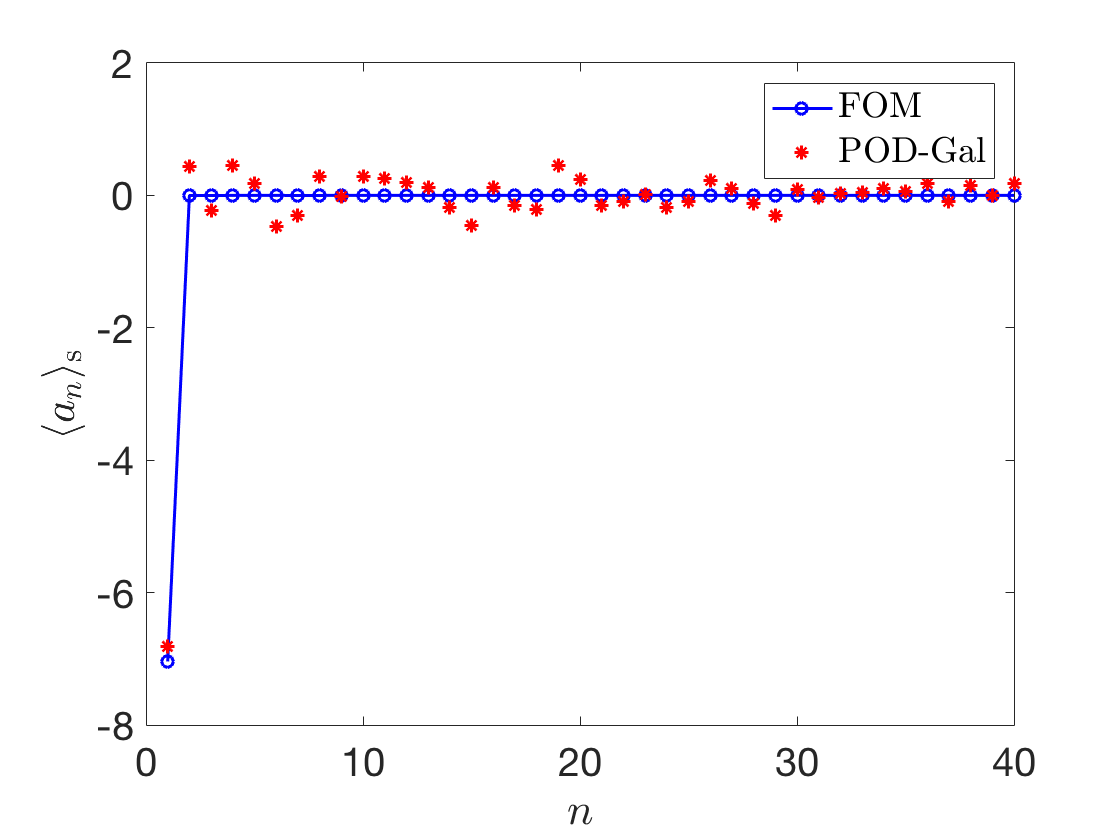}} 
  ~~
  \subfloat[ $N=60$]  
{  \includegraphics[width=0.30\textwidth]
 {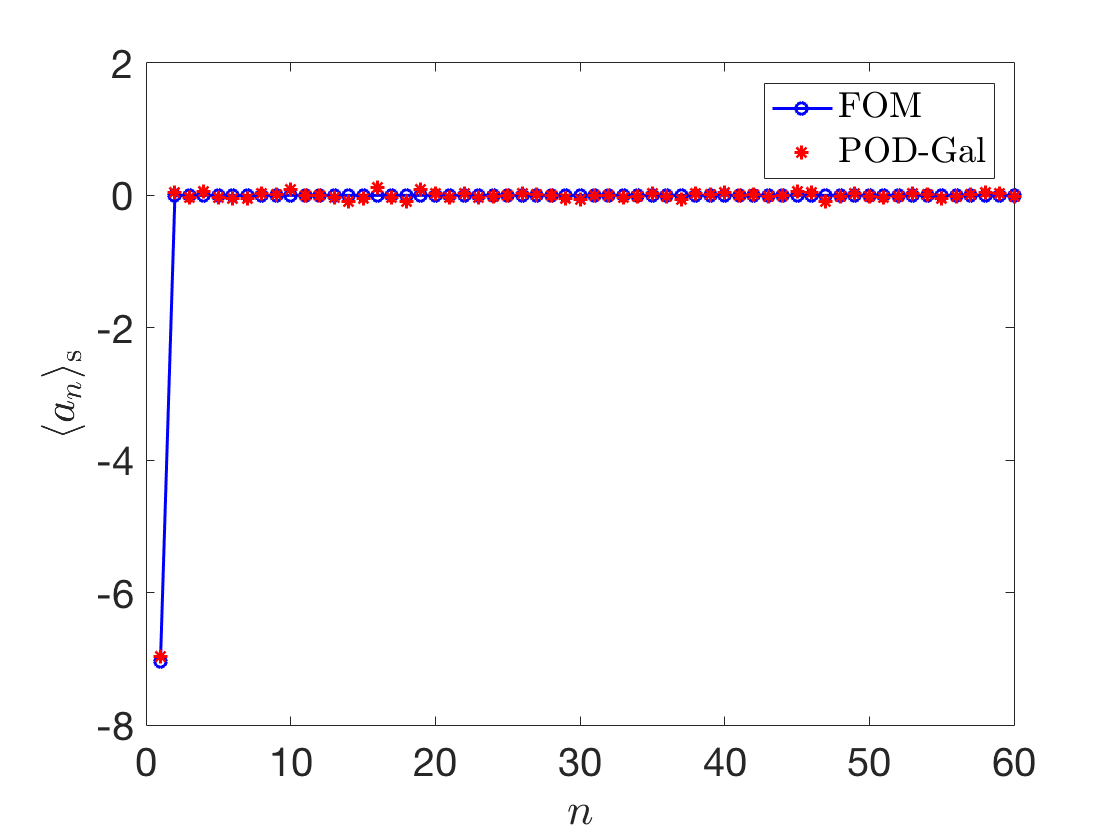}} 
 
 \subfloat[ $N=20$]  
{  \includegraphics[width=0.30\textwidth]
 {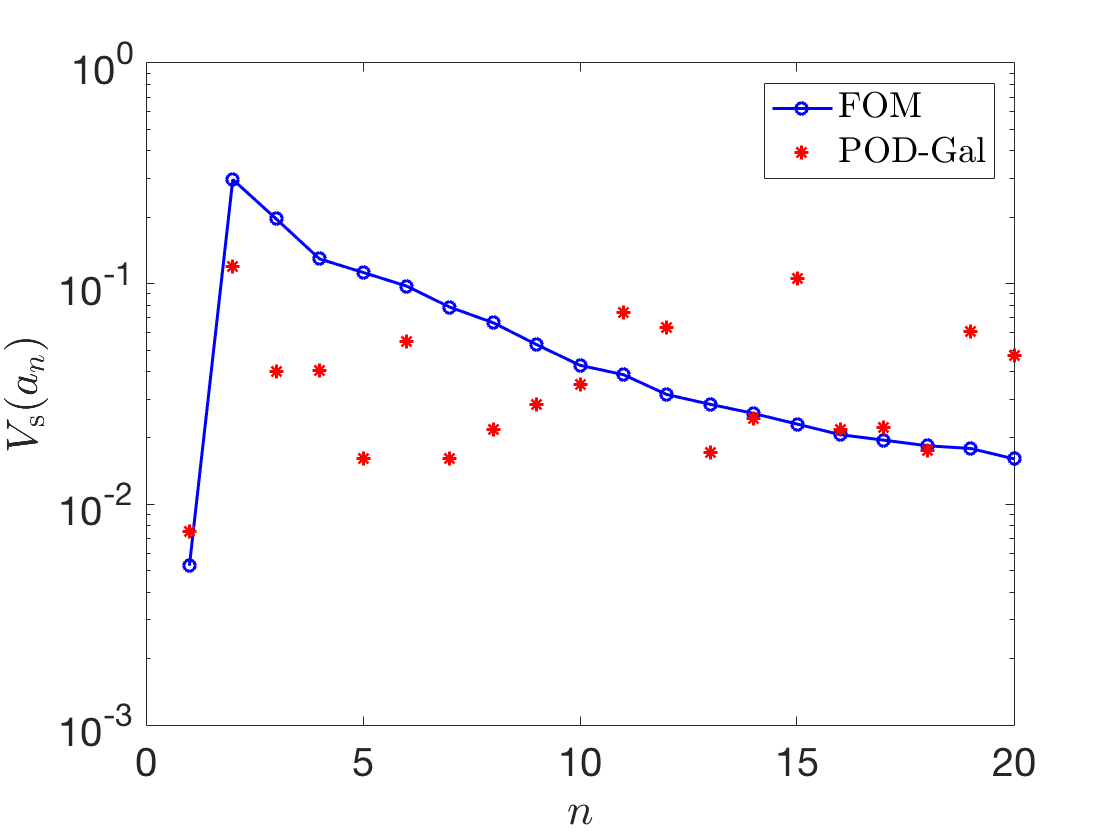}}
  ~~
\subfloat[ $N=40$]  
{  \includegraphics[width=0.30\textwidth]
 {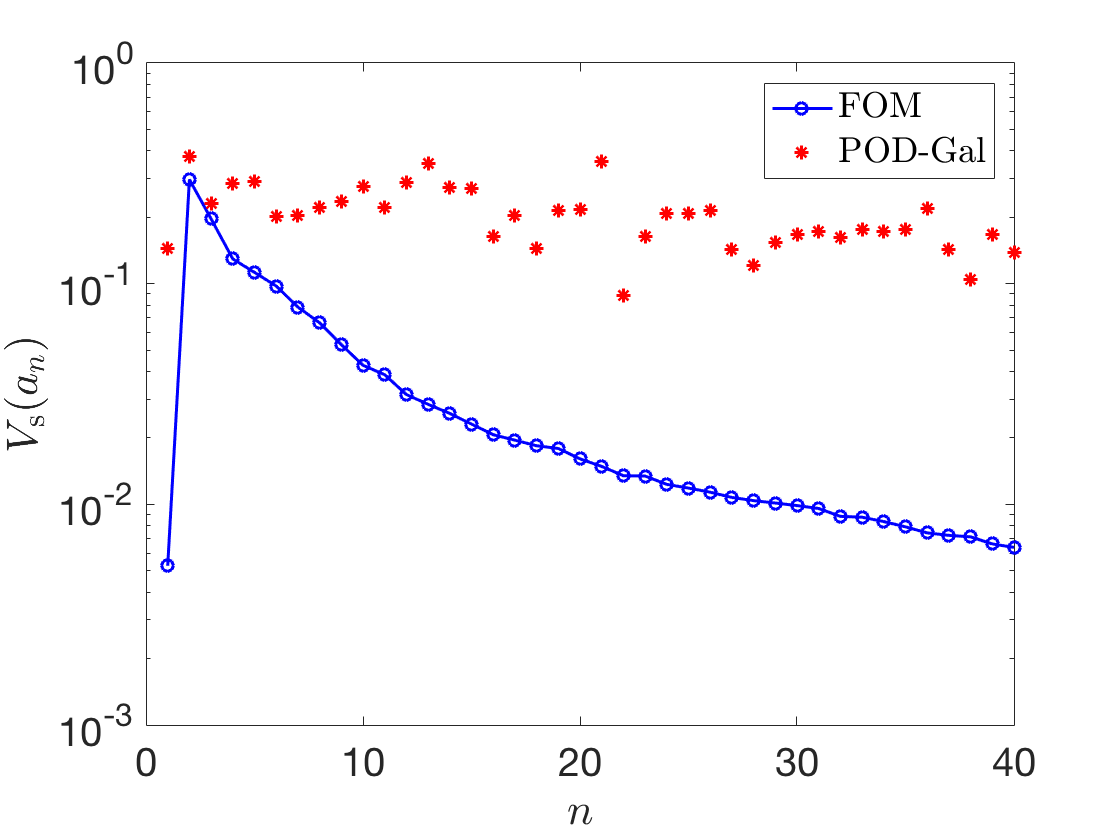}}
  ~~
\subfloat[ $N=60$]  
{  \includegraphics[width=0.30\textwidth]
 {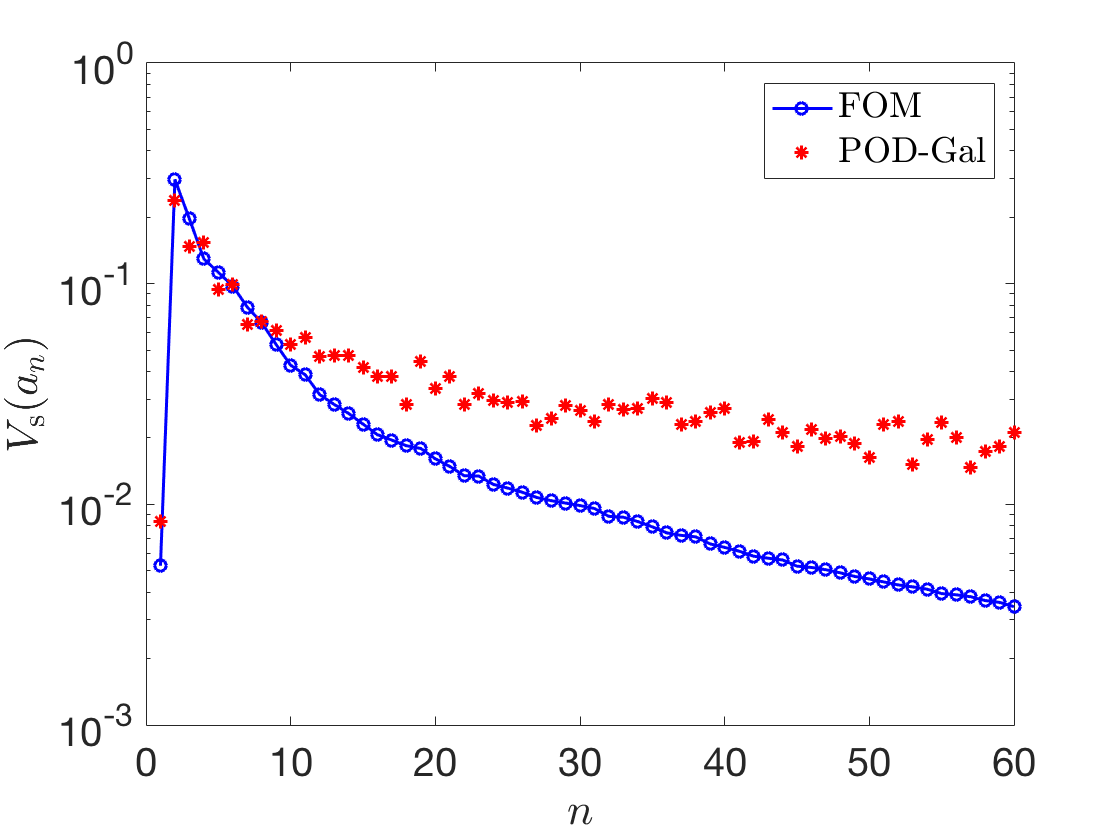}} 
 
 \caption{The solution reproduction problem; POD-Galerkin.
 Behavior of the sample mean and sample variance of the coefficients 
 $\{  a_n^j \}_j$.
 (${\rm Re} = 15000$).
  }
 \label{fig:pod_gal_numerics_2}
  \end{figure}

\begin{figure}[h!]
\centering
\subfloat[$N=20$]  
{  \includegraphics[width=0.30\textwidth]
 {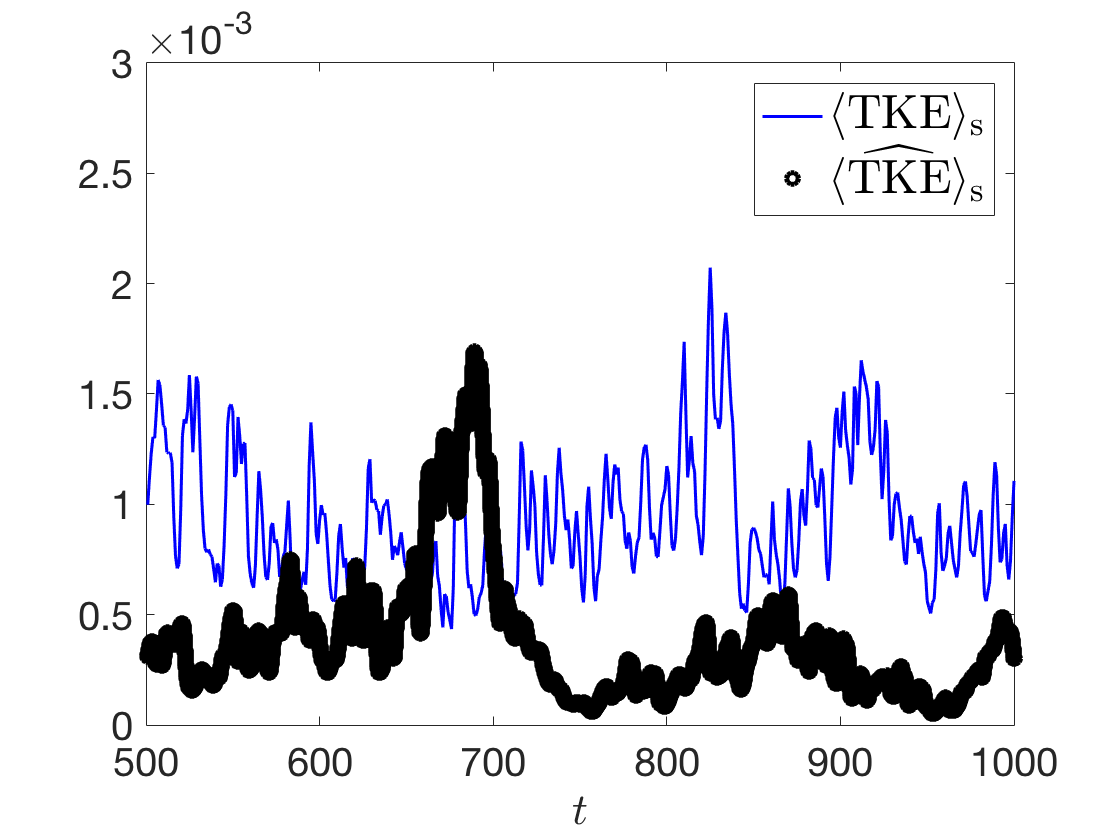}}
 ~~
\subfloat[$N=40$]  
{  \includegraphics[width=0.30\textwidth]
 {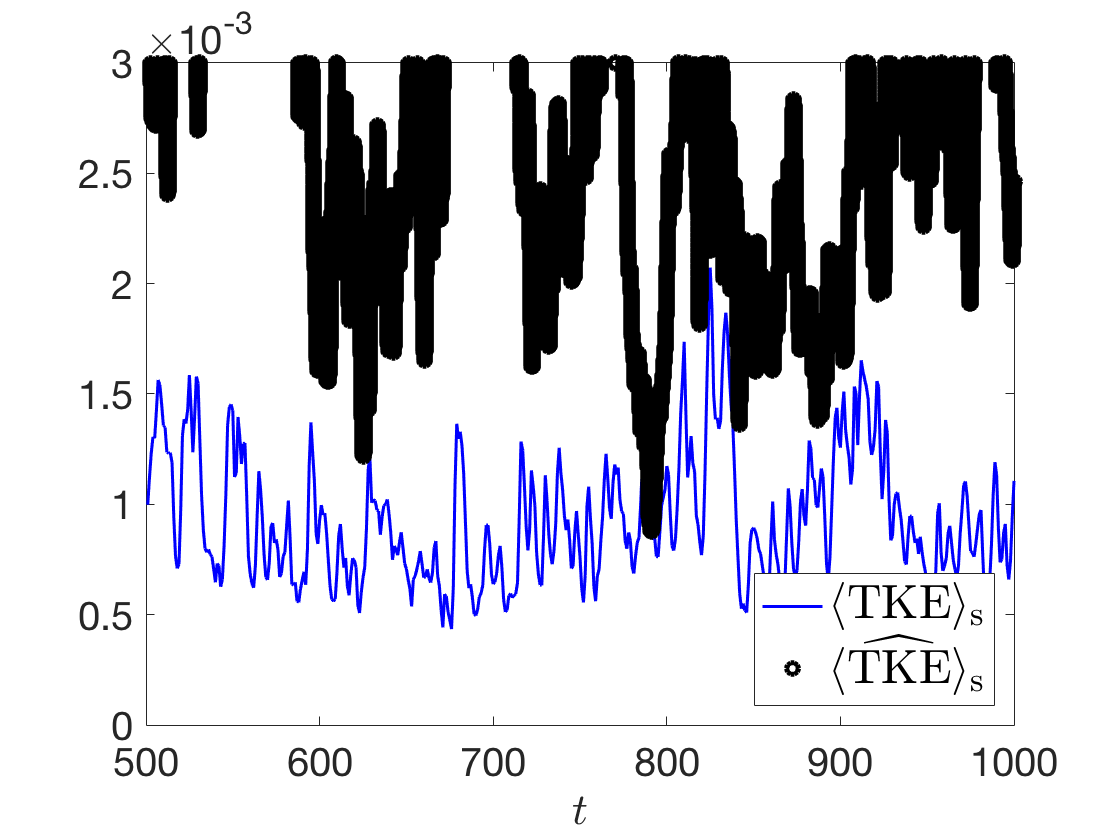}}
 ~~
\subfloat[$N=60$]  
{  \includegraphics[width=0.30\textwidth]
 {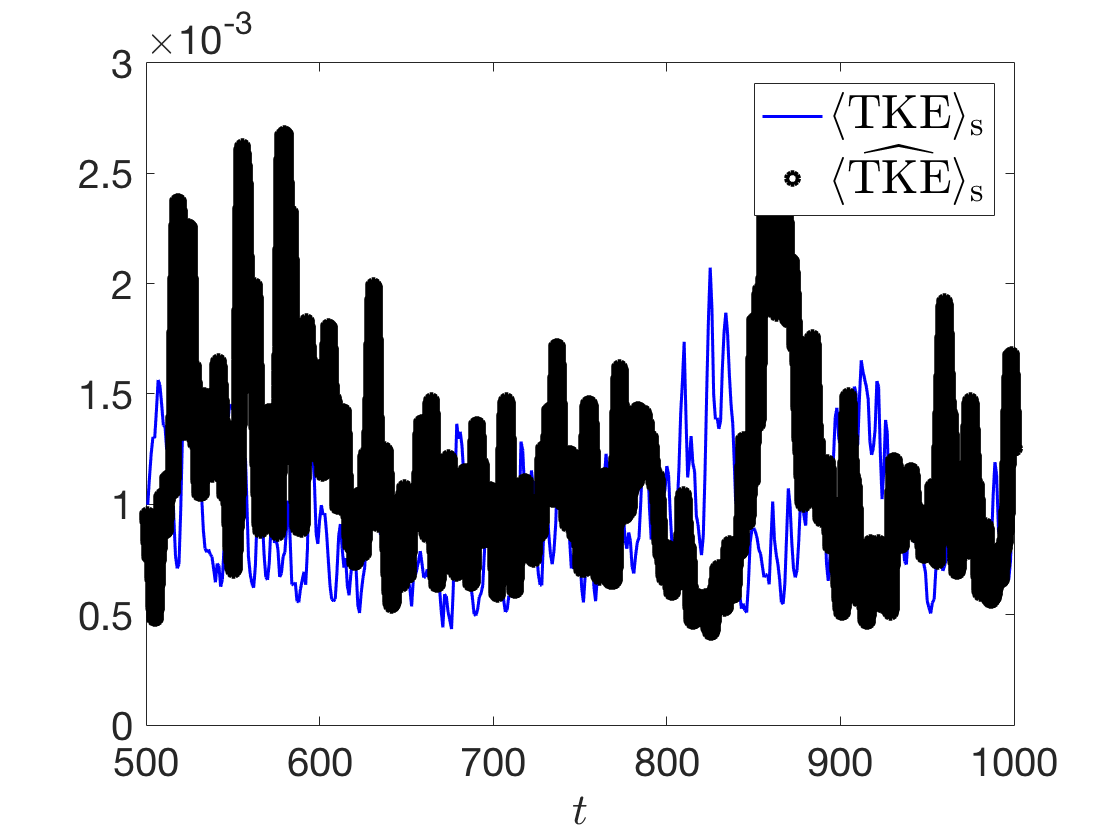}} 
 
 \caption{The solution reproduction problem; 
POD-Galerkin.
Behavior of the TKE as a function of time for three values of $N$.
$\langle \widehat{ {\rm TKE}} \rangle_{\rm s} = 3.8 \cdot 10^{-4}$ ($N=20$),
$3.5 \cdot 10^{-3}$ ($N=40$),
$1.1 \cdot 10^{-3}$ ($N=60$).
$V_{\rm s} (\widehat{ {\rm TKE}} )  = 8.8 \cdot 10^{-8}$ ($N=20$),
$6.5 \cdot 10^{-6}$ ($N=40$),
$1.9 \cdot 10^{-7}$ ($N=60$).
($\langle  {\rm TKE} \rangle_{\rm s} = 9.4 \cdot 10^{-4}$, 
$V_{\rm s} (  {\rm TKE} ) = 8.5 \cdot 10^{-8}$)
 (${\rm Re} = 15000$).
  }
 \label{fig:POD_Gal_energy_time}
  \end{figure}

Figure \ref{fig:pod_gal_numerics_2} shows the behavior for different values of $N$ of the sample mean and sample variance of the coefficients $\{ a_n^j \}_{j}$:
$$
\langle a_n \rangle_{\rm s}
=
\frac{1}{K}
\sum_{k=1}^K
\,
a_n(t_{\rm s}^k),
\qquad
V_{\rm s} (a_n)
=
\frac{1}{K-1}
\sum_{k=1}^K
\,\left(
a_n(t_{\rm s}^k) - 
\langle a_n \rangle_{\rm s}
\right)^2,
$$
for  the Full Order Model (FOM) and for the POD Galerkin ROM (POD-Gal).
Figure \ref{fig:POD_Gal_energy_time} shows the behavior of the TKE as a function of time for three values of $N$;
predictions of first and second order moments --- based on sampling times --- are reported in the caption of the Figure.
Results are consistent with the results in Figure \ref{fig:pod_gal_numerics}. For small-to-moderate values of $N$, we observe several spurious behaviors,
namely convergence to false stable steady flows, and overly unstable flows.
As $N$ increases, the accuracy of the Galerkin ROM appears to increase.

Interestingly, the behavior of the ROM observed here is qualitatively similar to the one observed for highly-truncated spectral approximations to turbulent flows (\cite{curry1984order}). We argue that the need for large reduced spaces might greatly reduce the benefit of Model Reduction:
if the value of  $N$ required to obtain sufficiently accurate results 
is too large, the resulting ROM might not lead to significant computational speed-ups, and might also not be beneficial in terms of memory.
This observation motivates the correction to the Galerkin formulation proposed in the next section.
We finally remark that the  results shown in this section suggest the need for a pragmatic definition of long-time stability: we address this issue in Appendix \ref{sec:ROM_stability}.

\subsection{The constrained POD-Galerkin formulation}
\label{sec:cgalerkin_ROM}

\subsubsection{Formulation}
\label{sec:cgalerkin_formulation}

Given the reduced space $\mathcal{Z}^{\rm u} = {\rm span} \{ \zeta_n \}_{n=1}^N \subset V_{\rm div}$, and the time grid $\{ t_{\rm g}^j \}_{j=0}^J$, we seek the coefficients $\{ \mathbf{a}^j \}_{j=0}^J \subset \mathbb{R}^N$ such that
\begin{equation}
\label{eq:constrained_galerkin}
\mathbf{a}^{j+1}:=
{\rm arg} \min_{\mathbf{a} \in \mathbb{R}^N}
\,
\|
\mathbb{A}(\mathbf{a}^j; {\rm Re})
\mathbf{a} 
-
\mathbf{F}(\mathbf{a}^j; {\rm Re})
\|_2^2,
\quad
{\rm s.t.}
\; \;
\alpha_n \leq a_n \leq \beta_n,
\; \;
n=1,\ldots,N;
\end{equation}
where $\{ \alpha_n \}_{n=1}^N$ and $\{ \beta_n \}_{n=1}^N$ are suitable hyper-parameters that will be specified later, and 
$\mathbb{A}(\cdot;  {\rm Re})$, $\mathbf{F}(\cdot;  {\rm Re})$ are defined in \eqref{eq:galerkin_algebraic}.
Formulation \eqref{eq:constrained_galerkin} reads as a constrained quadratic programming problem where the objective function corresponds to the Euclidean norm of the error in the reduced Galerkin formulation, while the constraints impose that each coefficient of the $N$-term expansion remains in the interval 
$[\alpha_n, \beta_n]$, $n=1,\ldots,N$. 
We refer to \eqref{eq:constrained_galerkin} as \emph{constrained (POD-)Galerkin formulation}.

The hyper-parameters $\{ \alpha_n \}_n$ and $\{ \beta_n \}_n$ are designed to embed in the ROM formulation information about the variation in time of the coefficients $\{ a_n^j \}_j$, for $n=1,\ldots,N$.
For each value of $n$, if we introduce the 
projection\footnote{We assume here that the POD eigenmodes $\{ \zeta_n \}_{n=1}^N$ are orthonormalized. }
 of the lifted field on the $n$-th POD mode at time $t_{\rm g}^j$,
$a_n^{\rm FOM,j}:= (\mathring{u}^j, \zeta_n)_V$,
we can interpret $\alpha_n$ and $\beta_n$  as lower and upper bounds  for the sequence  $\{ a_n^{\rm FOM,j} \}_{j=J_0}^J$,
where $J_0>0$ is introduced in \eqref{eq:mean_flow_def} to discard the transient dynamics.
The hyper-parameters $\alpha_n$ and $\beta_n$ are not directly related to the POD eigenvalues $\lambda_n$: the latter are  --- up to a multiplicative constant --- 
estimates of the squared $\ell^2$-norm of the coefficients,
$\lambda_n = \sum_k \left( a_n^{\rm FOM,k}  \right)^2  
\approx \frac{K}{J} \sum_j  \left( a_n^{\rm FOM, j}  \right)^2$.

Based on the interpretation of  the hyper-parameters, 
 we propose to estimate  $\{ \alpha_n \}_{n=1}^N$ and $\{ \beta_n \}_{n=1}^N$ based on the sample minima and the sample maxima associated with  the snapshots $\{ \mathring{u}^k \}_{k=1}^K$: 
\begin{subequations} 
\label{eq:choice_coefficients}
\begin{equation}
\alpha_n:= m_n^{\rm u} - \epsilon \Delta_n^{\rm u},
\qquad
\beta_n:= M_n^{\rm u} + \epsilon \Delta_n^{\rm u},
\end{equation}
where  $m_n^{\rm u}$ and $M_n^{\rm u}$ are sample minimum and sample maximum associated with the projection of the lifted field on the $n$-th POD mode,
\begin{equation}
\label{eq:sample_min_max}
m_n^{\rm u}
=
\min_k \, a_n^{\rm FOM, k} :=
(\mathring{u}^k, \zeta_n)_V,
\quad
M_n^{\rm u}
=
\max_k \, a_n^{\rm FOM, k};
\end{equation}
$\Delta_n^{\rm u}$ is the sample estimate of the difference between maximum and minimum,
\begin{equation}
\Delta_n^{\rm u}
:=
M_n^{\rm u} - m_n^{\rm u};
\end{equation}
and the constant $\epsilon>0$ takes into account the fact that sample minima and sample maxima 
in \eqref{eq:sample_min_max} are upper and lower bounds for the true minima and true maxima, respectively.
We emphasize that in our framework $K \ll J$; therefore, $\{ a_n^{\rm FOM, k} \}_{k=1}^K$ should be interpreted as a (deterministic) sample from the population $\{  a_n^{\rm FOM, j} \}_{j=J_0}^J$.
Given the special features of the learning task at hand ---
the estimation of minima and maxima of a population ---
we expect that we can estimate the hyper-parameters based on \emph{sparse} DNS data (i.e., data that are not dense in any specific region of the time interval). 
\end{subequations}

Accurate estimates of the hyper-parameters of the formulation based on sparse DNS data represent an important feature of our constrained formulation. 
As observed by many authors, low-frequency features of the turbulent flow
---  which largely contribute to long-time flow averages ---
are well-represented by the snapshots $\{ \mathring{u}^k \}_{k=1}^K$
 and consequently by the POD space only if the sampling times $\{ t_{\rm s}^k \}_{k=1}^K$ are not clustered in any specific region of the time interval.
 This implies that both the ingredients  of the ROM 
 --- the space $\mathcal{Z}^{\rm u}$ and the hyper-parameters $\{ \alpha_n \}_n$ and 
$\{ \beta_n \}_n$ ---
require the same sampling strategy for the construction of the snapshot set. 
Therefore, the same dataset used to generate the POD space is well-suited to estimate the hyper-parameters of the ROM. This observation allows us to limit the size  $K$ of the snapshot set, and ultimately leads to a reduction of  the offline memory cost.

Unlike the standard POD-Galerkin ROM,  we here use DNS data twice: first, to build the space $\mathcal{Z}^{\rm u}$; second,  to estimate the hyper-parameters 
$\{ \alpha_n \}_{n=1}^N$ and $\{ \beta_n \}_{n=1}^N$.
Furthermore, while  POD-Galerkin is independent of the particular basis 
$\zeta_1,\ldots,\zeta_N$ chosen for $\mathcal{Z}^{\rm u}$, the box constraints in \eqref{eq:galerkin_algebraic} depend on the choice of the basis.
We emphasize that by choosing $\{ \zeta_n\}_n$ as basis of  $\mathcal{Z}^{\rm u}$ we   explicitly incorporate (prior) information about the decay of the POD coefficients directly in the formulation.

We observe that if the solution to Galerkin ROM
\eqref{eq:galerkin_algebraic} ---
$\mathbf{a}_{\rm Gal}^{j+1} = \mathbb{A} (\mathbf{a}^j; {\rm Re})^{-1}$ $\mathbf{F}(\mathbf{a}^j; {\rm Re})$ ---
satisfies the box constraints in \eqref{eq:constrained_galerkin}, then $\mathbf{a}^{j+1}= \mathbf{a}_{\rm Gal}^{j+1}$. Therefore, our constrained formulation \emph{corrects} the Galerkin formulation only if 
$\mathbf{a}_{\rm Gal}^{j+1}$ does not satisfy the prescribed bounds. This represents the main difference between our approach and the other stabilized ROMs proposed in the literature and briefly mentioned in the introduction:
rather than introducing artificial dissipation in the Galerkin model, we exploit prior information about the attractor to correct the ROM. 

We finally comment on time discretization.
In this work, we employ the first-order semi-implicit time-discretization introduced in \eqref{eq:galerkin_rom_discretized}. However, the approach can be trivially extended to other time discretizations: first, we derive the discrete Galerkin ROM from \eqref{eq:galerkin_rom_non_discretized},
then we substitute the resulting algebraic formulation in the objective function of \eqref{eq:constrained_galerkin}. For explicit and semi-implicit single-step time integrators, the resulting constrained formulation corresponds to a quadratic programming problem, which can be solved using interior-point methods  (see, e.g., \cite{nocedal2006numerical}).
For fully-implicit single-step methods, the constrained formulation reads as a nonlinear constrained optimization problem, which again can be solved using interior-point methods or sequential quadratic programming.
We envision that the extension to multistep methods  might require some additional care since the solution is not guaranteed to be smooth in time when the constraints are active. A thorough analysis of different time integrators is beyond the scope of this paper.

\subsubsection{Performance of the constrained POD-Galerkin ROM}
\label{sec:cgalerkin_performance}
We present numerical results for ${\rm Re}=15000$.
Time grid $\{ t_{\rm g}^j \}_{j=0}^J$ and sampling times 
$\{ t_{\rm s}^k \}_{k=1}^K$ are the same considered for POD-Galerkin, if not specified otherwise. As for the previous test  the ROM is implemented in Matlab;  the quadratic programming problem is solved using the routine \texttt{quadprog} based on an interior-point algorithm. 
We here set $\epsilon=0.01$ in \eqref{eq:choice_coefficients}.

\begin{figure}[h!]
\centering
 \subfloat[]
{  \includegraphics[width=0.30\textwidth]
 {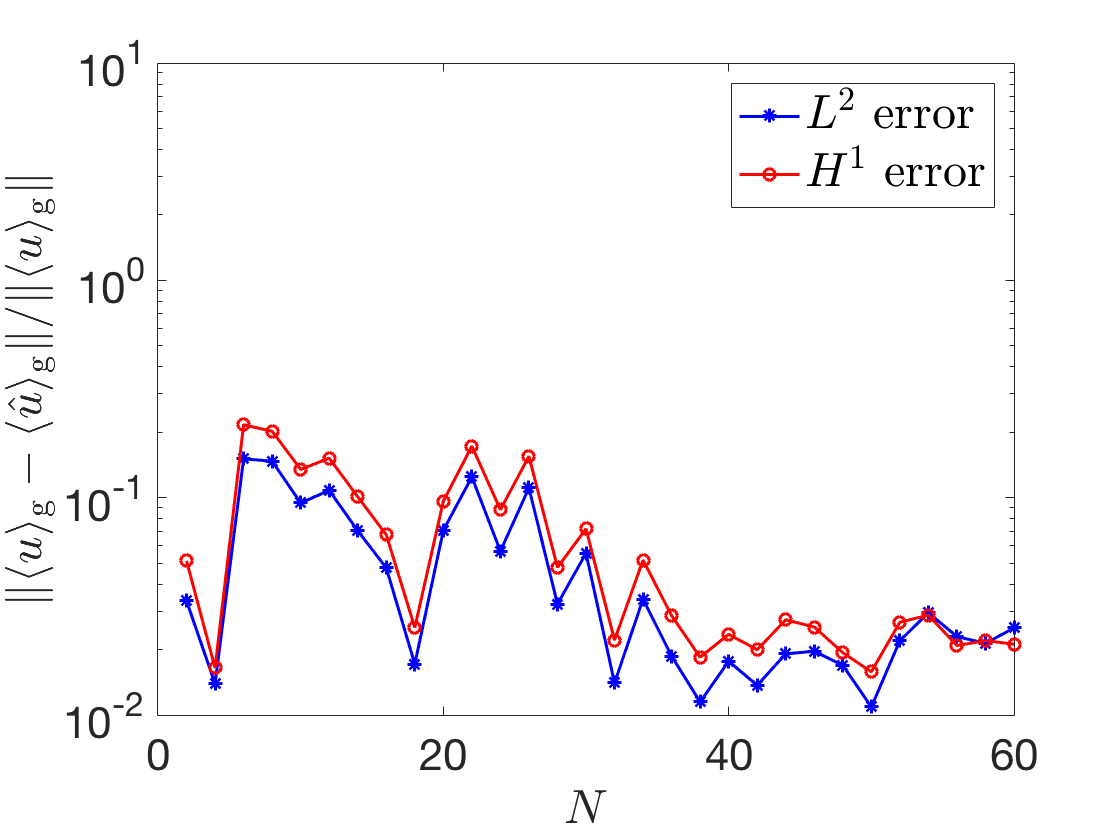}}
   ~~
 \subfloat[]
{  \includegraphics[width=0.30\textwidth]
 {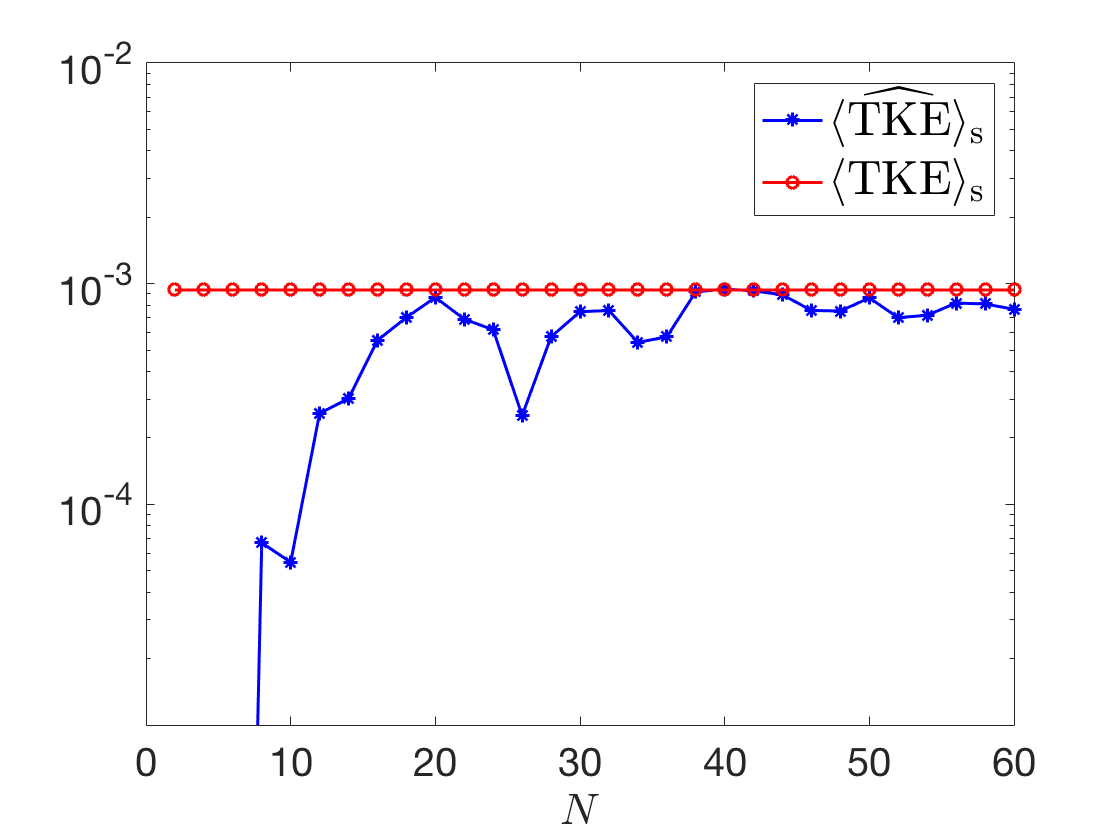}}
 
 \caption{The solution reproduction problem; constrained POD-Galerkin.
 Figure (a): behavior of the relative $L^2$ and $H^1$ errors in mean flow prediction with $N$.
  Figure (b): behavior of the mean TKE with $N$.
 (${\rm Re} = 15000$, $\epsilon=0.01$).
  }
 \label{fig:cgal_numerics}
  \end{figure}

\begin{figure}[h!]
\centering
 \subfloat[]
{  \includegraphics[width=0.30\textwidth]
 {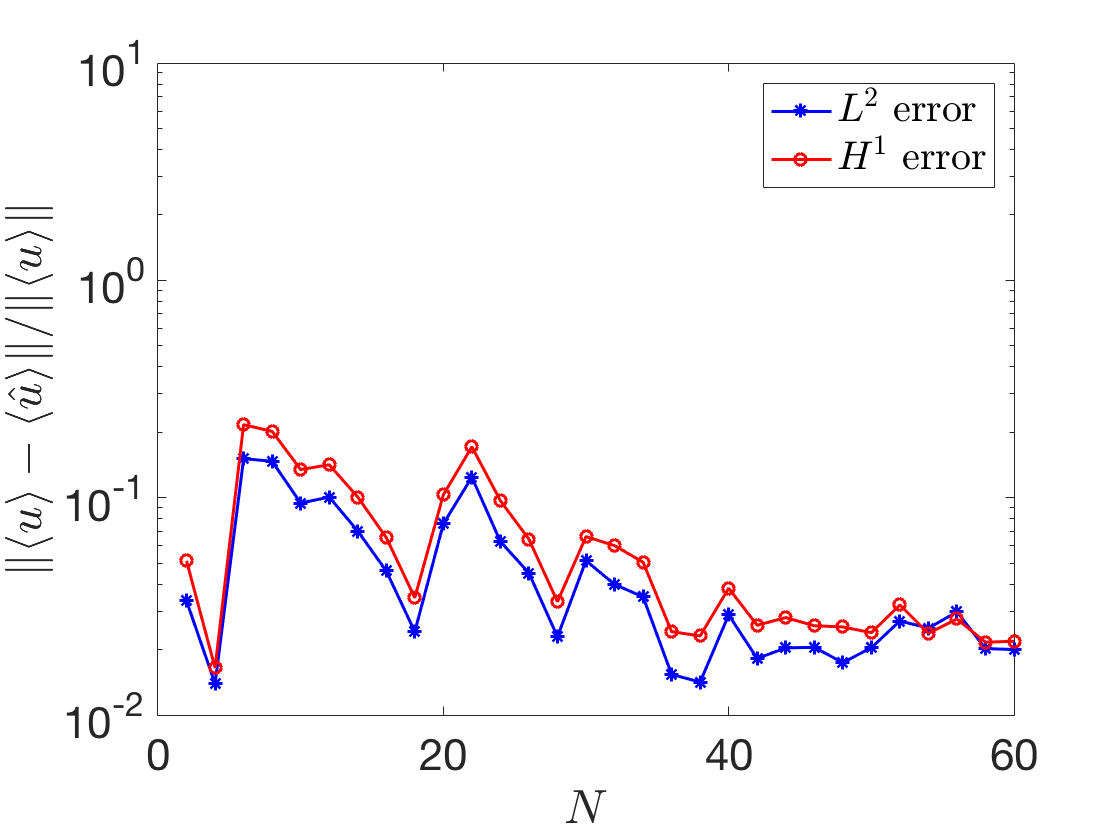}}
   ~~
 \subfloat[]
{  \includegraphics[width=0.30\textwidth]
 {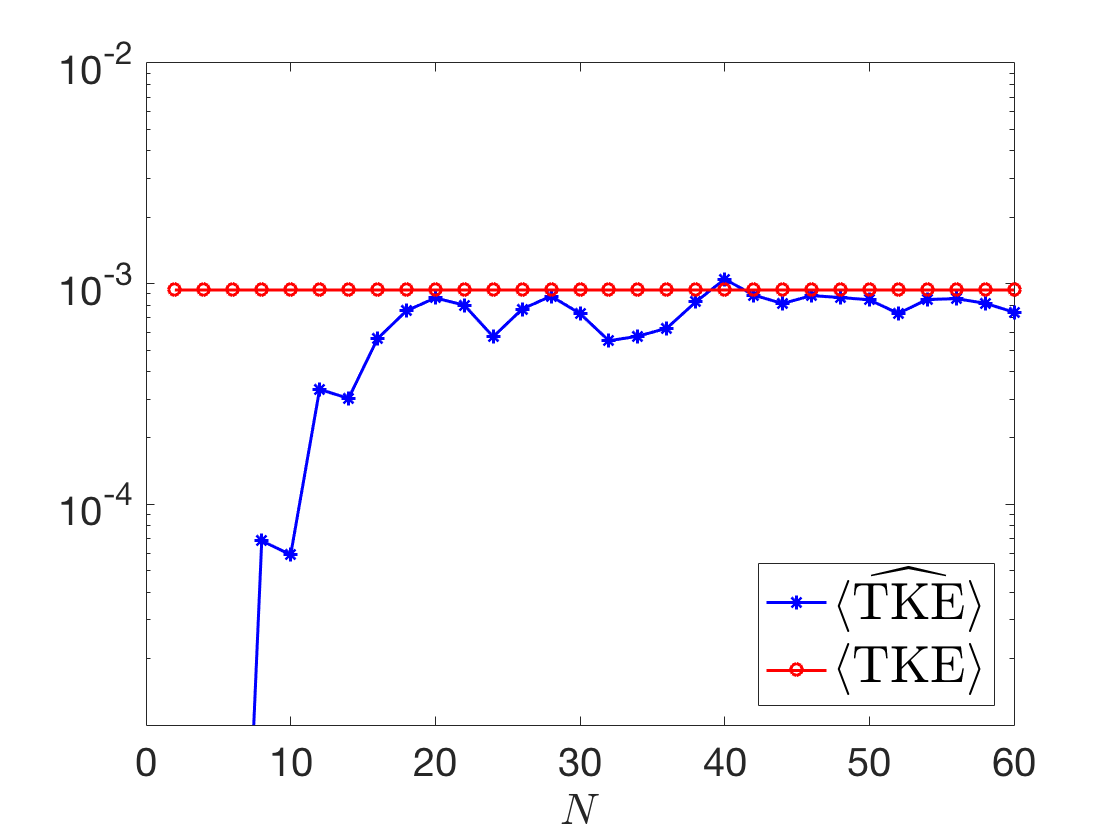}}
 
 \caption{The solution reproduction problem; constrained POD-Galerkin for a finer time grid.
 Figure (a): behavior of the relative $L^2$ and $H^1$ errors in mean flow prediction with $N$.
  Figure (b): behavior of the mean TKE with $N$.
 (${\rm Re} = 15000$, $\epsilon=0.01$, $\Delta t = 2.5 \cdot 10^{-3}$).
  }
 \label{fig:cgal_numerics_time}
  \end{figure}

\begin{figure}[h!]
\centering
\subfloat[ $N=20$]  
{  \includegraphics[width=0.30\textwidth]
 {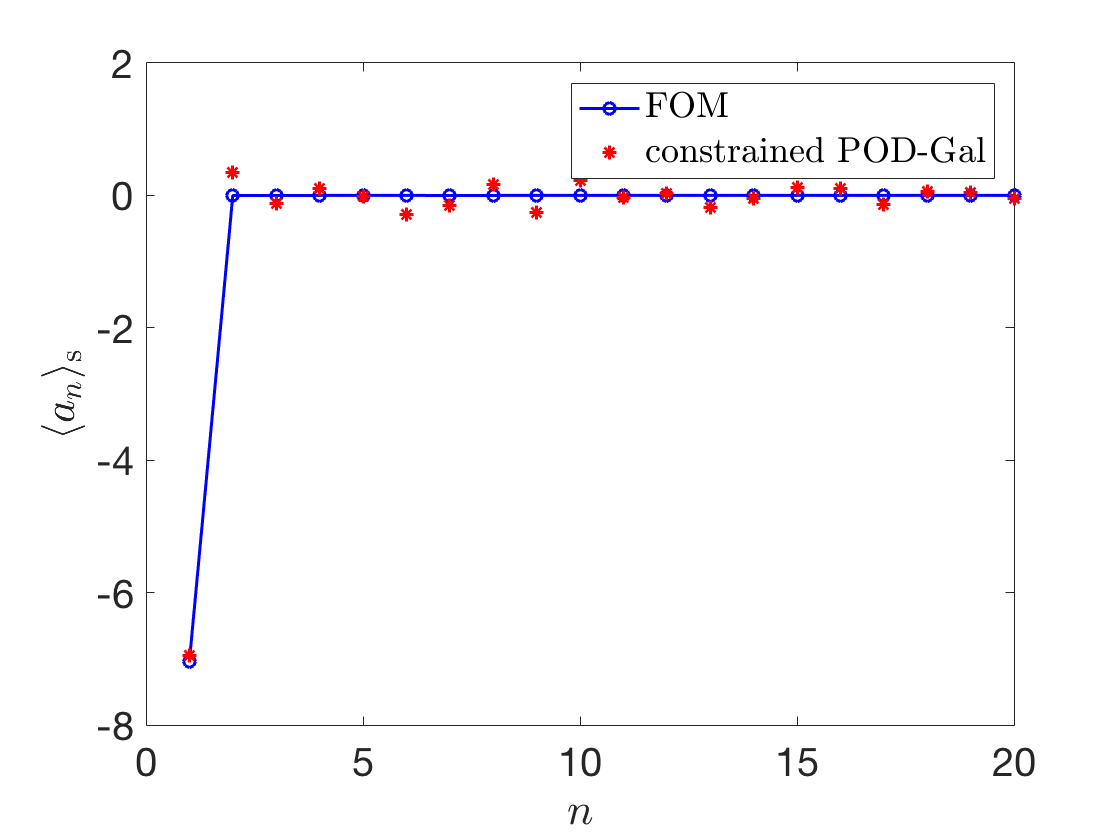}}
 ~~
  \subfloat[ $N=40$]  
{  \includegraphics[width=0.30\textwidth]
 {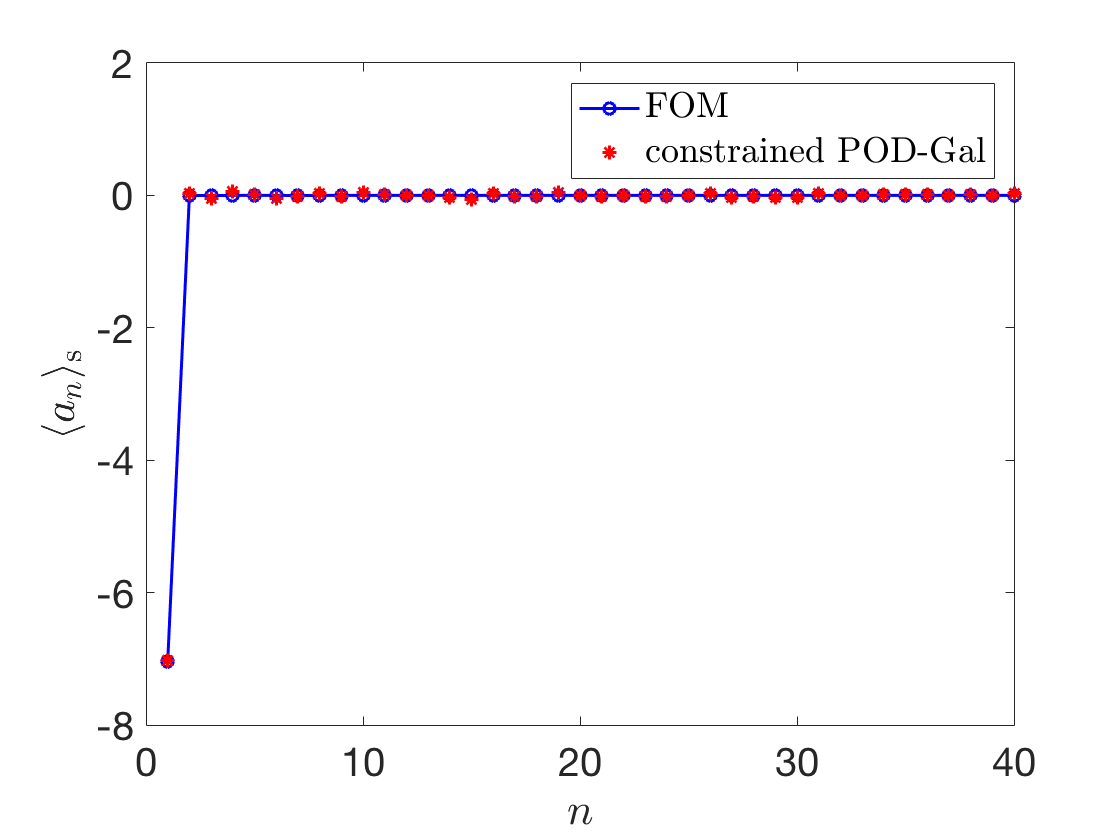}} 
  ~~
  \subfloat[ $N=60$]  
{  \includegraphics[width=0.30\textwidth]
 {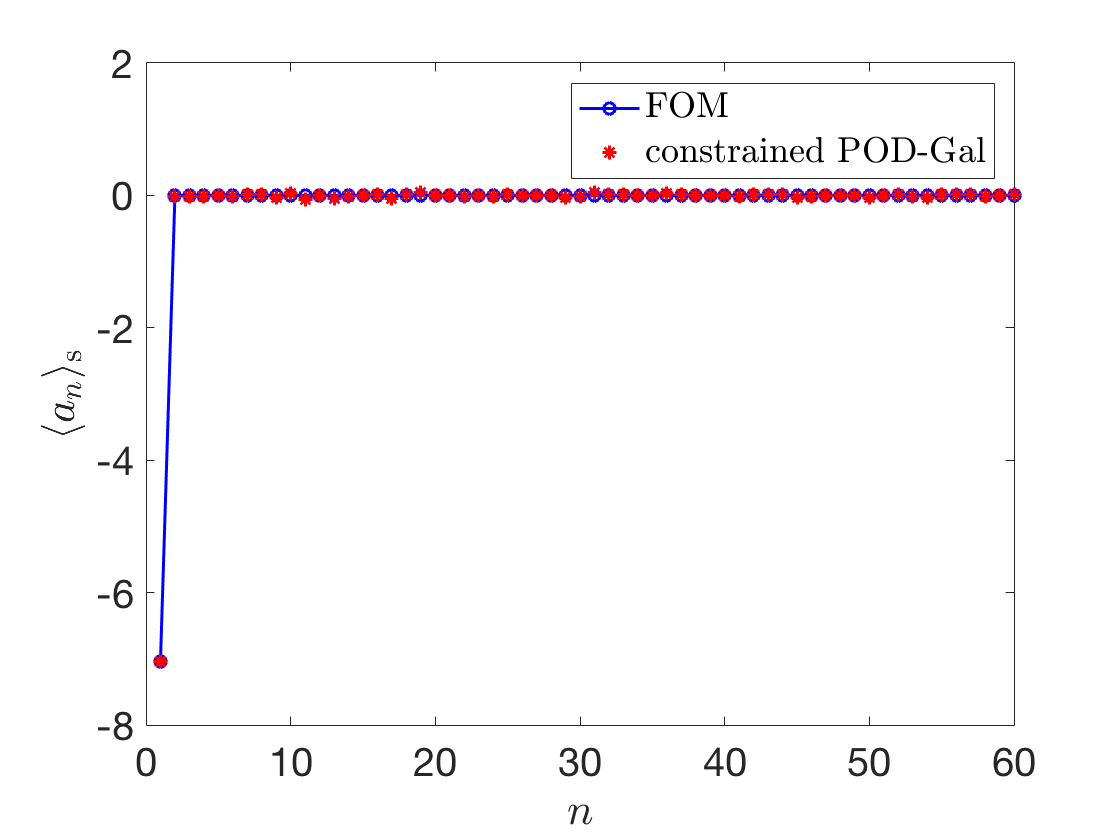}} 
 
 \subfloat[ $N=20$]  
{  \includegraphics[width=0.30\textwidth]
 {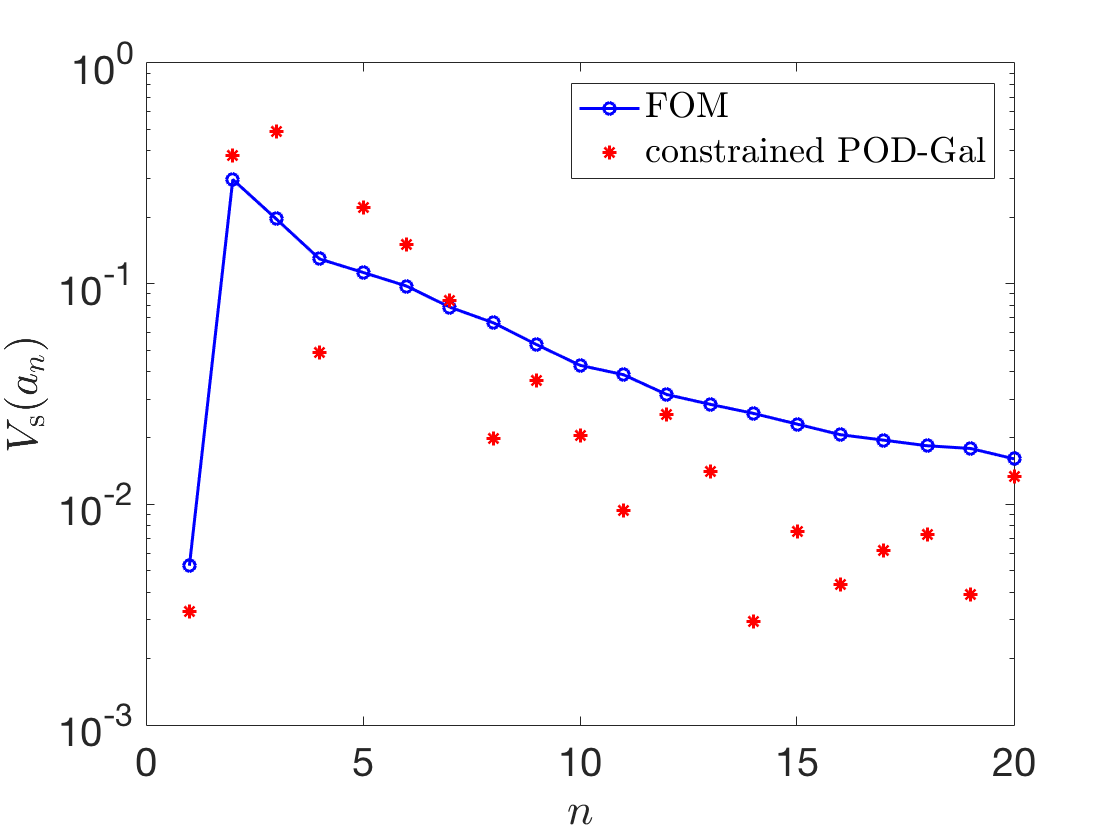}}
  ~~
\subfloat[ $N=40$]  
{  \includegraphics[width=0.30\textwidth]
 {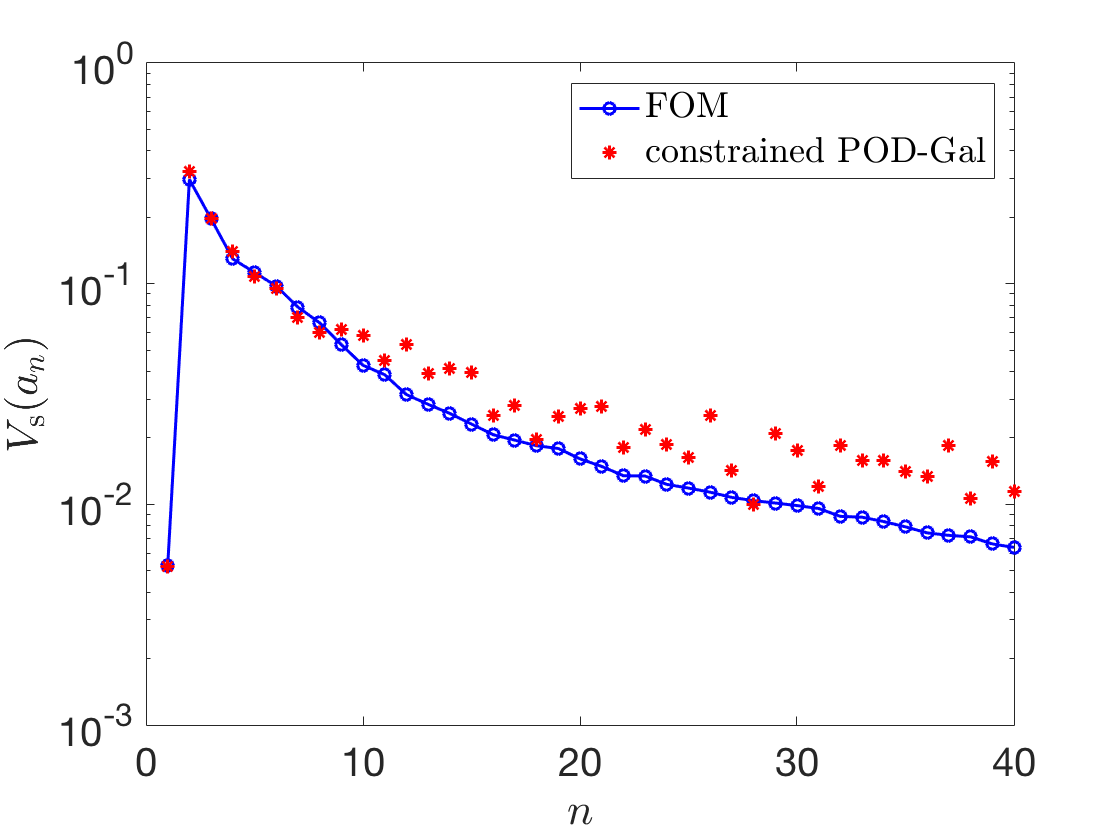}}
  ~~
\subfloat[ $N=60$]  
{  \includegraphics[width=0.30\textwidth]
 {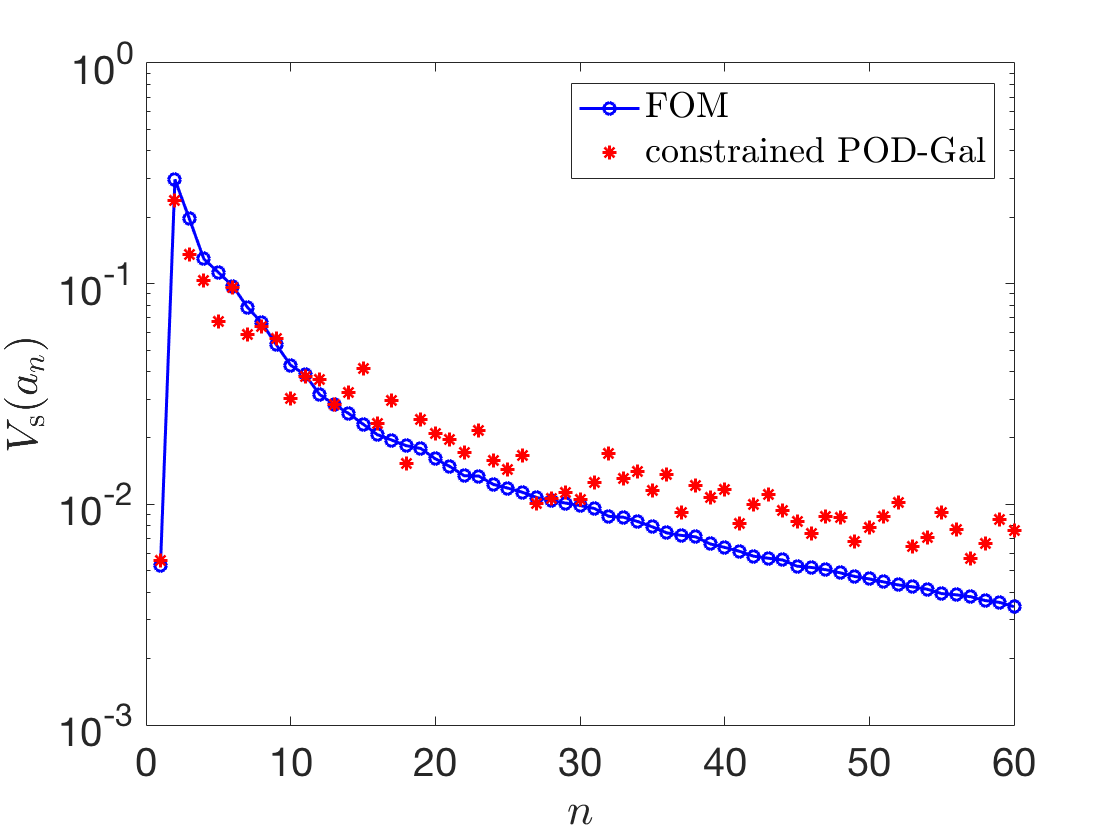}} 
 
 \caption{The solution reproduction problem; constrained POD-Galerkin.
 Behavior of the sample mean and sample variance of the coefficients 
 $\{  a_n^j \}_j$.
 (${\rm Re} = 15000$, $\epsilon=0.01$).
  }
 \label{fig:cgal_numerics_2}
  \end{figure}

\begin{figure}[h!]
\centering
\subfloat[$N=20$]  
{  \includegraphics[width=0.30\textwidth]
 {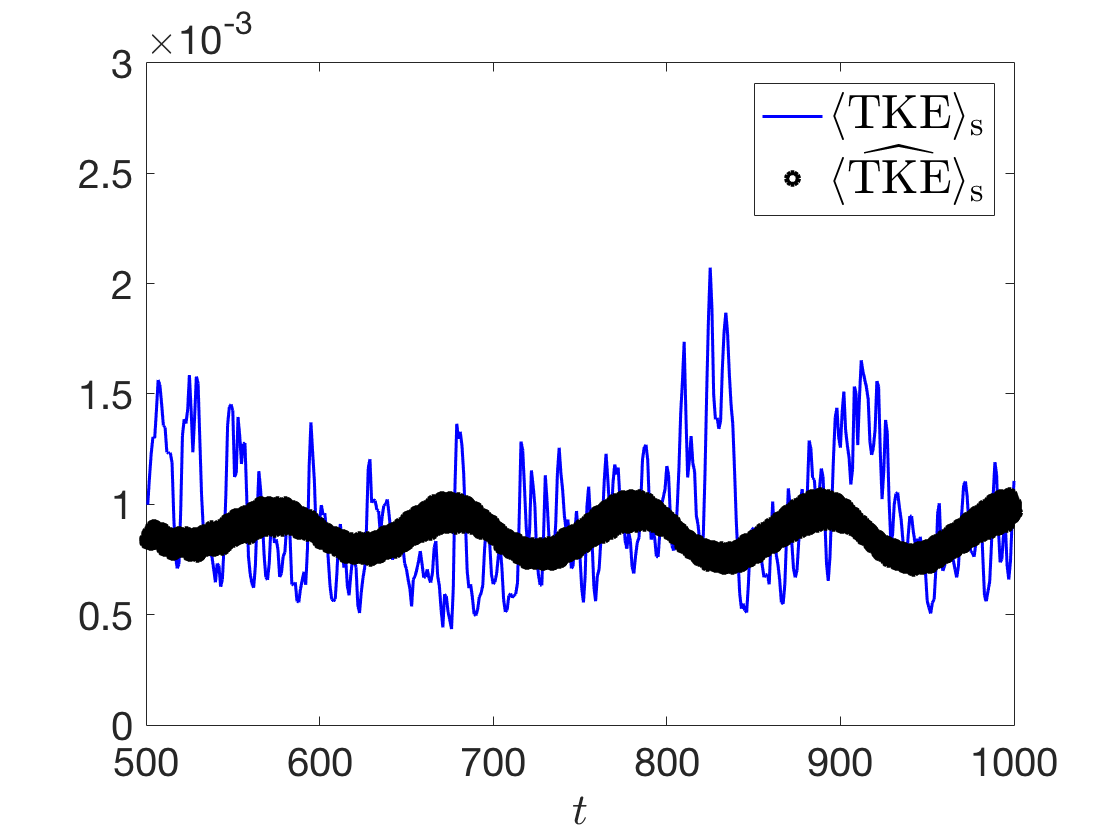}}
 ~~
\subfloat[$N=40$]  
{  \includegraphics[width=0.30\textwidth]
 {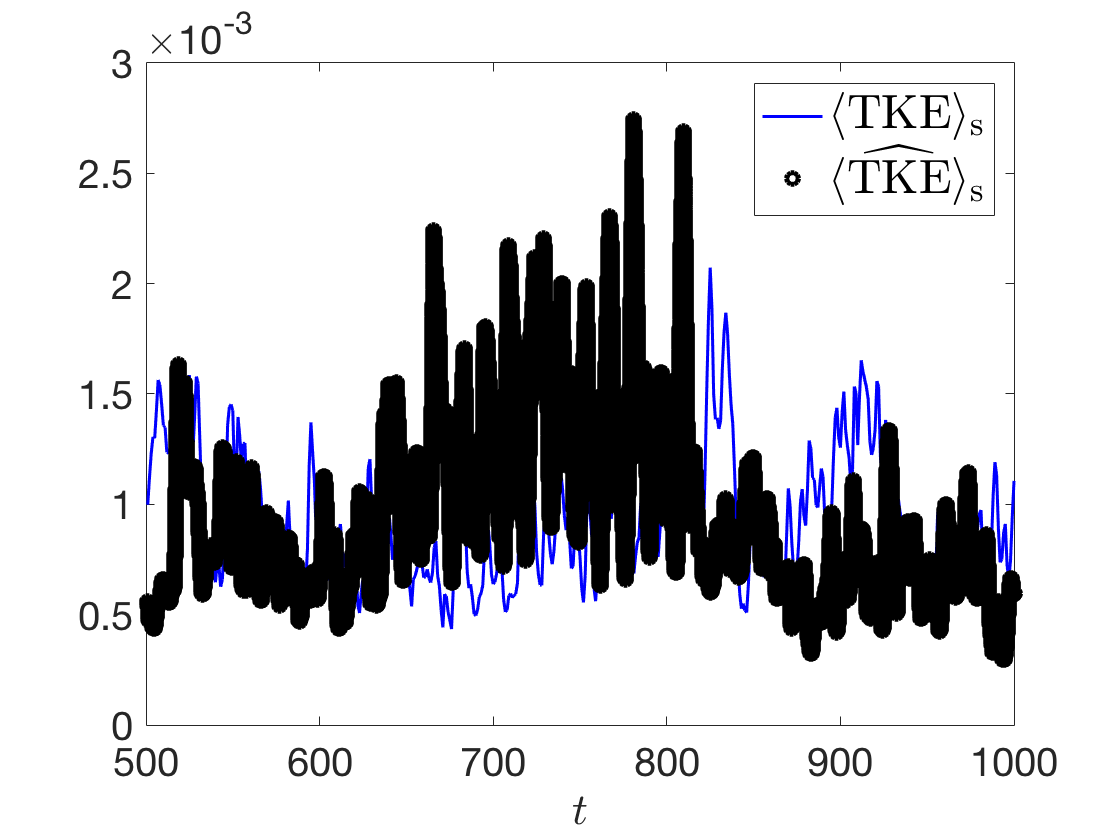}}
 ~~
\subfloat[$N=60$]  
{  \includegraphics[width=0.30\textwidth]
 {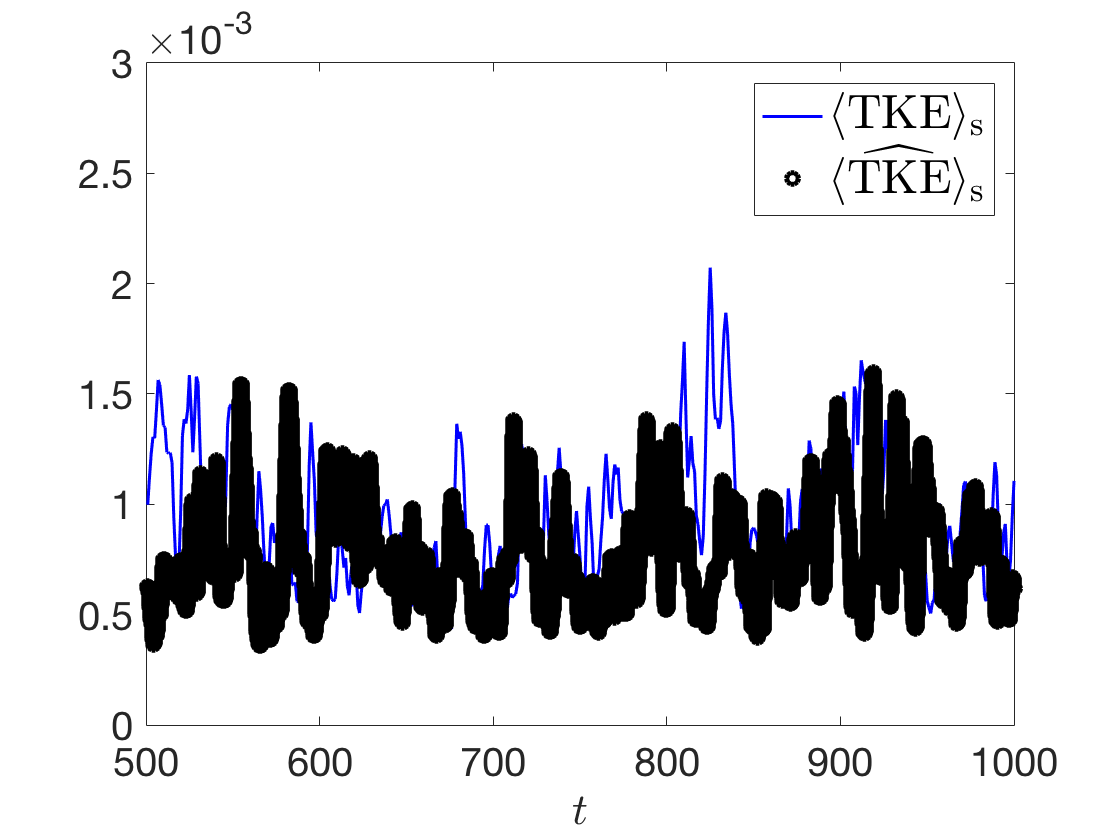}} 
 
 \caption{The solution reproduction problem; 
 constrained POD-Galerkin.
Behavior of the TKE as a function of time for three values of $N$.
$\langle \widehat{{\rm TKE}} \rangle_{\rm s} = 8.6 \cdot 10^{-4}$ ($N=20$),
$9.4 \cdot 10^{-4}$ ($N=40$),
$7.7 \cdot 10^{-4}$ ($N=60$).
$V_{\rm s} (\widehat{{\rm TKE}} )  = 5.5 \cdot 10^{-9}$ ($N=20$),
$1.7 \cdot 10^{-7}$ ($N=40$),
$5.8 \cdot 10^{-8}$ ($N=60$).
($\langle {\rm TKE} \rangle_{\rm s} = 9.4 \cdot 10^{-4}$, 
$V_{\rm s} ( {\rm TKE} ) = 8.5 \cdot 10^{-8}$)
 (${\rm Re} = 15000$, $\epsilon=0.01$).
  }
 \label{fig:cGal_energy_time}
  \end{figure}

Figure \ref{fig:cgal_numerics}(a) shows the behavior of the relative $L^2$ and $H^1$ errors in the mean flow prediction with $N$,
while Figure \ref{fig:cgal_numerics}(b)  shows the behavior of the mean TKE with $N$.
We observe that the constrained formulation leads to a substantial improvement in performance compared to the standard POD-Galerkin method (cf. Figures \ref{fig:pod_gal_numerics}(b) and \ref{fig:pod_gal_numerics}(c)):
for $N \gtrsim 40$ the relative error in the mean is less than $2 \%$, while the predicted mean TKE is bounded from above by $\langle {\rm TKE} \rangle_{\rm s}$ for all values of $N$. 
Furthermore, we observe that the TKE of our constrained Galerkin formulation is larger than the one predicted by the Galerkin ROM for certain values of $N$, and is smaller  for other values of $N$:
this empirically proves that our approach does not necessarily add dissipation to the Galerkin ROM.
In Figure \ref{fig:cgal_numerics_time}, we repeat the tests of Figure \ref{fig:cgal_numerics} for $\Delta t' = 0.5 \Delta t = 2.5 \cdot 10^{-3}$.
We observe that results are consistent with the results shown in Figure  \ref{fig:cgal_numerics}:
this provides empirical evidence for the stability of our constrained formulation under time-step refinement.
Figure \ref{fig:cgal_numerics_2} shows the behavior  of the sample mean and sample variance of the coefficients 
 $\{  a_n^j \}_j$ for three different values of $N$.  Also in this case, we observe a substantial improvement in performance compared to POD-Galerkin, particularly for high modes.   Finally, Figure \ref{fig:cGal_energy_time} shows the behavior of the TKE as a function of time for three values of $N$.   We observe that   for $N=40$ and $N=60$ the predicted TKE is in good qualitative agreement with the truth;   in addition, predictions of first-  and second-order moments (reported in the caption) are accurate.
In Appendix  \ref{sec:constrained_stability}, we present additional results to demonstrate the efficiency of the constrained formulation, and also the robustness with respect to the choice of $\epsilon$.

\section{The parametric problem}
\label{sec:parametric_problem}
We  consider the extension of our MOR approach to the parametric context.
For the purpose of exposition, we focus our discussion on the lid-driven cavity problem presented in section \ref{sec:lid_driven_cavity_problem}:
we wish to estimate the solution to \eqref{eq:lid_driven_cavity} for ${\rm Re} \in \mathcal{P}=[15000,25000]$.
In view of the $h$-refinement, we introduce the partition of $\mathcal{P}$, 
$\mathcal{I}_1,\ldots,\mathcal{I}_M$ such that
$\bigcup_{m=1}^M \mathcal{I}_m = \mathcal{P}$, 
$\mathcal{I}_m \cap \mathcal{I}_{m'}= \emptyset$.
We seek an estimate of the lifted velocity field $\mathring{u}= \mathring{u}(x,t; {\rm Re})$ of the form 
$\hat{u}(x,t; {\rm Re}) = \sum_{n=1}^N a_n^m(t; {\rm Re}) \zeta_n^m(x)$ for all ${\rm Re} \in \mathcal{I}_m$.
The approach can be trivially extended to other parametric problems that do not involve geometric variations; as already mentioned in section \ref{sec:reproduction_problem}, 
the extension to the latter case is beyond the scope of the present work.
Algorithm \ref{parametric_problem} summarizes the general offline/online paradigm for the parametric problem.
We highlight that, for the sake of generality, in Algorithm \ref{parametric_problem} we distinguish between $L$ (number of offline solves) and $M$ (number of reduced spaces). However, in this work we consider the case $L=M$.

In order to tackle the parametric problem outlined above, we should address two challenges: first, we should extend the constrained formulation to the parametric case; 
second, we should develop a Greedy strategy for the proper selection of the parameters ${\rm Re}_1^{\star},\ldots,{\rm Re}_L^{\star}$,
and the partition $\{ \mathcal{I}_m \}_m$.
We emphasize that a proper selection of the parameters reduces the number of offline full order solves, and is thus crucial for the feasibility of the approach.
In order to address the first challenge, we propose an actionable procedure for the selection of the hyper-parameters
$\{  \alpha_n \}_n$ and $\{  \beta_n \}_n$ associated with the constrained formulation \eqref{eq:constrained_galerkin} in the parametric case.
On the other hand, the Greedy approach relies on an inexpensive error indicator, which corresponds to the dual norm of the residual associated with the time-averaged momentum equation.

\begin{algorithm}[H]                      
\caption{Offline/online paradigm for the parametric problem}     
\label{parametric_problem}                           

\begin{flushleft}
\textbf{Task:}
find an estimate of the lifted velocity field $\mathring{u}= \mathring{u}(x,t; {\rm Re})$ of the form
$\hat{u}(x,t; {\rm Re}) = \sum_{n=1}^N a_n^m(t; {\rm Re}) \zeta_n^m(x)$ for all ${\rm Re} \in \mathcal{I}_m$, $m=1,\ldots,M$.
\end{flushleft}  

\textbf{Offline stage}
\medskip

\begin{algorithmic}[1]

\State
Generate the DNS data $\{ \mathring{u}^k({\rm Re_{\ell}}^{\star}) := \mathring{u}(t_{\rm s}^k;  {\rm Re_{\ell}^{\star}}) \}_{k=1}^K \subset V$,
and ${\rm Re}_1^{\star},\ldots,{\rm Re}_L^{\star} \in \mathcal{P}$.
\vspace{4pt}

\State
Generate the 
partition $\{ \mathcal{I}_m \}_{m}$ of $\mathcal{P}$, 
and the 
reduced spaces $\mathcal{Z}_m^{\rm u} = {\rm span} \{  \zeta_n^m \}_{n=1}^N$, $m=1,\ldots,M$.
\vspace{4pt}

\State
Formulate the Reduced Order Models for each subregion.
\vspace{4pt}
\end{algorithmic}

\medskip

\textbf{Online stage}
\medskip

\begin{algorithmic}[1]

\State
Given ${\rm Re} \in \mathcal{P}$, 
find $m \in \{ 1,\ldots, M \}$ such that ${\rm Re} \in \mathcal{I}_m$.

\State
Estimate the coefficients $\{ a_n^{m,j}({\rm Re})= a_n^m(t_{\rm g}^j; {\rm Re})  \}_{n=1}^N$
for $j=0,1,\ldots,J$.
\vspace{4pt}

\State
Compute the QOIs 
(e.g., mean flow, TKE,...).
\end{algorithmic}
\end{algorithm}

The section is organized as follows.
In section \ref{sec:pod_greedy} we present the POD-$h$Greedy approach,  in section \ref{sec:cGal_parametric}  we present the ROM formulation, and we discuss the choice of the 
hyper-parameters $\{  \alpha_n \}_n$ and $\{  \beta_n \}_n$.
Then,  in section \ref{sec:error_indicator} we propose the  time-averaged error indicator.
 Finally, in section \ref{sec:numerical_results_pod_greedy}, we present the numerical results for the lid-driven cavity problem.

 \subsection{POD-Greedy algorithm}
\label{sec:pod_greedy}

We first present the POD-$h$Greedy algorithm for the construction of the reduced spaces $\{ \mathcal{Z}_{m}^{\rm u} \}_{m=1}^M$, and the partition $\{ \mathcal{I}_m \}_{m=1}^M$ of $\mathcal{P}$, based on the results of $L$ FOM simulations associated with the parameters
${\rm Re}_1^{\star}, \ldots, {\rm Re}_L^{\star}$. 
The approach is a simplified version of the $h$-refinement procedure for parabolic problems  proposed in \cite{eftang2011hp}. In particular, we here consider $L=M$: this implies that each reduced space $\mathcal{Z}_{\ell}^{\rm u}$ is based on the POD of a single full-order solve.

In view of the presentation of the algorithm, we introduce the discretized parameter space
$\mathcal{P}_{\rm train} = \{ {\rm Re}_i \}_{i=1}^{n_{\rm train}}$, 
${\rm Re}_1 \leq \ldots \leq {\rm Re}_{n_{\rm train}}$,
the integers $L$ and $N$, which fix the maximum number of offline solves and the size of the reduced space $\mathcal{Z}^{\rm u}$, 
the integer $n_{\rm cand}<L$, which is the number of ROM evaluations performed online for a given value of the parameters, 
and the \emph{a posteriori} error indicator
$\Delta^{\rm u}: \bigotimes_{j=0}^J V \times \mathcal{P} \to \mathbb{R}_+$.
The error indicator takes as input a sequence $\{ w^j \}_{j=0}^J \subset V$ and the value of the parameter, and returns an estimate of the error in the prediction of the mean flow; we formally present the indicator in section \ref{sec:error_indicator}. We further introduce the functions
$$
\begin{array}{l}
\displaystyle{
[\{  \zeta_n \}_{n=1}^N]
=
\texttt{POD}_V
\left(
\mathcal{S},
N
\right);
\quad
[\{  \mathring{u}^k({\rm Re}) \}_{k=1}^K ]
=
\texttt{DNS-solver}
\left(
{\rm Re},
\{ t_{\rm s}^k \}_{k=1}^K
\right);
} \\[4mm]
\displaystyle{
[\{  \hat{u}^j({\rm Re}) \}_{j=0}^J ]
=
\texttt{ROM-solver}
\left(
{\rm Re},
\mathcal{Z}^{\rm u}
\right).
}
\\
\end{array}
$$
$\texttt{POD}_V$ takes as input the set of snapshots
$\mathcal{S}= \{  w^i  \}_{i=1}^{|\mathcal{S}| }$ and an integer $N>0$, and returns the
orthonormalized
 first $N$ POD  eigenmodes  (see section \ref{sec:galerkin_POD});
on the other hand, $\texttt{DNS-solver}$ takes as input the value of the Reynolds number and the sampling times 
$\{ t_{\rm s}^k \}_{k=1}^K$, and returns the instantaneous velocity at times $\{ t_{\rm s}^k \}_{k=1}^K$;
finally, 
$\texttt{ROM-solver}$ takes as input the value of the Reynolds number 
and the reduced space $\mathcal{Z}^{\rm u}$, 
and returns the ROM solution
$\hat{u}^j = \sum_{n=1}^N \, a_n^j \zeta_n$ for each time step of the grid 
$\{ t_{\rm g}^j \}_{j=0}^J$.
Algorithm \ref{POD_hGreedy} presents the computational procedure for both offline and online stage.
With some abuse of notation, we use  $\Delta_{\ell}^{\rm u}(\cdot)$ to refer to the error estimate associated with the $\ell$-th model. 

\begin{algorithm}[H]
 \caption{POD-$h$Greedy algorithm for the construction of  $\{  \mathcal{Z}_{\ell}^{\rm u},  \mathcal{I}_{\ell}  \}_{\ell}$}
  \label{POD_hGreedy}
 
\textbf{Offline stage:} $[ \{ \mathcal{Z}_{\ell}^{\rm u} \}_{\ell=1}^L  ] = $ \texttt{Offline} $(\mathcal{P}_{\rm train}, N,L,\Delta^{\rm u}, \{ t_{\rm s}^k \}_{k=1}^K)$. 

\small 
\emph{Inputs:}  $\mathcal{P}_{\rm train} = \{ {\rm Re}_i \}_{i=1}^{n_{\rm train}}=$ discretized parameter space, $N=$ dimension of each reduced space,
$L=$ maximum number of offline solves, $\Delta^{\rm u}=$ error indicator, $\{ t_{\rm s}^k \}_{k=1}^K=$ sampling times.
 
 \emph{Output:} $\{ (\mathcal{Z}_{\ell}^{\rm u}, {\rm Re}_{\ell}^{\star} \}_{\ell=1}^L=$  reduced space/anchor point pairs.
 
 \normalsize 
 
 \begin{algorithmic}[1]
 \State
${\rm Re}_1^{\star} = {\rm rand} (  \mathcal{P}_{\rm train}   )$
\smallskip

\For{$\ell=1,\ldots,L$ }
\State
$[\{  \mathring{u}^k({\rm Re}_{\ell}^{\star}) \}_{k=1}^K ]
= \texttt{DNS-solver}
\left( {\rm Re}_{\ell}^{\star}, \{ t_{\rm s}^k \}_{k=1}^K \right)$
\smallskip

\State
$[\{  \zeta_n^{\ell} \}_{n=1}^{N} ]
=
\texttt{POD}_V
\left(
\{  \mathring{u}^k({\rm Re}_{\ell}^{\star}) \}_{k=1}^K,
N
\right)$
 \smallskip

\State
Define $\mathcal{Z}_{\ell}^{\rm u}= {\rm span} \{ \zeta_n^{\ell} \}_{n=1}^{N}$,
build the ROM  structures (cf. sections \ref{sec:galerkin_ROM_formulation} and \ref{sec:error_indicator}).

\For { $i=1,\ldots,n_{\rm train}$ }
\State
$[\{  \hat{u}_{\ell}^j({\rm Re}_i) \}_{j=0}^J ]
=
\texttt{ROM-solver}
\left(
{\rm Re}_i,
\mathcal{Z}_{\ell}^{\rm u}
\right)$;

\State
Compute the error estimate
$\Delta_{\ell}^{\rm u}( {\rm Re}_i )$
\EndFor

\State
${\rm Re}_{\ell+1}^{\star}
=
{\rm arg} \max_{{\rm Re} \in \mathcal{P}_{\rm train}} \,
\min_{\ell'=1,\ldots,\ell}
\Delta_{\ell'}^{\rm u}(  {\rm Re} )$.
\EndFor
\end{algorithmic}
  
  \textbf{Online stage:} $[ \{ \hat{u}^j \}_j    ] = $ 
  \texttt{Online} 
  $(\{ (\mathcal{Z}_{\ell}^{\rm u}, {\rm Re}_{\ell}^{\star} \}_{\ell=1}^L, \Delta^{\rm u}, n_{\rm cand}, {\rm Re})$. 

\small 
\emph{Inputs:}  
$\{ (\mathcal{Z}_{\ell}^{\rm u}, {\rm Re}_{\ell}^{\star} \}_{\ell=1}^L = $  reduced space/anchor point pairs,
$n_{\rm cand}=$ online ROM evaluations, 
${\rm Re}=$ input parameter.

\emph{Output:}  
$\{ \hat{u}^j \}_j =$  solution estimate.

 \normalsize 
  
\begin{algorithmic}[1]
\State
Find the $n_{\rm cand}$ nearest anchors to ${\rm Re}$:
  ${\rm Re}_{(1)}^{\star}, \ldots, {\rm Re}_{(n_{\rm cand})}^{\star}$
\smallskip
\For{$i=1,\ldots,n_{\rm cand}$}
\State
$[ \{  \hat{u}_{(i)}^j  \}_{j=0}^J ]
=
\texttt{ROM-solver}
\left(
{\rm Re},
\mathcal{Z}_{(i)}^{\rm u}
\right)$

\State
Compute the error estimate
$\Delta_{(i)}^{\rm u}( {\rm Re})$
\EndFor
  
\State  
  Return 
  $\{ \hat{u}^j = \hat{u}_{(i^{\star})}^j \}_j$,
  where  $i^{\star}$ is the minimizer of $\{ \Delta_{(i)}^{\rm u}( {\rm Re}) \}_i$.
\end{algorithmic}
   \end{algorithm}

Algorithm \ref{POD_hGreedy} combines a POD in time with a Greedy in parameter.
As explained in the introduction, Greedy techniques are crucial to allow efficient parameter explorations at an affordable offline computational cost.
We emphasize that our approach is different  from the POD-Greedy strategy proposed in \cite{haasdonk2008reduced}:
rather than building $L$ different reduced spaces, the authors of \cite{haasdonk2008reduced} combine data from different parameters to generate  a single reduced space.

We also  remark that our definition of the partition --- see  Algorithm \ref{POD_hGreedy}, Online stage --- 
can be formally expressed as follows:
\begin{equation}
\label{eq:definition_partition}
\mathcal{I}_{\ell}
=
\left\{
{\rm Re} \in \mathcal{P}:
\, \ell \in I({\rm Re}),
\;
\Delta_{\ell}^{\rm u}({\rm Re}) < \Delta_{\ell'}^{\rm u}({\rm Re}), 
\ell'  \in I({\rm Re}), 
\ell' \neq \ell
\right\},
\end{equation}
where $I({\rm Re}) \subset \{1,\ldots,L\}$ is the set of indices associated with the $n_{\rm cand}$ nearest anchor points
to ${\rm Re}$.

\subsection{Constrained Galerkin formulation}
\label{sec:cGal_parametric}

Given the reduced spaces
$\{ \mathcal{Z}_{\ell}^{\rm u}= {\rm span} \{ \zeta_n^{\ell} \}_{n=1}^N    \}_{\ell=1}^L$,   we consider the constrained Galerkin formulation proposed in section \ref{sec:reproduction_problem}:
given ${\rm Re} \in \mathcal{P}$, 
and the time grid $\{  t_{\rm g}^j \}_{j=0}^J$, 
the $\ell$-th ROM seeks the coefficients $\{ \mathbf{a}_{\ell}^j \}_{j=0}^J \subset \mathbb{R}^N$
such that
\begin{equation}
\label{eq:constrained_galerkin_param}
\begin{array}{l}
\displaystyle{
\mathbf{a}_{\ell}^{j+1}:=
{\rm arg} \min_{\mathbf{a} \in \mathbb{R}^N}
\,
\|
\mathbb{A}_{\ell}(\mathbf{a}^j; {\rm Re})
\mathbf{a} 
-
\mathbf{F}_{\ell}(\mathbf{a}^j; {\rm Re})
\|_2^2,
\; \; \;
{\rm s.t.}
\; \;
\alpha_n^{\ell}({\rm Re}) \leq a_n \leq \beta_n^{\ell}({\rm Re}),
}
\\
\hfill
n=1,\ldots,N;
\\
\end{array}
\end{equation}
where  $\mathbb{A}_{\ell}$ and $\mathbf{F}_{\ell}$ can be computed by exploiting 
\eqref{eq:galerkin_algebraic} for the reduced space
$\mathcal{Z}_{\ell}^{\rm u}$, and the constraints 
$\{ \alpha_n^{\ell} \}_{n,\ell}$ and $\{ \beta_n^{\ell} \}_{n,\ell}$ are based on the DNS data for the anchor point ${\rm Re}_{\ell}^{\star}$. In greater detail, given 
$\ell \in \{ 1,\ldots,L \}$, 
assuming that $\zeta_1^{\ell},\ldots,\zeta_N^{\ell}$ are orthonormal in $V$, we define 
$\{ \alpha_n^{\ell} \}_{n}$ and $\{ \beta_n^{\ell} \}_{n}$
such that
\begin{subequations}
\label{eq:alpha_beta_param}
\begin{equation}
\alpha_n^{\ell}  := m_{n,\ell}^{\rm u} - \epsilon \Delta_{n,\ell}^{\rm u},
\quad
\beta_n^{\ell} := M_{n,\ell}^{\rm u} + \epsilon \Delta_{n,\ell}^{\rm u},
\end{equation}
where
\begin{equation}
m_{n,\ell}^{\rm u}:= \min_k \, a_{n,\ell}^{\rm FOM, k},
\quad
M_{n,\ell}^{\rm u}:= \max_k \, a_{n,\ell}^{\rm FOM, k},
\quad
\Delta_{n,\ell}^{\rm u} :=
M_{n,\ell}^{\rm u} - m_{n,\ell}^{\rm u},
\end{equation}
and $a_{n,\ell}^{\rm FOM, k}:=
(\mathring{u}^k({\rm Re}_{\ell}^{\star}), \zeta_n^{\ell}     )_V$.
\end{subequations}
The offline/online decomposition is equivalent to the one described in section \ref{sec:galerkin_ROM_formulation}.  We omit the details.

We observe that our choices of $\alpha_n^{\ell}$ and $\beta_n^{\ell}$ correspond to a  constant approximation of the functions
\begin{equation}
\label{eq:golden_choice_parameters}
m_{n,\ell}^{\rm FOM, u}({\rm Re}):=
\min_{j=J_0,\ldots,J}
\,
\left(
\mathring{u}^j({\rm Re}),
\zeta_n^{\ell}
\right)_V,
\; \;
M_{n,\ell}^{\rm FOM, u}({\rm Re}):=
\max_{j=J_0,\ldots,J}
\,
\left(
\mathring{u}^j({\rm Re}),
\zeta_n^{\ell}
\right)_V;
\end{equation}
where $\mathcal{I}_{\ell} \subset \mathcal{P}$ is defined in \eqref{eq:definition_partition}. 
We observe that the piece-wise constant approximations 
of $m_{n,\ell}^{\rm FOM, u}$ and $M_{n,\ell}^{\rm FOM, u}$
are justified
by our  Greedy algorithm, which adaptively determines the partition of $\mathcal{P}$ based on the error indicator.
For practical parametrizations, and practical values of $L$ (i.e., 
number of offline solves)
we expect that accurate estimates of 
$m_{n,\ell}^{\rm FOM, u}$ and $M_{n, \ell}^{\rm FOM, u}$
over $\mathcal{P}$ might be out of reach.
Therefore, 
we  here effectively rely on 
(i) the robustness of our constrained approach to perturbations in the value of the hyper-parameters, and 
(ii) the weak sensitivity of the functions $m_{n,\ell}^{\rm ROM, u}$ and $M_{n, \ell}^{\rm ROM, u}$ with respect to the parameter. 
For the lid-driven cavity problem considered in this work, we provide numerical evidence to support these two assumptions in Appendix 
\ref{sec:constrained_stability}.



\subsection{A time-averaged error indicator}
\label{sec:error_indicator}

Given the sequence $\{ w^j \}_{j=0}^J \subset V$ and ${\rm Re} \in \mathcal{P}$, we define the discrete time-averaged residual
$\langle R  \rangle: \bigotimes_{j=0}^J V \times V_{\rm div} \times \mathcal{P} \to \mathbb{R}$ associated with \eqref{eq:galerkin_rom_discretized}:
\begin{subequations}
\label{eq:time_avg_residual}
\begin{equation}
\langle R  \rangle \left(
\{ w^j \}_{j=0}^J, v; \, {\rm Re}
\right)
=
\frac{\Delta t}{T -T_0}
\sum_{j=J_0}^{J-1}
\,
e(\hat{u}^j, \hat{u}^{j+1}, {\rm Re}  )
\end{equation}
where $T=t_{\rm g}^J$,  $T_0=t_{\rm g}^{J_0}$, and
\begin{equation}
\begin{array}{l}
\displaystyle{
e(\hat{u}^j, \hat{u}^{j+1}, {\rm Re}  )
:=
\left(
\frac{ \hat{u}^{j+1} - \hat{u}^j   }{\Delta t}, v
\right)_{L^2(\Omega)}
\,
+
\frac{1}{\rm Re}
( \hat{u}^{j+1} + R_g, v    )_V}
\\[3mm]
\hfill
+
c(\hat{u}^{j} + R_g,  \hat{u}^{j+1} + R_g, v    ),
\,
j=J_0,\ldots,J-1
\\
\end{array}
\end{equation}
\end{subequations}
Then, we define the error indicator 
$\Delta^{\rm u}
: \bigotimes_{j=0}^J V \times   \mathcal{P} \to \mathbb{R}_+$ as follows:
\begin{equation}
\label{eq:time_avg_error_indicator}
\Delta^{\rm u}
\left(
   \{ w^j \}_{j=0}^J; \, {\rm Re}  \right)
:=
\big\|
 \langle R  \rangle \left(
\{ w^j \}_{j=0}^J, \cdot ; \, {\rm Re}
\right)
\big\|_{V_{\rm div}'}
\end{equation}
where 
$ \| \cdot  \|_{V_{\rm div}'}$ denotes the norm of  the dual space
$V_{\rm div}'$.

 In our numerical tests, as in \eqref{eq:mean_flow_def},   we consider $J_0$ such that $t_{\rm g}^{J_0} = T_0=500$:
this choice is designed to limit the effect of the transient  dynamics.
It is easy to verify that the solution to the FOM 
for any initial condition
--- provided that the same time discretization is employed ---
 satisfies $\Delta^{\rm u}  \equiv 0$.
This implies that two sequences
$\{ w^j \}_{j=0}^J,  \{ \tilde{w}^j \}_{j=0}^J   \subset V$ satisfying $\Delta^{\rm u}  \equiv 0$ might be far from each other at each time step
(i.e., $\| \tilde{w}^j - w^j  \|_V$ is large for  any $j \geq 0$).
However, for sufficiently large values of $J$, we expect 
 $\Delta^{\rm u}$  to be highly-correlated with the error in the mean flow prediction; for this reason, we can exploit  $\Delta^{\rm u}$ to guide the Greedy algorithm presented in section \ref{sec:pod_greedy}.
 We empirically investigate the correlation between $\Delta^{\rm u}$ and the error in the mean flow prediction in the numerical experiments at the end of the section. A theoretical justification of the error indicator is beyond the scope of the present work.

The error indicator $\Delta^{\rm u}$ can be computed efficiently for sequences in  $\mathcal{Z}^{\rm u}$ exploiting an offline/online computational decomposition; the procedure is standard in the Reduced Basis literature, and is reported in
Appendix  \ref{sec:error_indicator_offline_online}.

\begin{remark}
\label{remark_corrected_estimator}
We do not expect that the residual indicator \eqref{eq:time_avg_error_indicator} is in good quantitative agreement with the error in mean flow prediction $\|   \langle u  - \hat{u} \rangle_{\rm g}  \|_V$. More precisely, 
if we define  the effectivity $\eta := \frac{\Delta^{\rm u}}{\|   \langle u  - \hat{u} \rangle_{\rm g}   \|_V }$ of the
residual error indicator, we do not expect that
$\eta$ is close to one.

In order to obtain a quantitative estimate of the error of the ROM anchored in ${\rm Re}^{\star}$, we can consider the corrected estimator
\begin{equation}
\label{eq:corrected_estimator}
\Delta^{\rm u,corr}( 
{\rm Re};  {\rm Re}^{\star}  )
:=
\frac{1}{\eta(  {\rm Re}^{\star}   )}
\,
\Delta^{\rm u}( 
{\rm Re}  ),
\end{equation}
where 
$\Delta^{\rm u}( 
{\rm Re}  ) $ is the error indicator associated with the ROM anchored in ${\rm Re}^{\star}$, and  
$\eta(  {\rm Re}^{\star}   )$ is the effectivity evaluated at $\rm Re =   {\rm Re}^{\star} $.
Note that the computation of $\eta(  {\rm Re}^{\star}   )$ does not require any additional call  to the DNS solver.
\end{remark}

\subsection{Numerical results}
\label{sec:numerical_results_pod_greedy}

\begin{figure}[h!]
\centering
\subfloat[ ]  
{  \includegraphics[width=0.33\textwidth]
 {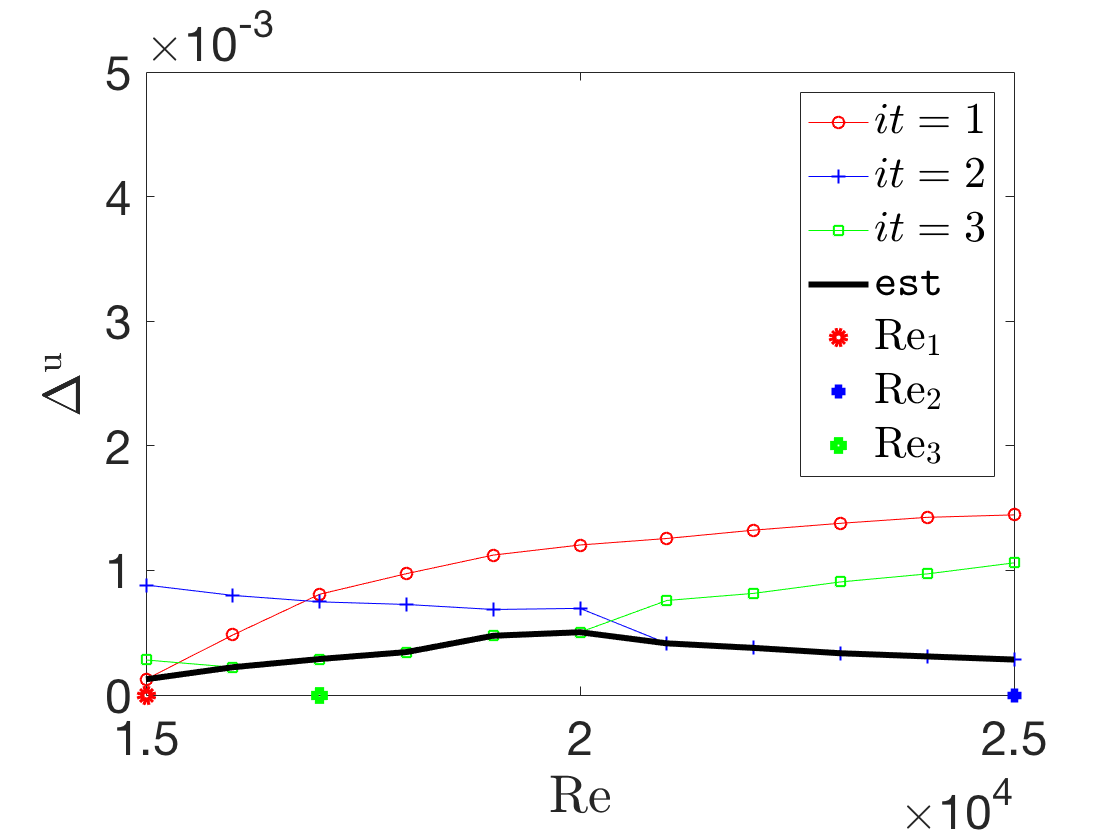}}
 ~~
 \subfloat[]
{  \includegraphics[width=0.33\textwidth]
 {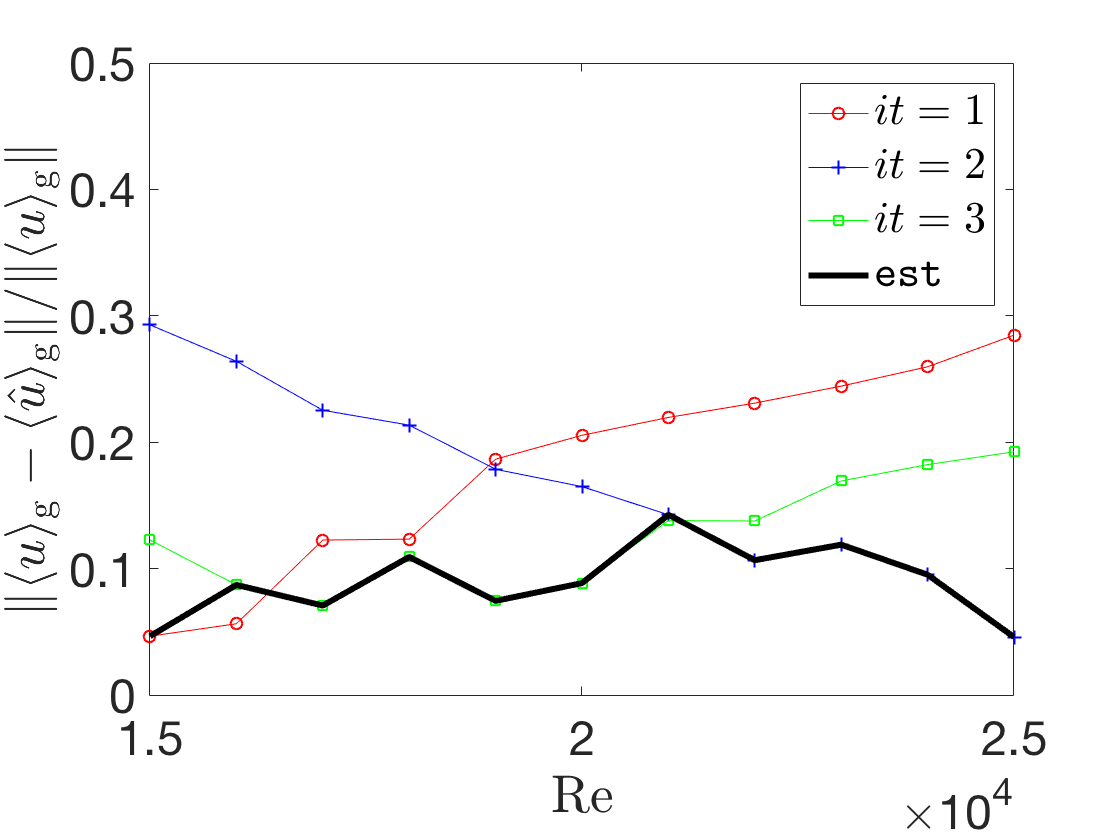}}
 
 \caption{The parametric problem; performance of POD-$h$Greedy.
Figure (a): behavior of $\Delta^{\rm u}$ with ${\rm Re}$ for three iterations.
Figure (b): behavior of the relative $H^1$ error  in mean flow prediction with ${\rm Re}$ for three iterations;
the black line (\texttt{est}) shows the performance of the reduced model which minimizes the error indicator (and thus is selected by the Greedy procedure).
 ($\epsilon = 0.05$, $N=80$,${\rm Re_1}=15000, {\rm Re_2}=25000,
 {\rm Re_3}=17000$). 
  }
 \label{fig:hGreedy_numerics}
  \end{figure}

Figure \ref{fig:hGreedy_numerics} shows the results of the application of Algorithm \ref{POD_hGreedy} for the construction of the ROM for the parametric problem. In order to assess performance, we generate DNS data for $t_{\rm g}^j \in \{0,\ldots,1500 \}$,
$\{  t_{\rm s}^k = 500 + k  \}_{k=1}^{K=1000}$ for 
${\rm Re}= 15000,16000, \ldots,25000$ ($n_{\rm train}=11$ datapoints).
Then, we apply Algorithm \ref{POD_hGreedy} with ${\rm Re}_1^{\star} = 15000$, $N=80$, and $\epsilon=0.05$.
We perform $L=3$ iterations of the Greedy procedure.
Figure \ref{fig:hGreedy_numerics}(a) shows the behavior of  $\Delta^{\rm u}$ with ${\rm Re}$ for the three iterations, while
Figure \ref{fig:hGreedy_numerics}(b) shows the  behavior of the relative $H^1$ error  in mean flow prediction with ${\rm Re}$.
The black continuous line denotes the performance of the reduced model which minimizes the error indicator, and thus is selected by the Greedy procedure (cf. Algorithm \ref{POD_hGreedy}, $n_{\rm cand}=2$).
We observe that the maximum relative error decreases at each iteration, and it is roughly $13\%$ after the third iteration.

\begin{figure}[h!]
\centering
 \subfloat[]
{  \includegraphics[width=0.33\textwidth]
 {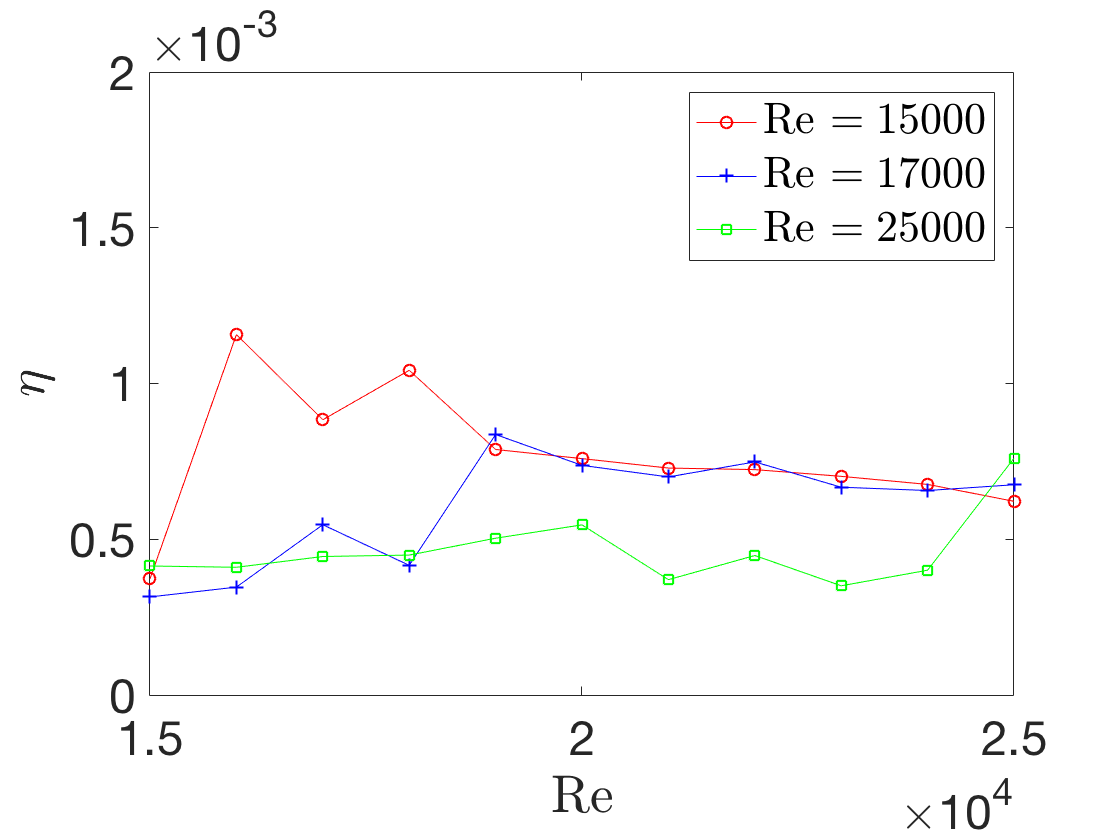}}
 ~~
 \subfloat[]
{  \includegraphics[width=0.33\textwidth]
 {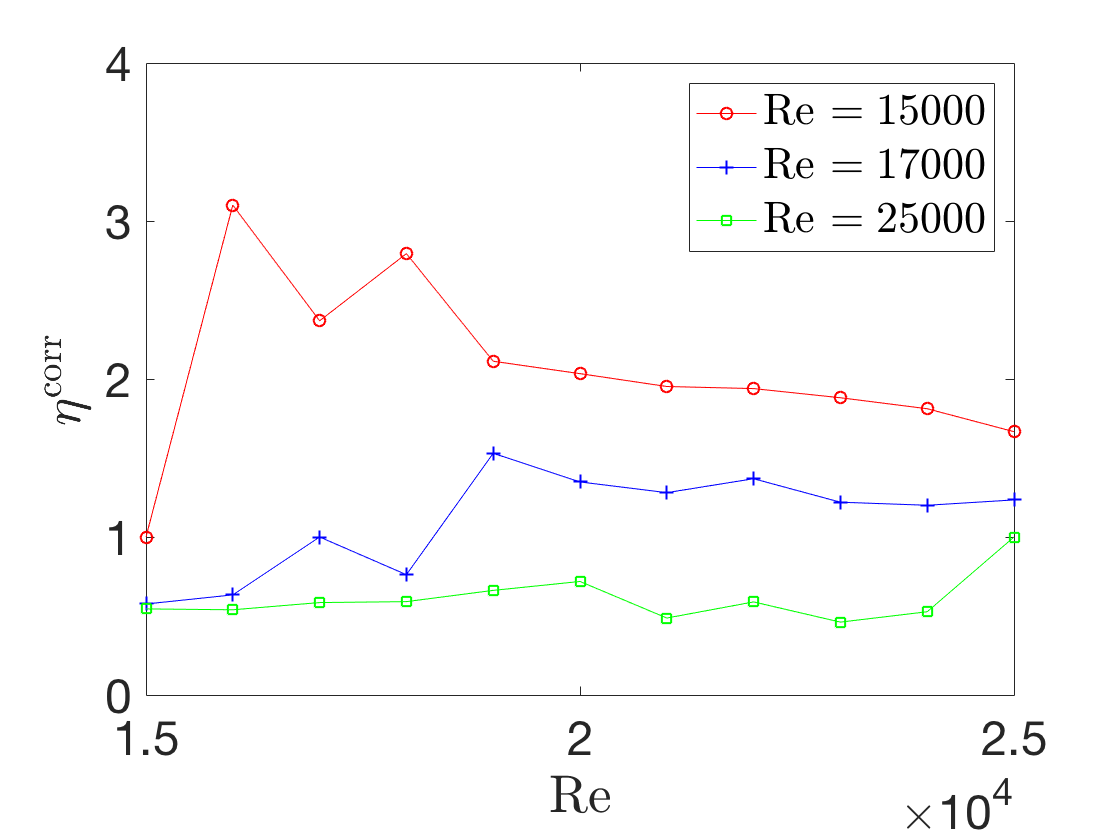}}
 
 \caption{The parametric problem; error estimator.
Figure (a): behavior of the effectivity
$\eta$ of the error indicator  $\Delta^{\rm u}$ \eqref{eq:time_avg_error_indicator}
for the three ROMs.
Figure (b): behavior of the effectivity
$\eta^{\rm corr}$ of the corrected error indicator $\Delta^{\rm u,corr}$ 
\eqref{eq:corrected_estimator}
for the three ROMs.
  }
 \label{fig:a_posteriori}
  \end{figure}

Results of Figure \ref{fig:hGreedy_numerics} show the importance of the error indicator $\Delta^{\rm u}$ in \eqref{eq:time_avg_error_indicator}	 to select the parameters ${\rm Re}_2$ and ${\rm Re}_3$, and also 
motivate  the choice of the partition $\{ \mathcal{I}_{\ell} \}_{\ell}$ in \eqref{eq:definition_partition}:
after the third iteration, for $10$ out of $11$ values of the Reynolds number, the reduced model that minimizes the error indicator (over all models) is the same that minimizes the true error. 
On the other hand, we   observe that the indicator is in poor quantitative agreement with the true error:
Figure \ref{fig:a_posteriori}(a) shows that
the effectivity $\eta$ of the error indicator
is $\mathcal{O}(10^{-3})$ for all three reduced order models and for all values of the Reynolds numbers considered. 
However, Figure \ref{fig:a_posteriori}(b) shows that the correction proposed in Remark \ref{remark_corrected_estimator} leads to an indicator that is in reasonable quantitative agreement with the  error in mean flow prediction.

\begin{figure}[h!]
\centering
\subfloat[ ${\rm Re} = 16000$]  
{  \includegraphics[width=0.33\textwidth]
 {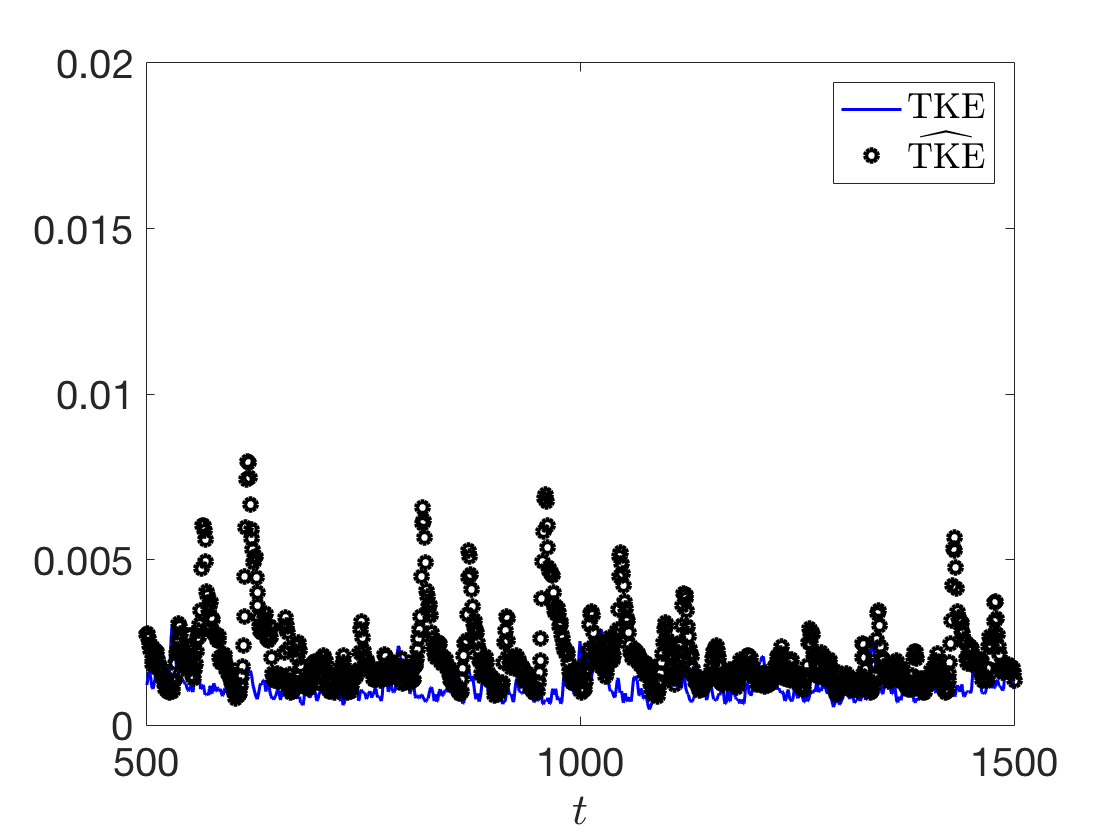}}
 ~~
 \subfloat[${\rm Re} = 20000$]
{  \includegraphics[width=0.33\textwidth]
 {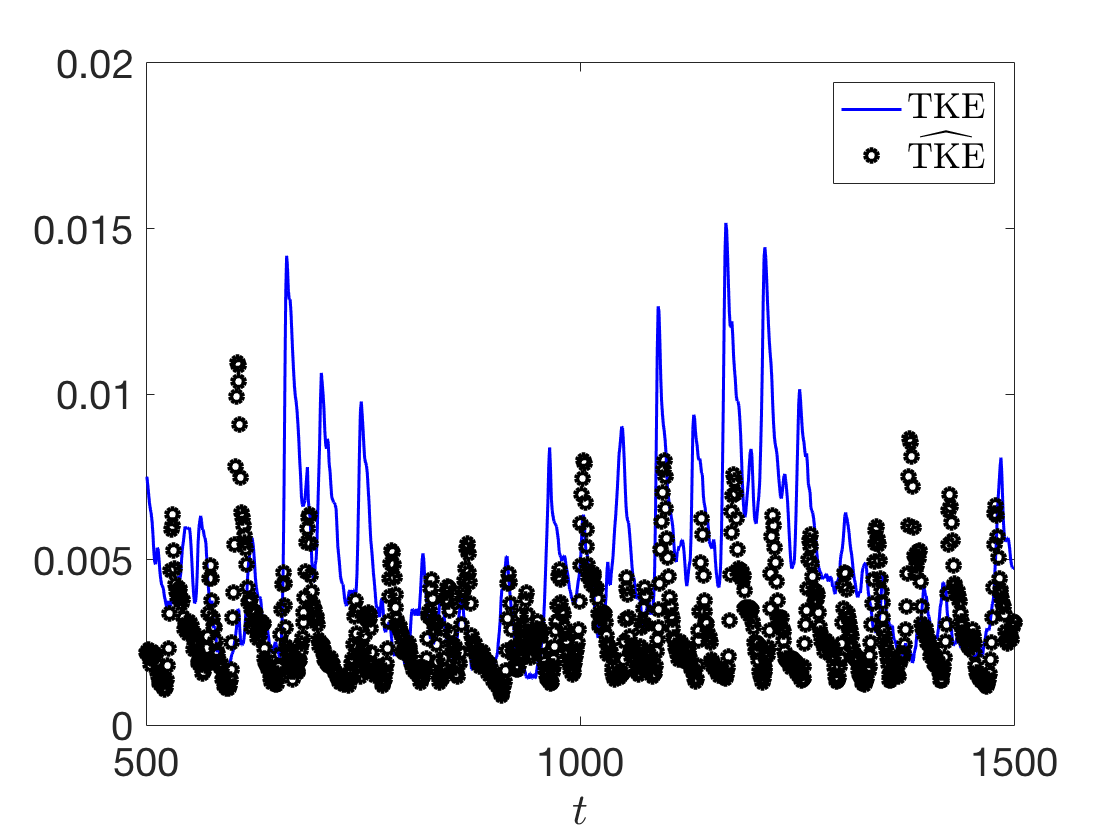}}
  ~~
 \subfloat[${\rm Re} = 23000$]
{  \includegraphics[width=0.33\textwidth]
 {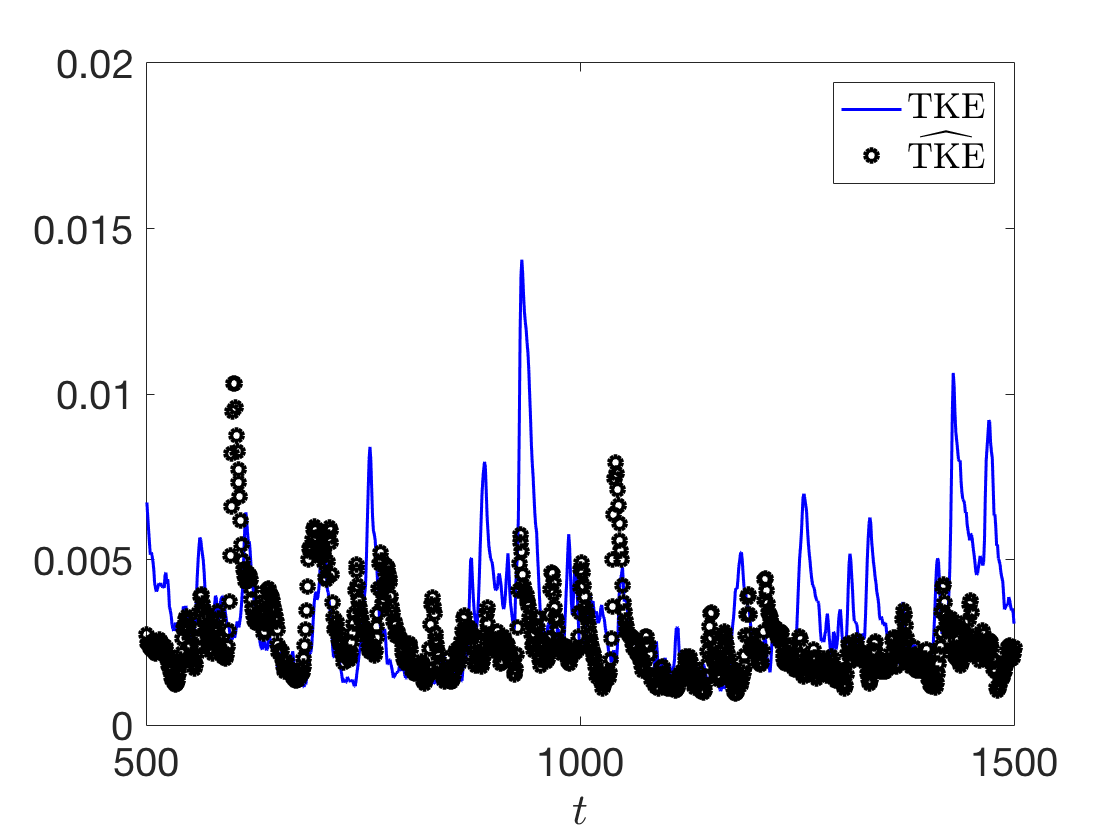}}
 
 \caption{The parametric problem; behavior of the TKE with time for three values of the Reynolds number.
  ($\epsilon = 0.05$, $N=80$, ${\rm Re_1}=15000, {\rm Re_2}=25000,
 {\rm Re_3}=17000$). 
 Predictions for ${\rm Re}_1$ and ${\rm Re}_2$ rely on the ROM anchored in ${\rm Re}^{\star}=17000$, while
 predictions for ${\rm Re}_3$ rely on the ROM anchored in ${\rm Re}^{\star}=25000$.
  }
 \label{fig:hGreedy_TKE}
  \end{figure}

Figure \ref{fig:hGreedy_TKE} shows the behavior of the TKE with time for three values of the Reynolds number, ${\rm Re=} 16000, 20000, 23000$, which have not been selected by the Greedy procedure.
Here, predictions are based on the ROM after three iterations of Algorithm \ref{parametric_problem}.

\section{Conclusions}
\label{sec:conclusions}
In this paper, we present a Reduced Basis technique for long-time integration of turbulent flows.
The three contributions of this work are
(i) a constrained Galerkin formulation that correct the Galerkin statement by incorporating prior information about the long-time attactor, 
(ii) an inexpensive time-averaged  indicator for the error in mean flow prediction,
and
(iii) a POD-$h$Greedy technique for the  construction of the ROM.
In order to assess performance, we apply our approach to a lid-driven cavity problem parametrized with respect to the Reynolds number: first,
we consider the solution reproduction problem (non-predictive case) to demonstrate the effectivity of our new constrained formulation;
second, we consider the parametric problem (predictive case) to validate our error indicator, and more broadly the POD-Greedy procedure.

Our constrained Galerkin formulation is able to accurately predict mean flow and also the TKE. 
The error indicator, despite it is not corroborated by a firm theoretical analysis, is found to be highly-correlated with the error in the prediction of the mean flow; hence, it is naturally suited to drive the offline Greedy.

In this paper, we also highlight a number of challenges, which are particularly relevant for turbulent flows, and that should be taken into consideration in the design of  MOR strategies for turbulent flows: first, the slow convergence of the Kolmogorov $N$-width  suggested by Figure 
\ref{fig:pod_gal_numerics}(a) which prevents us from accurately representing the full dynamics; 
second, the difficulty to combine modes associated with different parameters (cf. Appendix  \ref{sec:p_refinement});
third, the large offline costs both in terms of computational time and required storage.
In this paper, we propose to address the first challenge by \emph{reducing our goal}: rather than trying  to  estimate the full trajectory, we develop a ROM uniquely for the prediction of first and second moments of the long-time dynamics.
Furthermore, we propose to address the second challenge by resorting to an $h$-refinement in parameter.
On the other hand, we here postulate that 
the snapshot set $\{ u^k \}_{k=1}^K$  is   rich enough to accurately estimate the first $N$ POD modes associated with the full trajectory $\{ \mathring{u}^j \}_{j=J_0}^J$,
and also that it is possible to compute and store the Riesz representers $\tilde{\xi}_1,\ldots,\tilde{\xi}_M$, $M = N^2 + 3N + 2$, for residual calculations.
In 
Appendix 
\ref{sec:cv_POD}, we review a 
computational strategy to  assess \emph{a posteriori} the representativity  of our snapshot set;
on the other hand, we refer to a future work for the development of computational strategies to reduce the offline costs related to residual calculations.

We finally outline a number of potential next steps that we wish to pursue in the future.
\begin{itemize}
\item
\emph{Constrained formulation}
Our constrained formulation minimizes the $\ell^2$ error in the reduced Galerkin statement subject to lower and upper bounds for the coefficients of the $N$-term expansion.
We wish to consider other choices both for the objective functions  and for the constraints.
In particular, we wish to minimize the residual in a suitable dual norm,  and we wish to design other constraints to take into account the properties of the attractor.
Furthermore, we also wish to consider the post-processing rectification method
proposed in \cite{Chakir_rectification}
to improve the accuracy 
of the mean flow.
Finally, we wish to consider alternative strategies for  writing the nonlinear term in the momentum equation, and also for  imposing strong boundary conditions.
\item
\emph{hp-Greedy}
In   Appendix
\ref{sec:p_refinement}, we discuss the limitation of the traditional POD-($p$)Greedy algorithm. However, in order to tackle complex parametrizations, we envision that the $h$-refinement strategy proposed in this paper might require an unfeasible number of offline simulations.
This is why we wish to consider more advanced sampling strategies that combine $h$-refinement and $p$-refinement.
\item
\emph{Extension to more challenging problems}
We wish to consider geometry variations, which are particularly relevant for applications.
As explained in the body of the paper, this might 
be accomplished by resorting to the Piola transform, or by considering  a two-field (velocity and pressure) formulation.
Furthermore, we wish to apply our approach to transient problems: in order to face this task, we envision that time-dependent constraints should be considered, and also the time-averaged error indicator should be modified based on the particular quantity of interest we wish to predict.
Finally, we wish to apply our approach to the reduction of LES/URANS flow simulations.  This would substantially increase the range of engineering applications we can tackle with our method.
\end{itemize}

\appendix

\section{Analysis of the solution to the lid-driven cavity problem}
\label{sec:lid_driven_appendix}

\begin{figure}[h!]
\centering
\subfloat[$t=501$]
{  \includegraphics[width=0.30\textwidth]
 {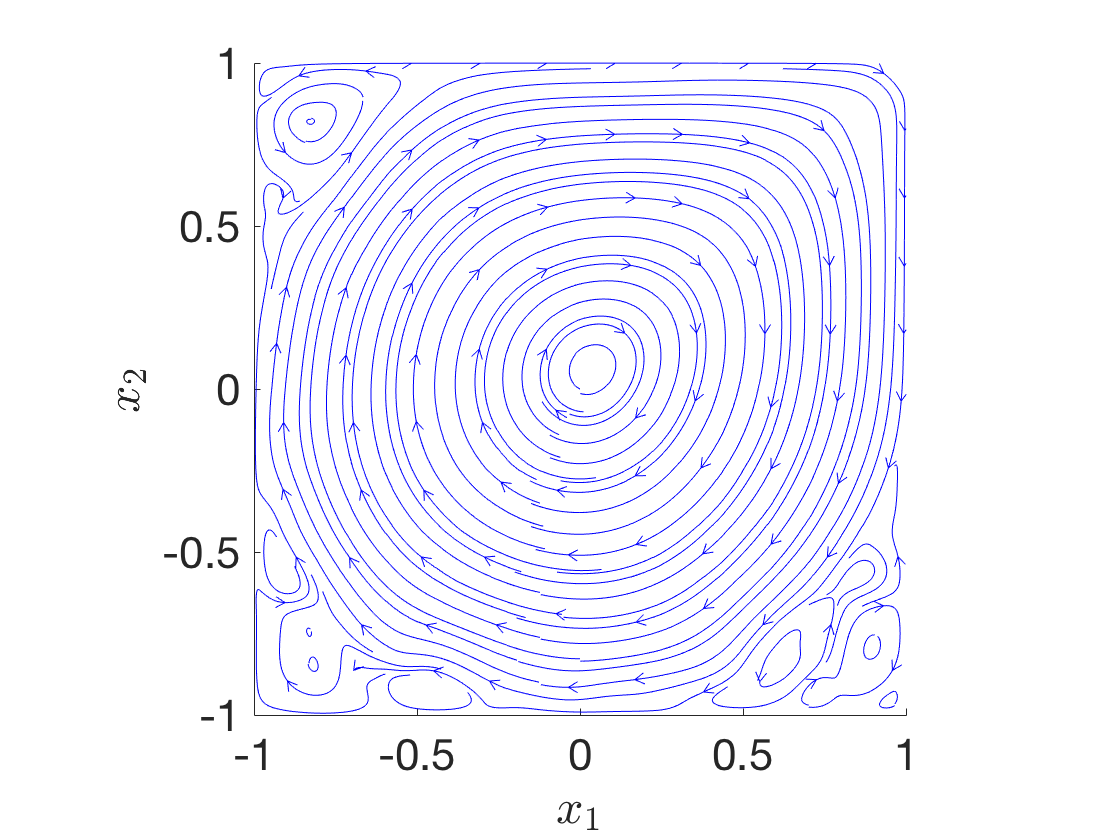}}
 ~~
 \subfloat[$t = 600$]
{  \includegraphics[width=0.30\textwidth]
 {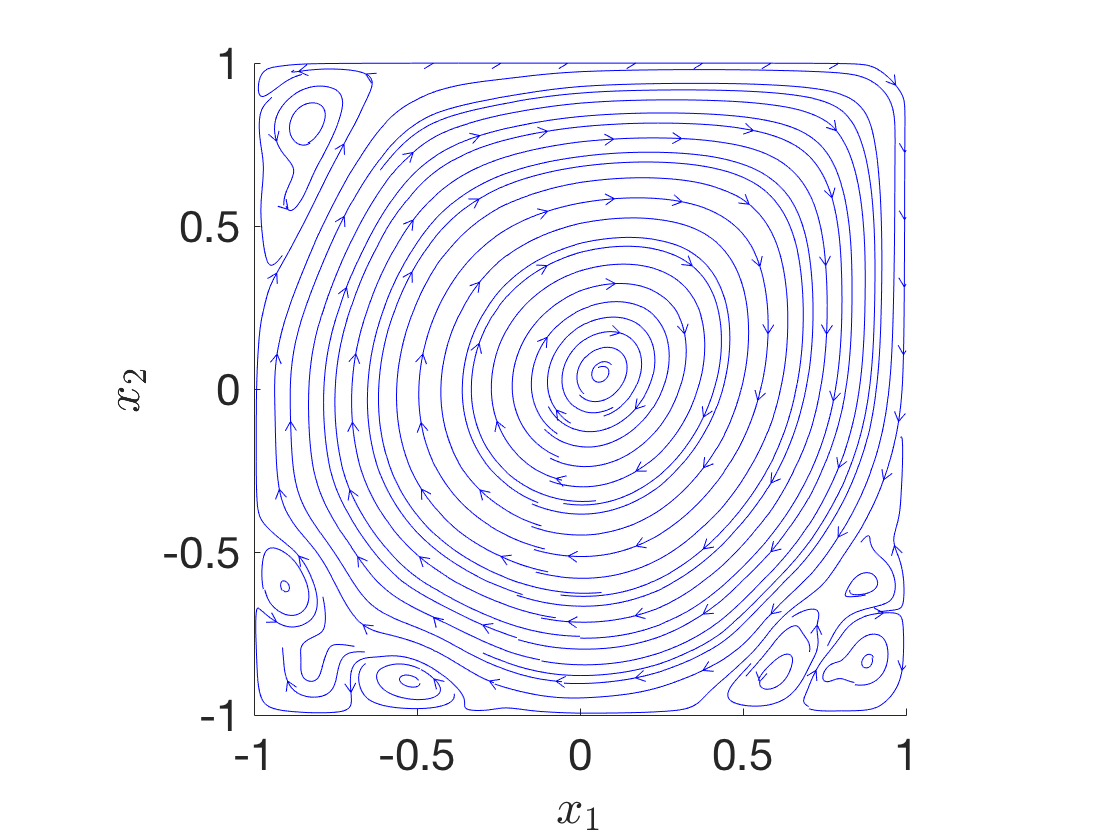}}
~~
 \subfloat[$t = 700$]
{  \includegraphics[width=0.30\textwidth]
 {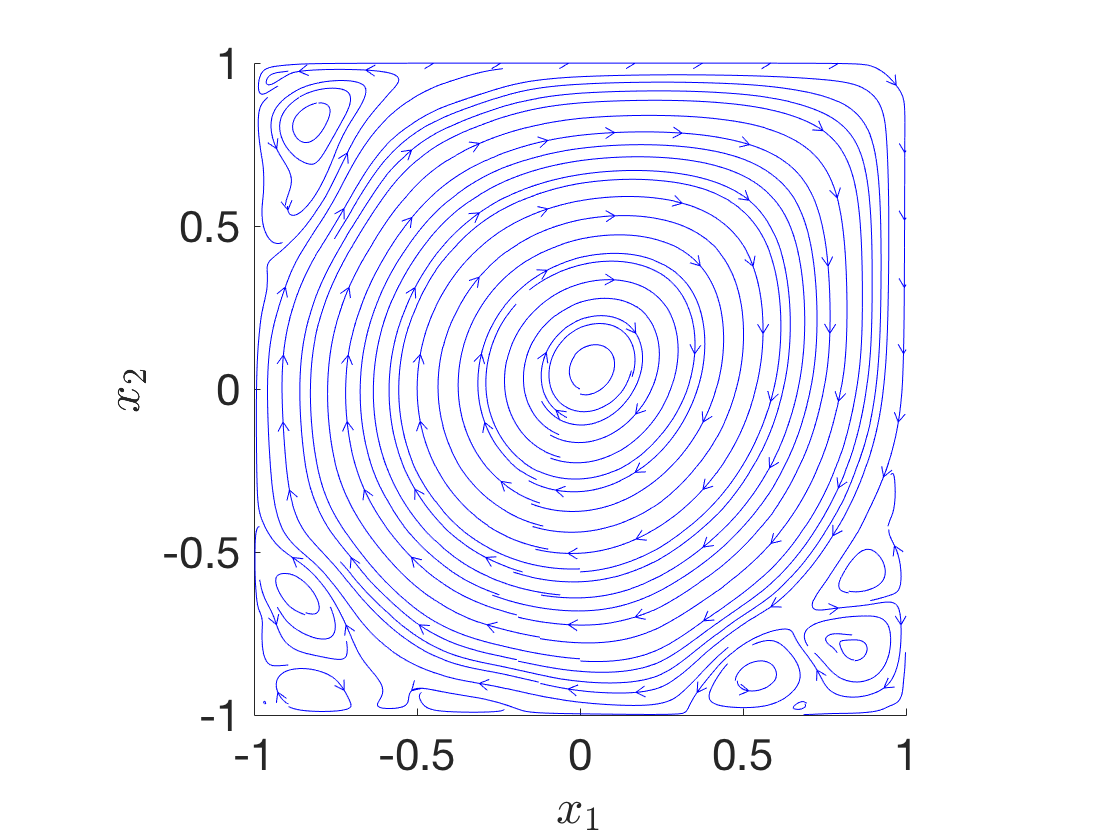}}  
    
 \caption{A lid-driven cavity problem.
 Velocity streamlines for ${\rm Re}= 15000$.
  }
 \label{fig:velocity_streamlines}
  \end{figure}

\begin{figure}[h!]
\centering
\subfloat[$t=1252$]
{  \includegraphics[width=0.30\textwidth]
 {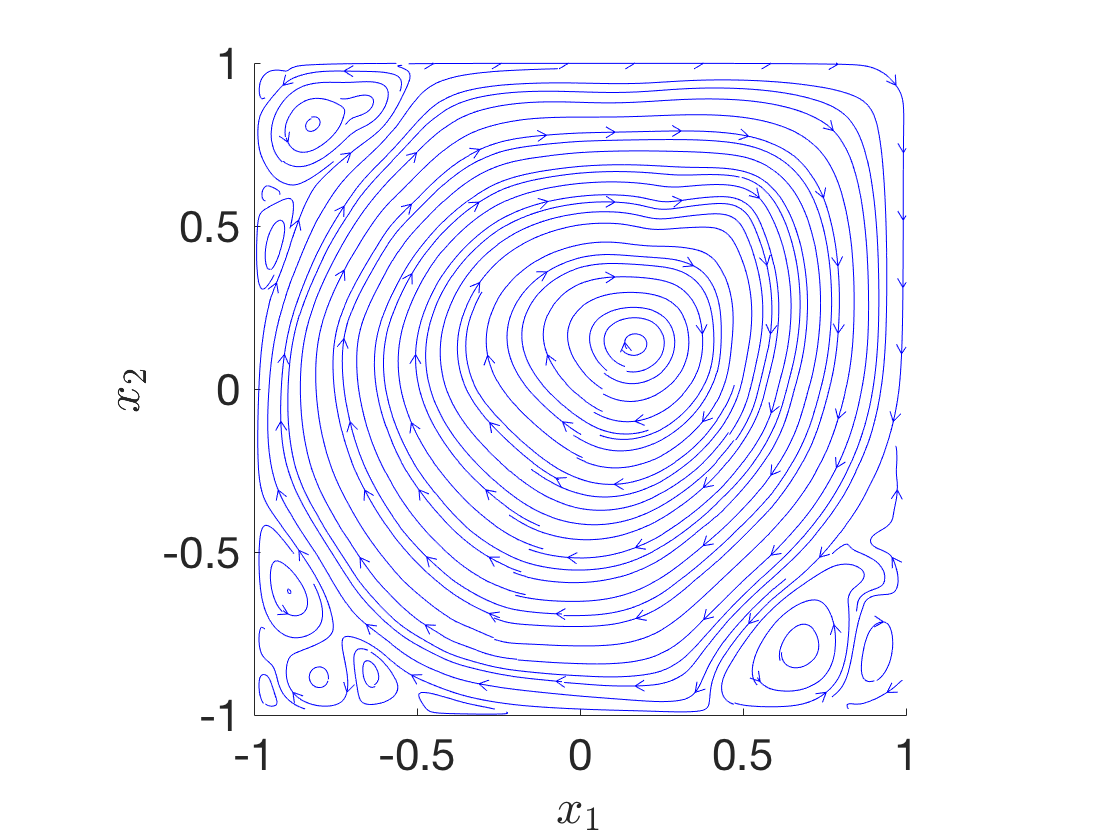}}
 ~~
 \subfloat[$t=1266$]
{  \includegraphics[width=0.30\textwidth]
 {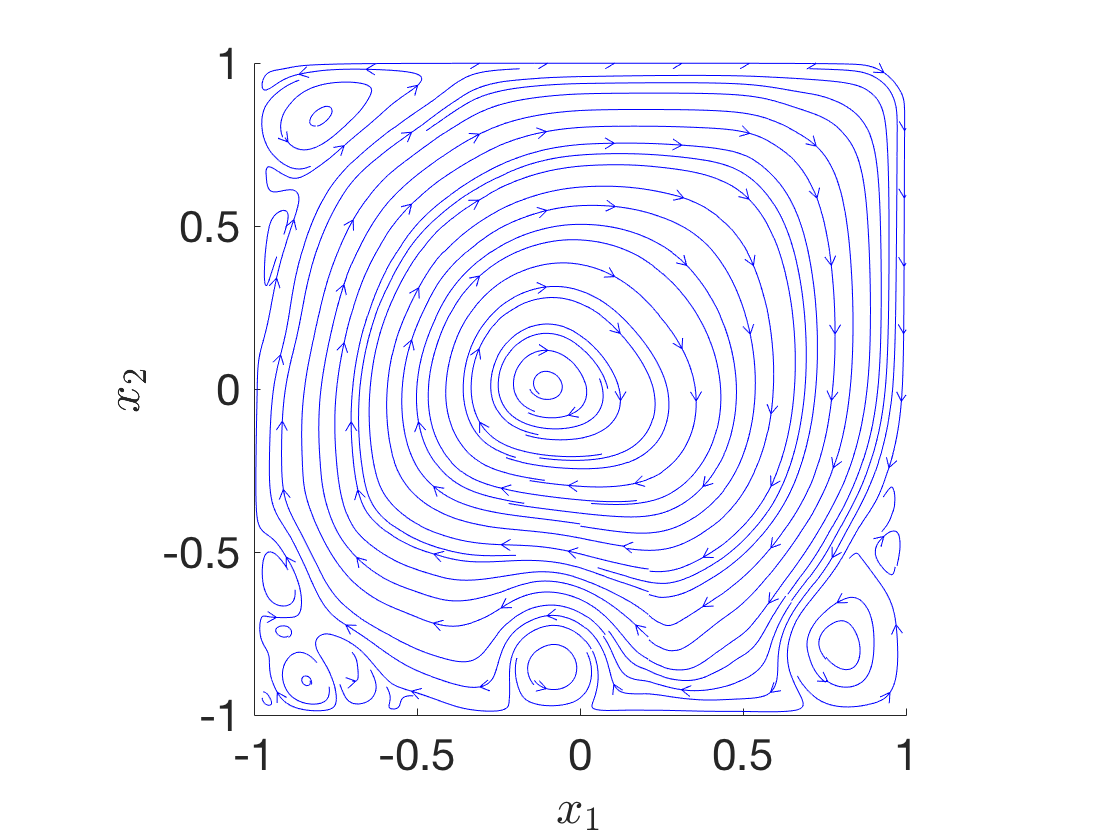}}
 ~~
 \subfloat[$t=1276$]
{  \includegraphics[width=0.30\textwidth]
 {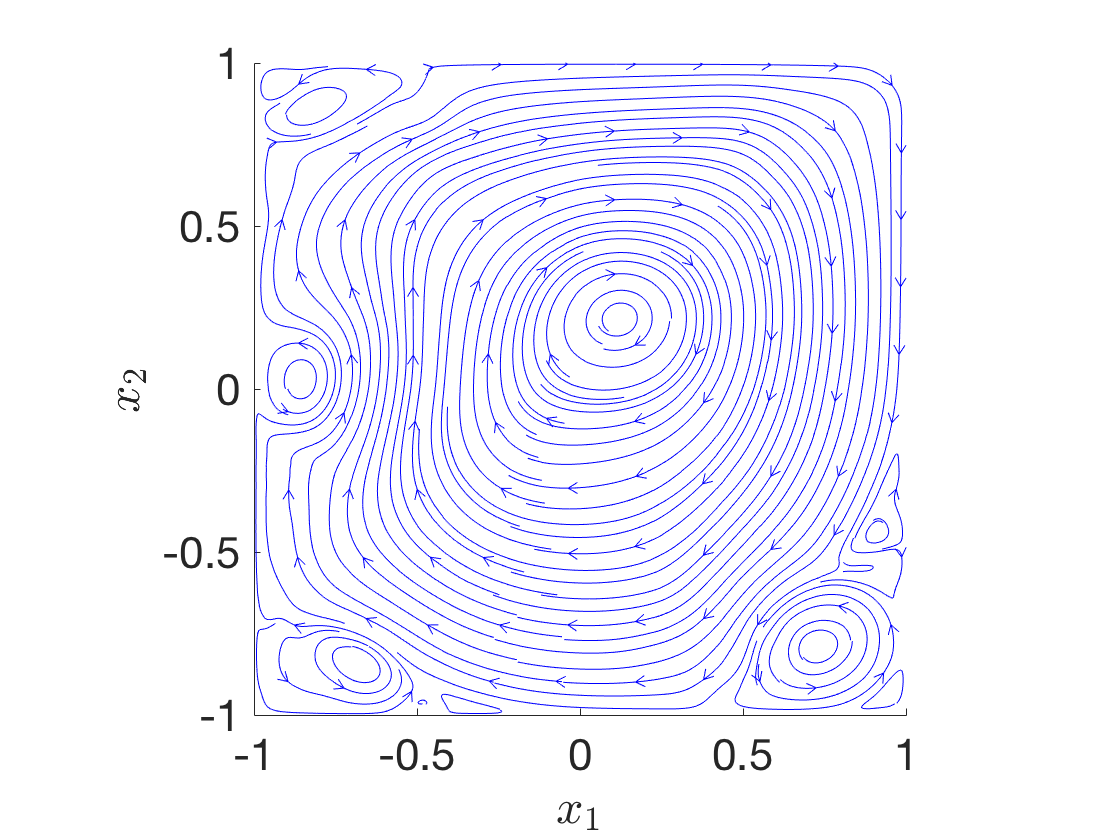}}

\subfloat[$t=1286$]
{  \includegraphics[width=0.30\textwidth]
 {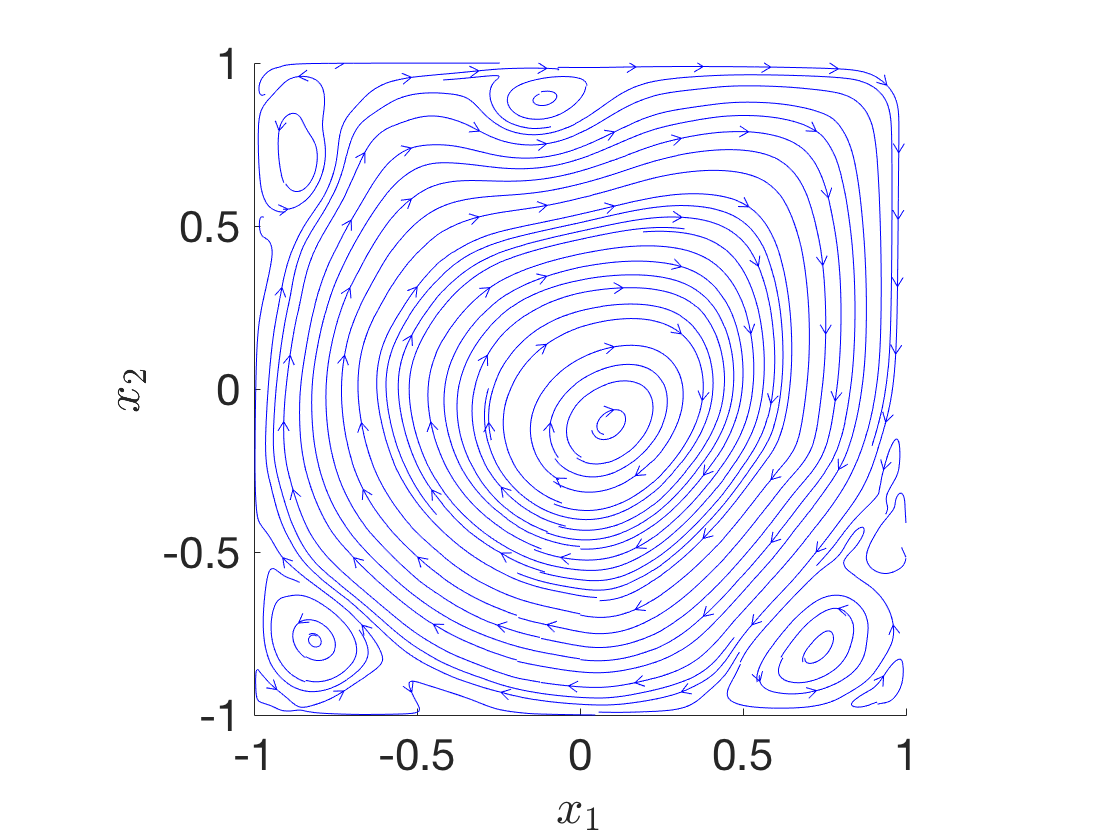}}
 ~~
 \subfloat[$t=1320$]
{  \includegraphics[width=0.30\textwidth]
 {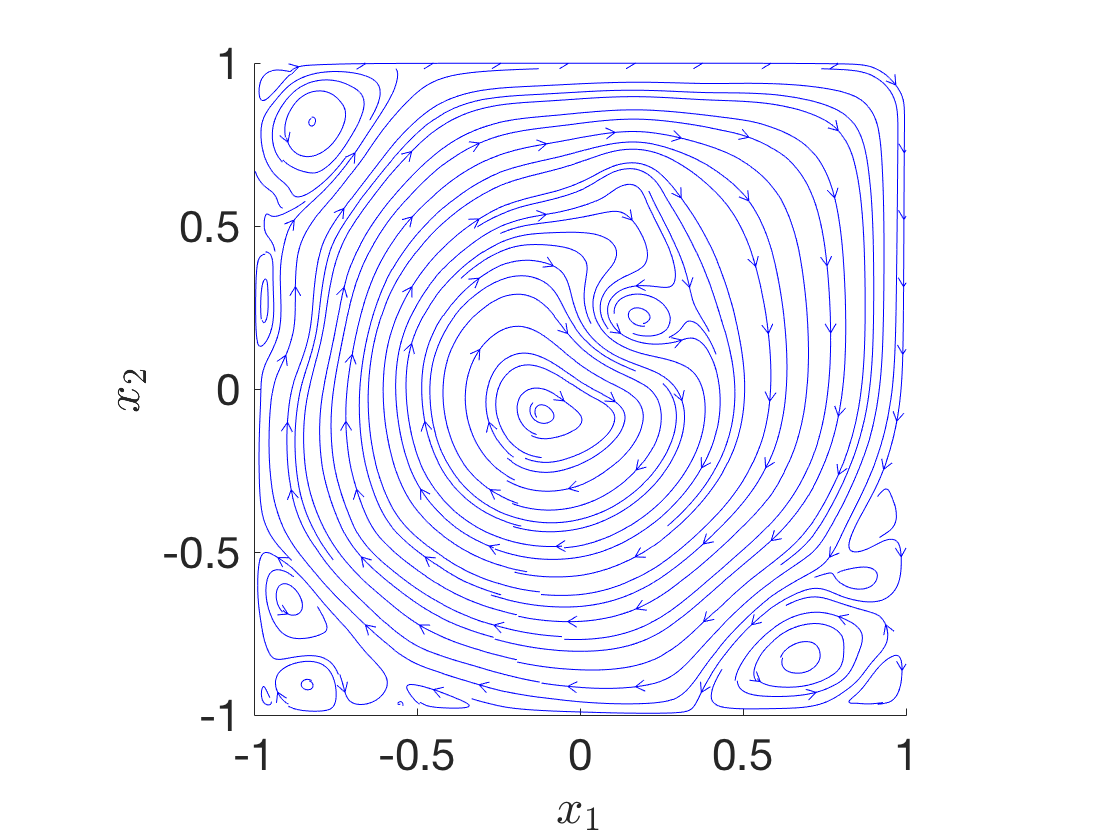}}
 ~~
 \subfloat[$t=1344$]
{  \includegraphics[width=0.30\textwidth]
 {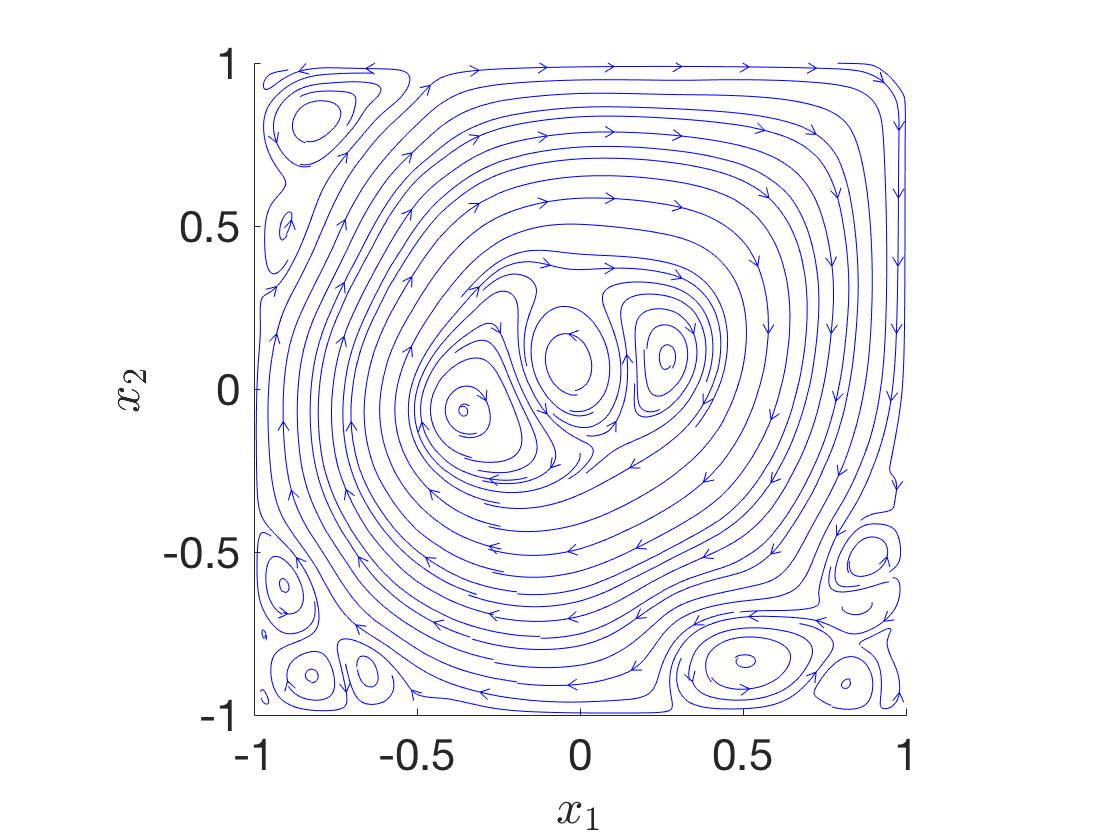}}

 \caption{A lid-driven cavity problem.
 Velocity streamlines for ${\rm Re}= 20000$ for several time steps.
  }
 \label{fig:velocity_streamlines_problem}
  \end{figure}

\begin{figure}[h!]
\centering
\subfloat[${\rm Re}= 15000$]
{  \includegraphics[width=0.30\textwidth]
 {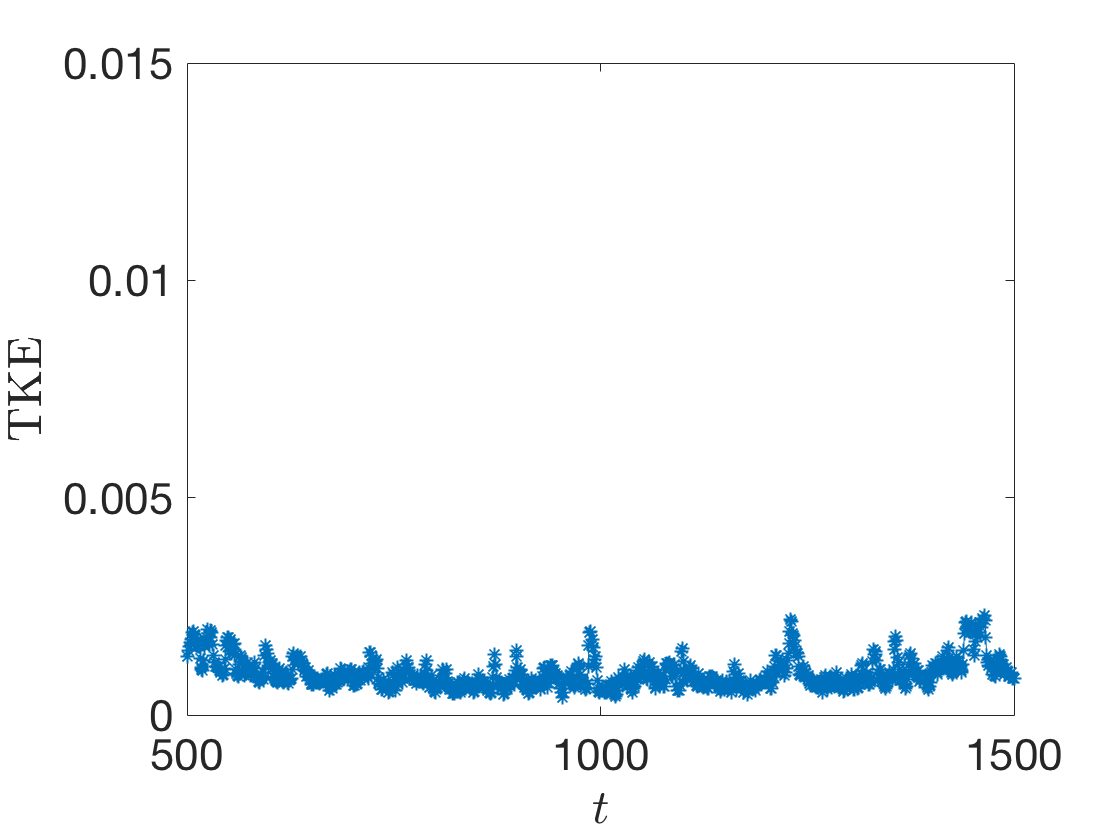}}
 ~~
 \subfloat[${\rm Re}= 20000$]
{  \includegraphics[width=0.30\textwidth]
 {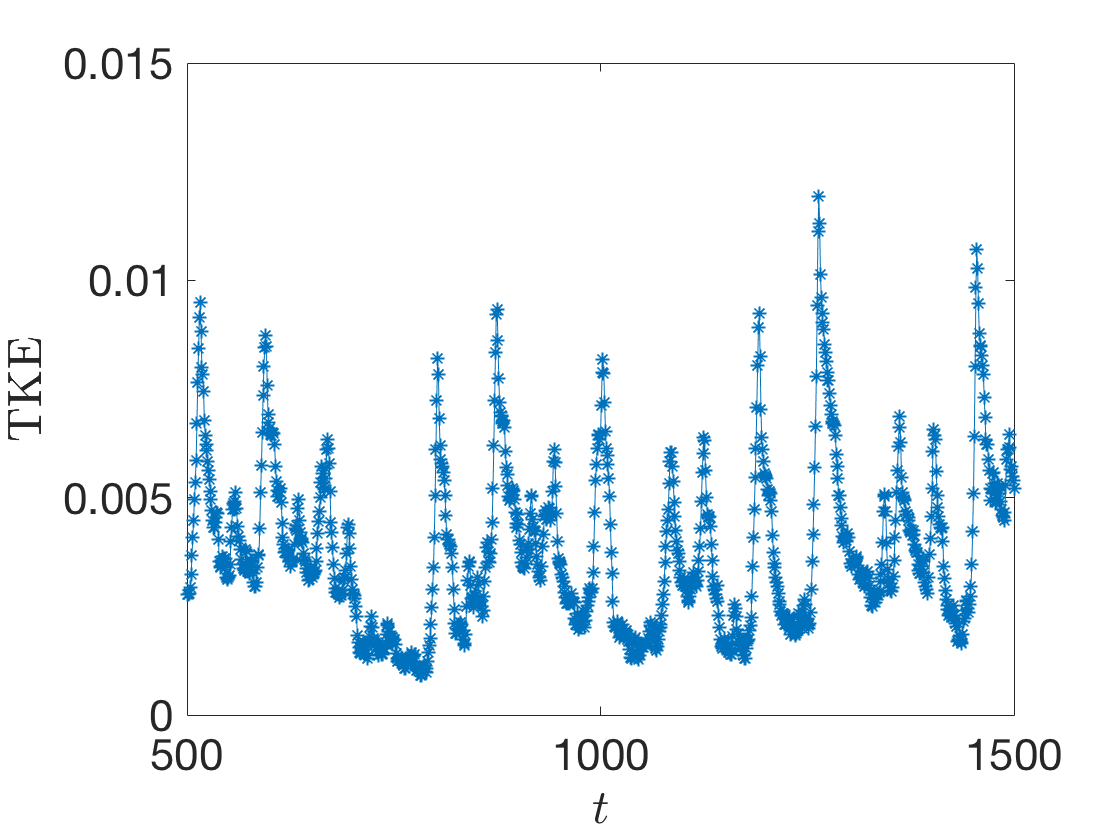}}
 ~~
 \subfloat[${\rm Re}= 25000$]
{  \includegraphics[width=0.30\textwidth]
 {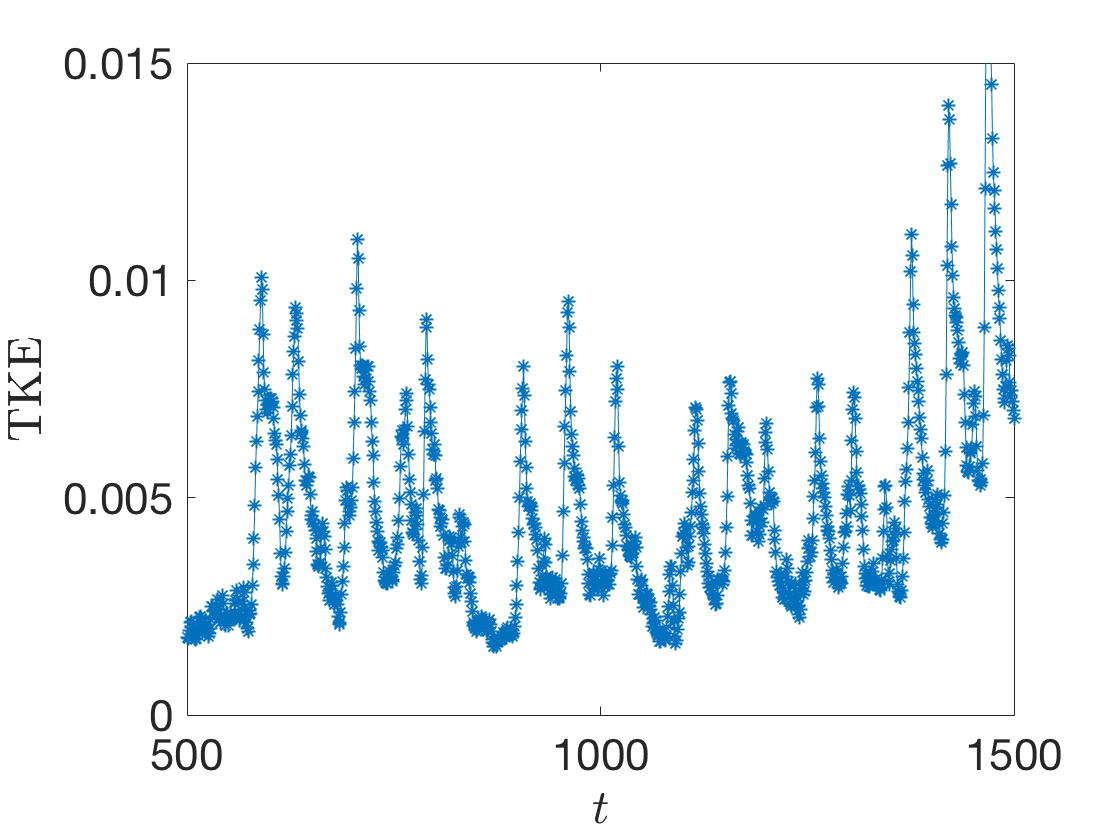}}
    
 \caption{A lid-driven cavity problem.
 Behavior of the turbulent kinetic energy ${\rm TKE}$ with time for three values of the Reynolds number.
  }
 \label{fig:TKE}
  \end{figure}

Figure \ref{fig:velocity_streamlines} shows the velocity streamlines  for three different times for  ${\rm Re}= 15000$, while 
Figure \ref{fig:velocity_streamlines_problem} shows the velocity streamlines for several times for 
${\rm Re}=20000$: for the latter value of the Reynolds number, 
we observe the presence of vortices along the edges  and in the center of the cavity.
 We remark that the same behavior has been observed by Cazemier et al. in \cite{cazemier1998proper} (cf. Figure 3, page 1687).
This rare behavior  makes estimates of long-time averages particularly difficult, especially for the TKE. 
Figure \ref{fig:TKE} shows the behavior of the 
turbulent kinetic energy ${\rm TKE}$ with time for three values of the Reynolds number.
We observe that for sufficiently large values of the Reynolds number we have significant peaks in the TKE. 
These peaks correspond to eddies that are ejected into the core region and cross the cavity. 

\begin{figure}[h!]
\centering
\subfloat[]
{  \includegraphics[width=0.25\textwidth]
 {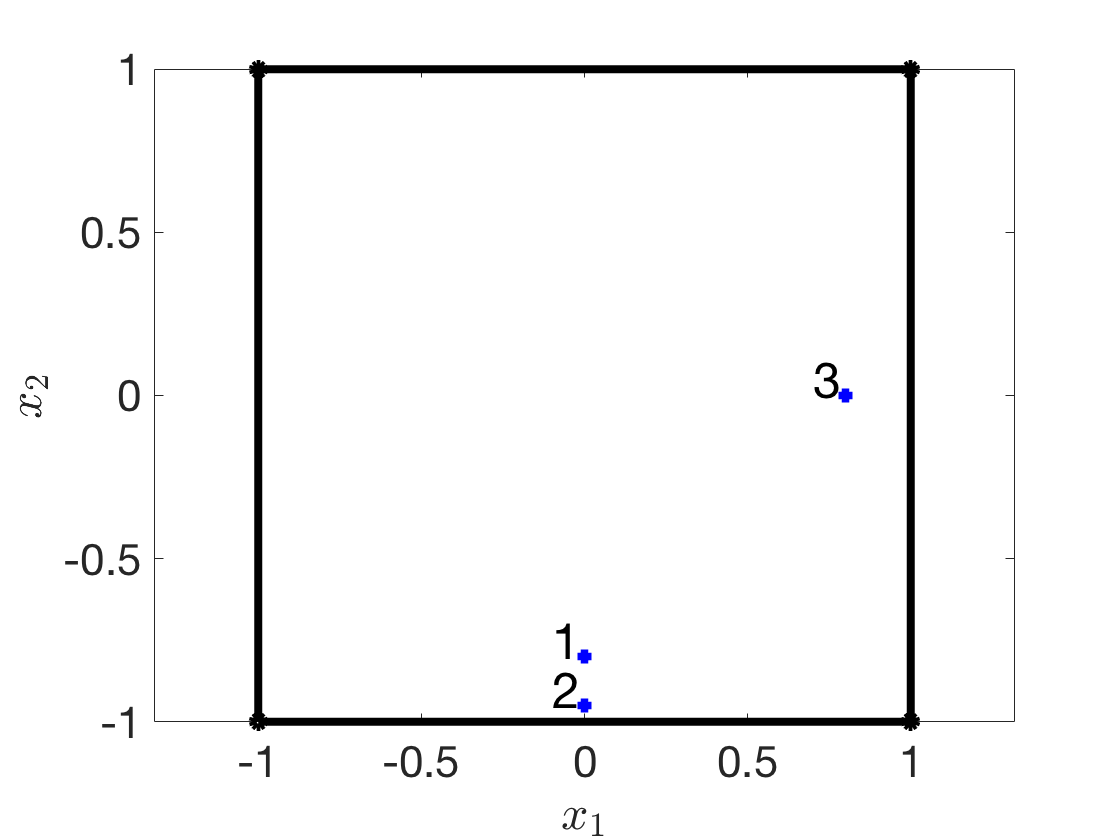}}
 ~~
 \subfloat[${\rm Re}= 15000$, $x_1^{\rm probe}$]
{  \includegraphics[width=0.25\textwidth]
 {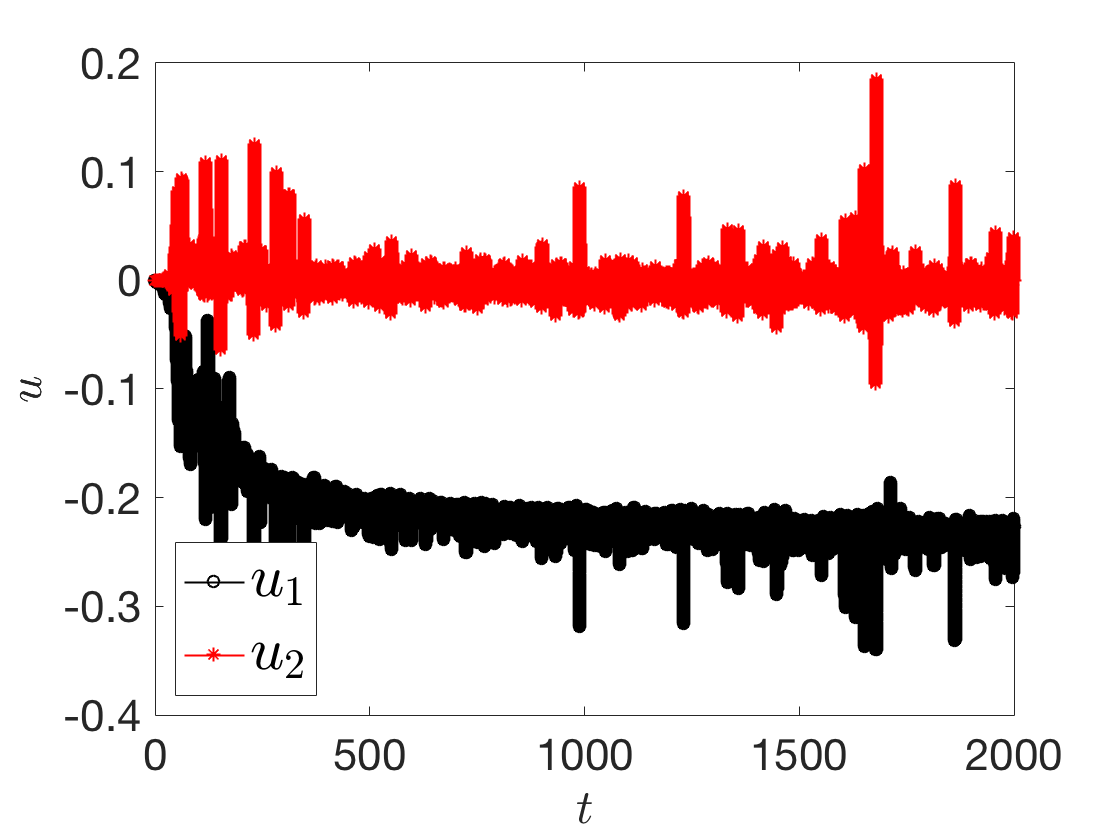}}
 ~~
 \subfloat[${\rm Re}= 15000$, $x_2^{\rm probe}$]
{  \includegraphics[width=0.25\textwidth]
 {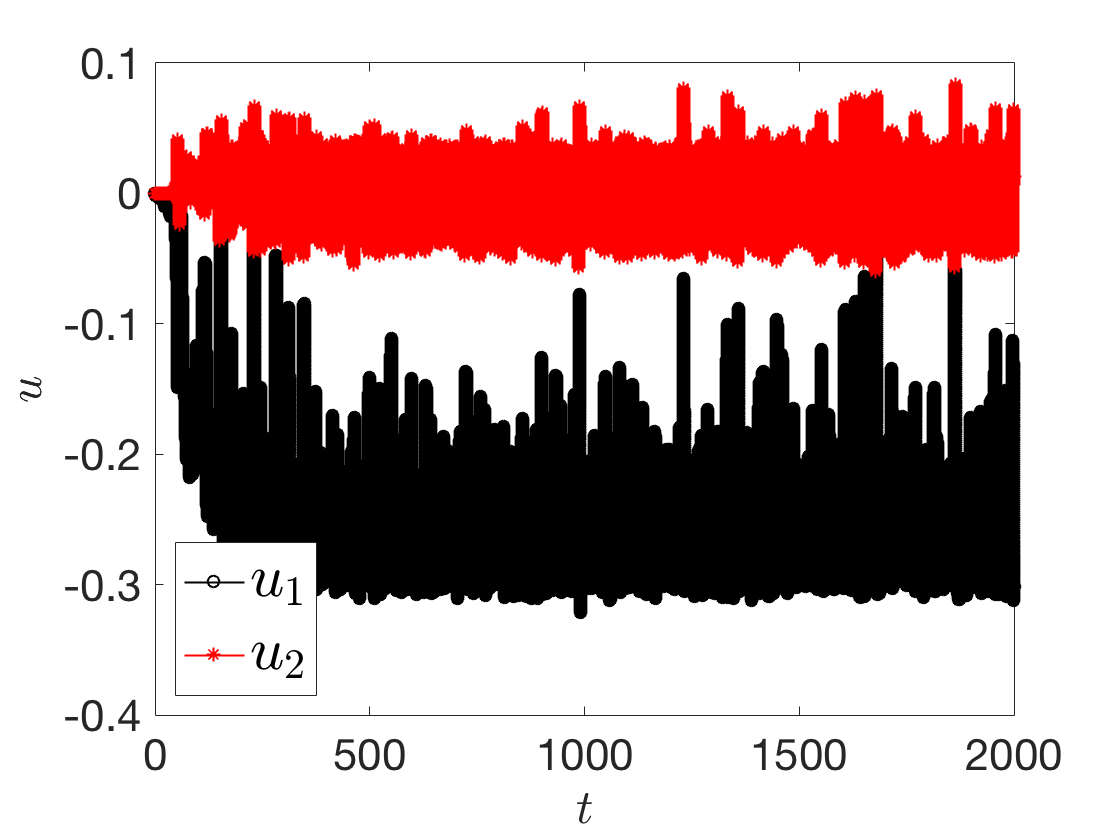}}
     ~~
 \subfloat[${\rm Re}= 15000$, $x_3^{\rm probe}$]
{  \includegraphics[width=0.25\textwidth]
 {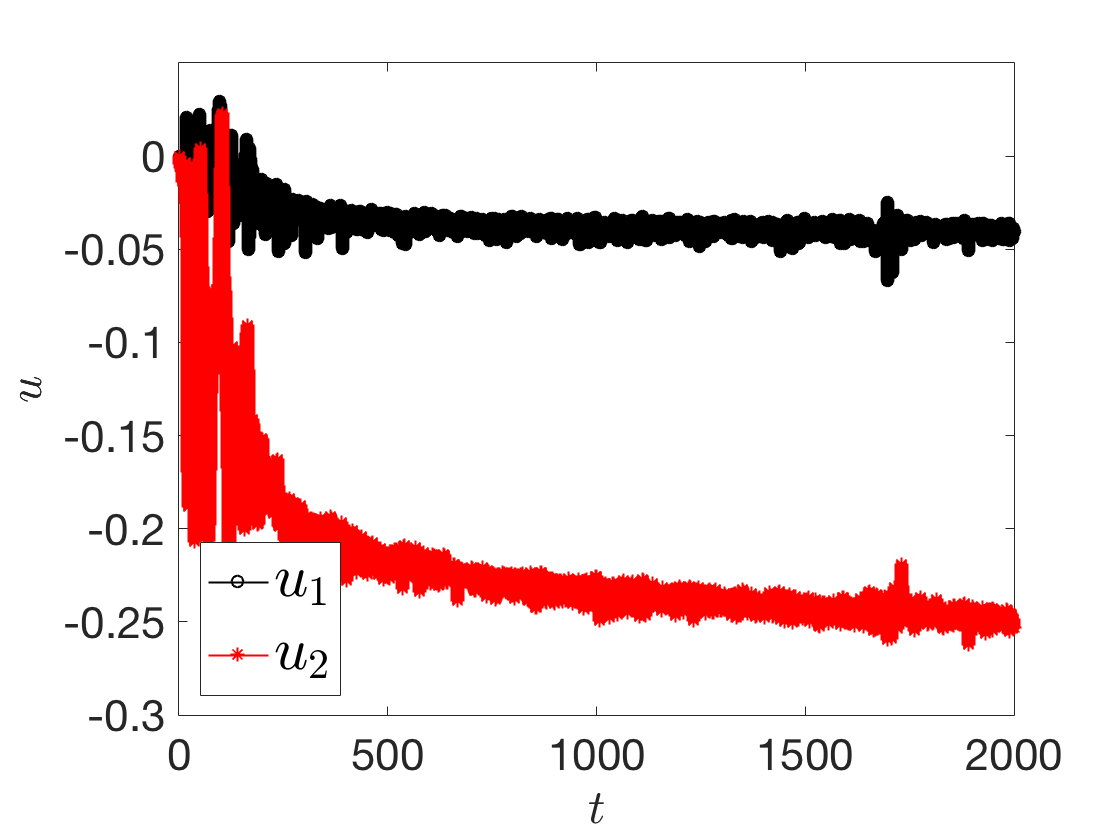}}
 
\subfloat[]
{  \includegraphics[width=0.25\textwidth]
 {imm/vis/time/points.png}}
 ~~
 \subfloat[${\rm Re}= 25000$, $x_1^{\rm probe}$]
{  \includegraphics[width=0.25\textwidth]
 {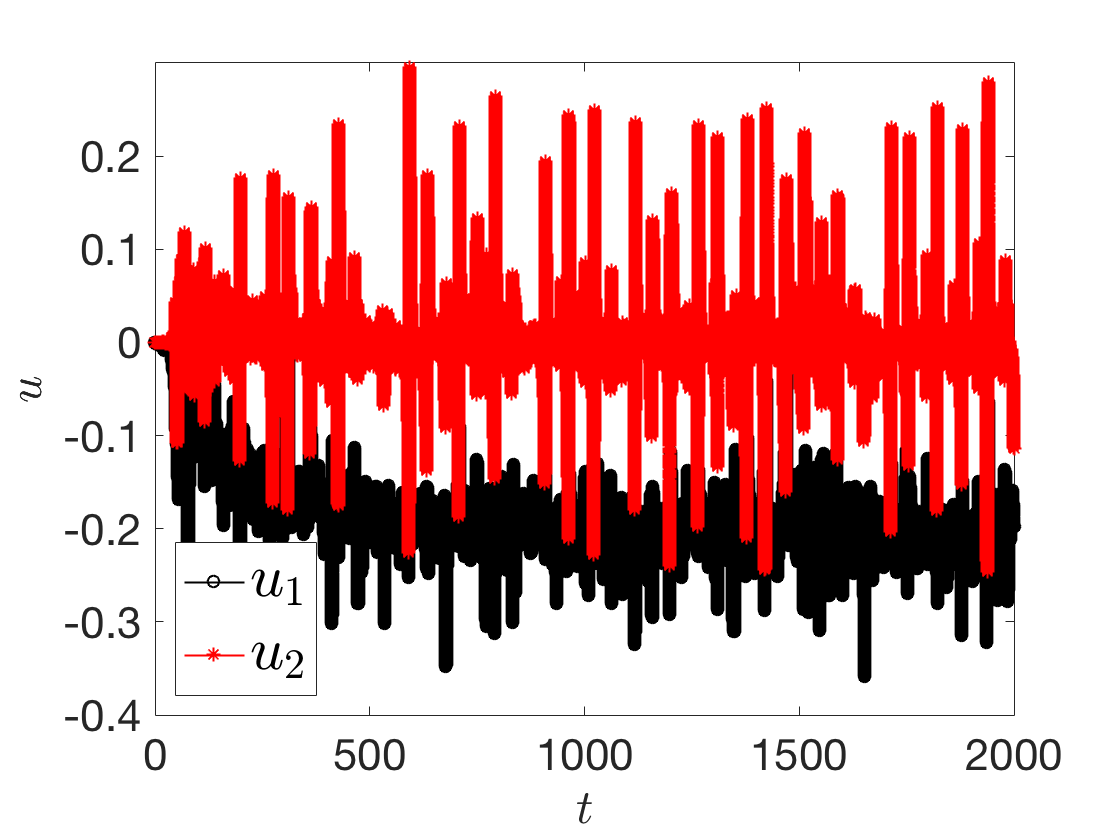}}
 ~~
 \subfloat[${\rm Re}= 25000$, $x_2^{\rm probe}$]
{  \includegraphics[width=0.25\textwidth]
 {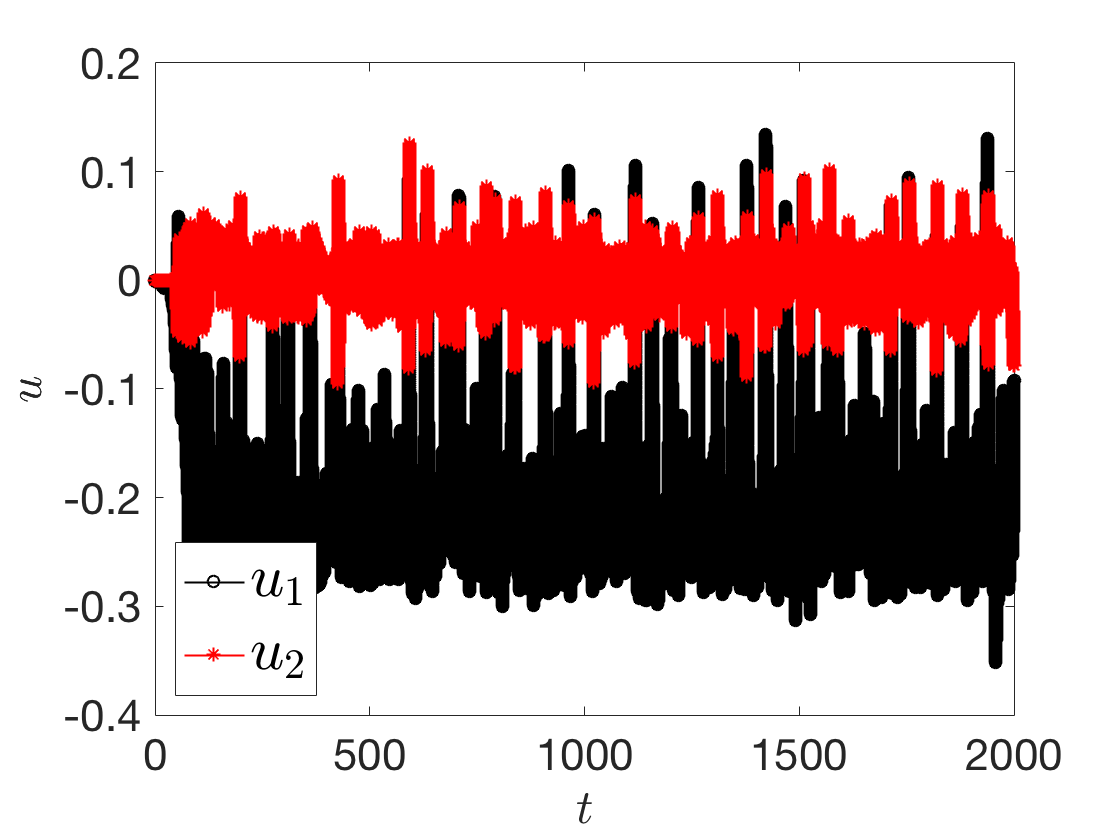}}
     ~~
 \subfloat[${\rm Re}= 25000$, $x_3^{\rm probe}$]
{  \includegraphics[width=0.25\textwidth]
 {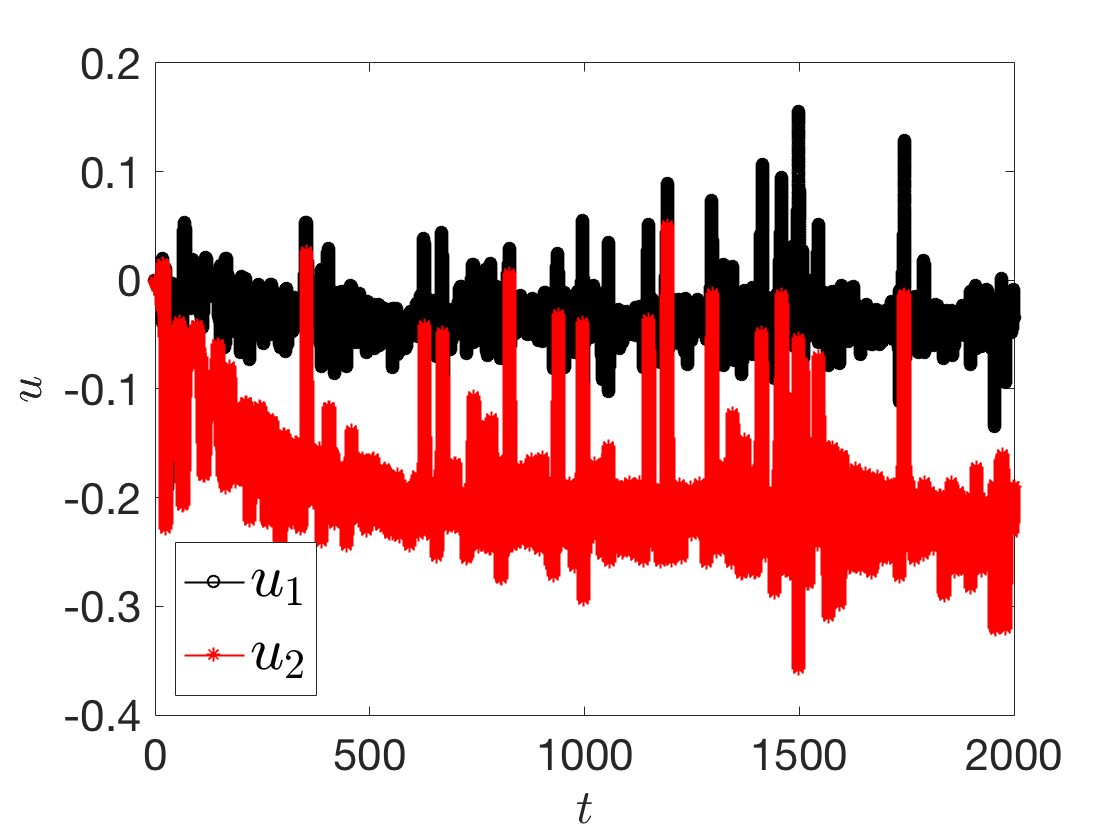}} 
  
 \caption{A lid-driven cavity problem.
 Behavior of the velocity components at three spatial locations, for two values of ${\rm Re}$
 ($x_1^{\rm probe} = [0,-0.8]$, $x_2^{\rm probe} = [0,-0.95]$, $x_3^{\rm probe} = [0.8,0]$).
  }
 \label{fig:time_plot}
  \end{figure}

\begin{figure}[h!]
\centering
\subfloat[${\rm Re}= 15 10^3$, $x_1^{\rm probe}$]
{  \includegraphics[width=0.33\textwidth]
 {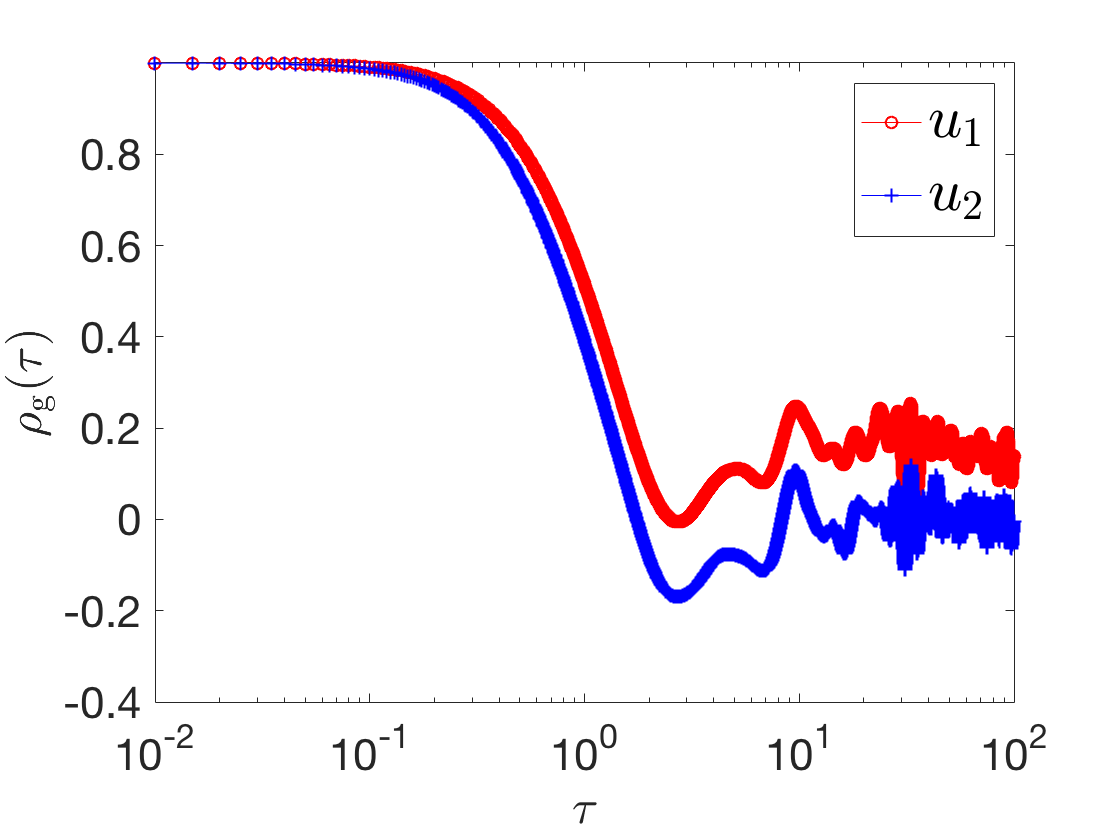}}
 ~~
\subfloat[${\rm Re}= 15 10^3$, $x_2^{\rm probe}$]
{  \includegraphics[width=0.33\textwidth]
 {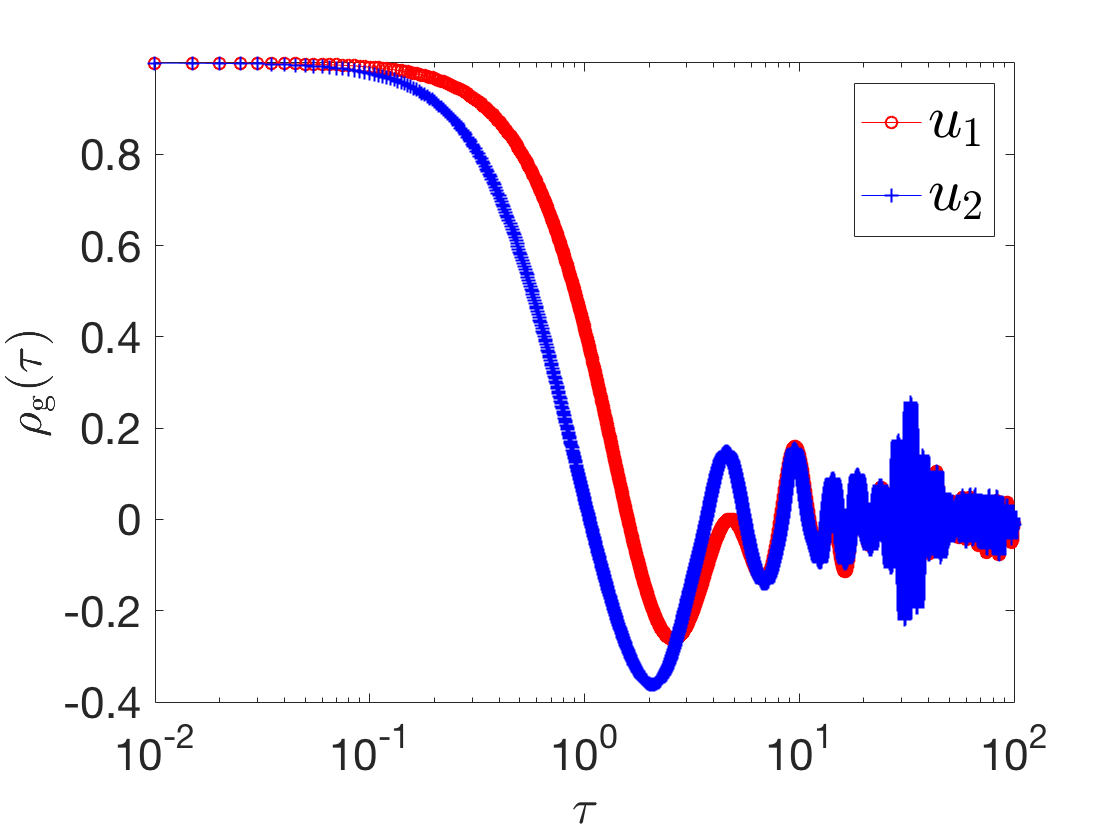}}
    
   \subfloat[${\rm Re}= 25 10^3$, $x_1^{\rm probe}$]
{  \includegraphics[width=0.33\textwidth]
 {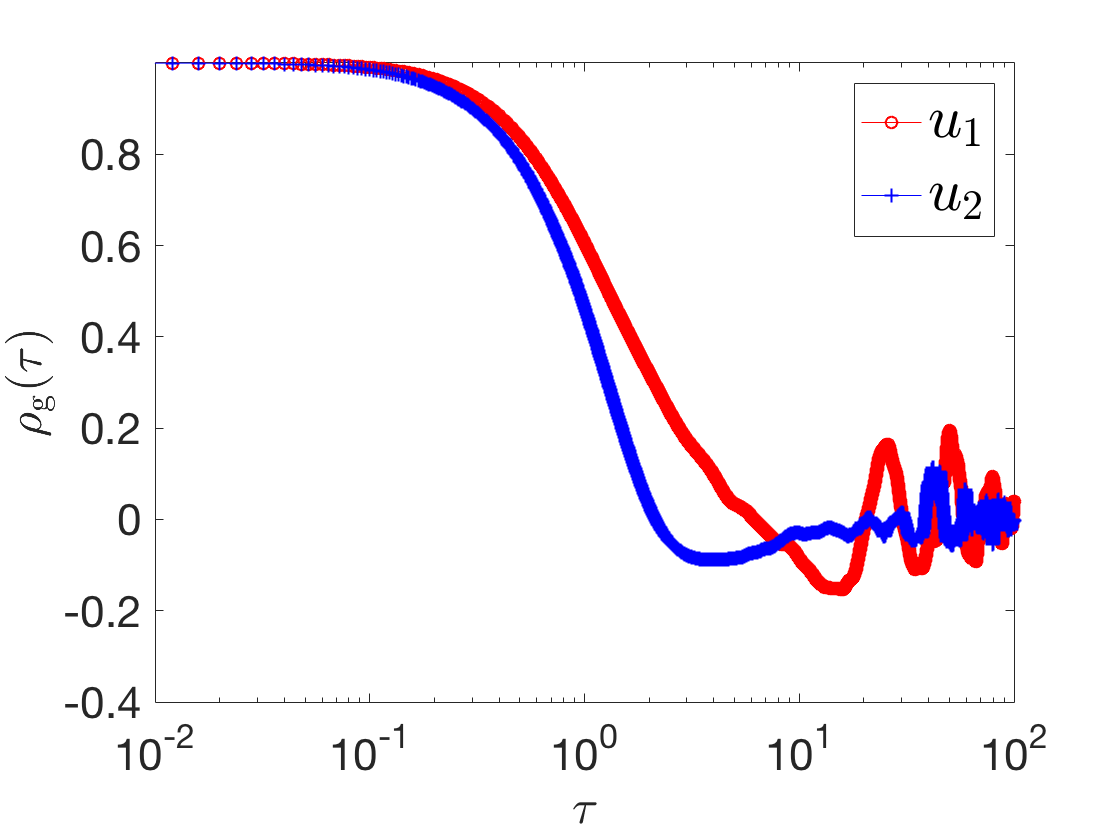}}
 ~~
\subfloat[${\rm Re}= 25 10^3$, $x_2^{\rm probe}$]
{  \includegraphics[width=0.33\textwidth]
 {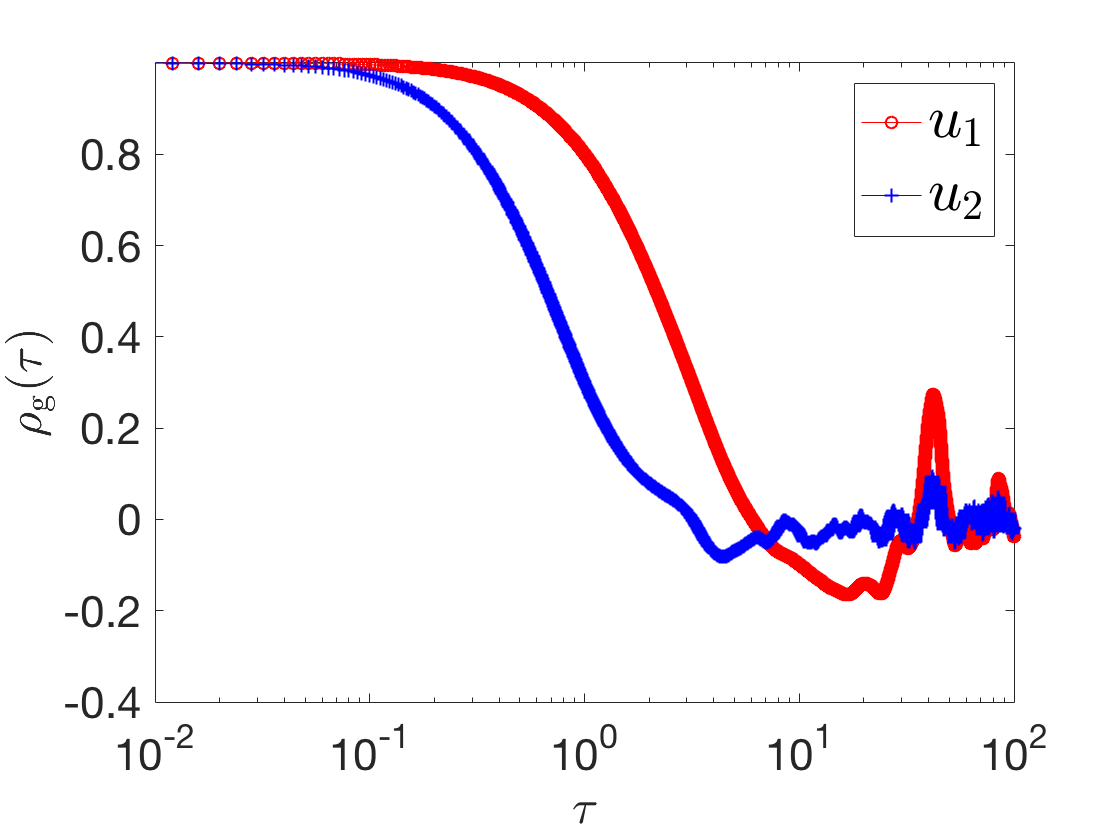}} 
    
\caption{A lid-driven cavity problem.
 Behavior of the autocorrelation for the velocity components at two spatial locations, for two values of ${\rm Re}$
 ($x_1^{\rm probe} = [0,-0.8]$, $x_2^{\rm probe} = [0,-0.95]$).
  }
 \label{fig:time_plot_autocorrelation}
  \end{figure}

Figure \ref{fig:time_plot} shows the behavior with time of the first and second components of the velocity field at three spatial locations for ${\rm Re}=15000$ and ${\rm Re}=25000$.
We observe that for $t \gtrsim T_0 = 500$ the effects of the transient dynamics are negligible.
Figure \ref{fig:time_plot_autocorrelation} shows the autocorrelation factors associated with the time series
$\{ u_i(x_{\ell}^{\rm probe}, t^j ,{\rm Re}) \}_{j=J_0}^J$ for $i,\ell=1,2$
and for $ {\rm Re}=15000, 25000$. 
We here define the autocorrelation factors
for a time sequence $\{ y^j \}_{j=0}^J$
 as follows:
$$
\rho_{\rm g}(\tau = \kappa \Delta t) =
\frac{
\displaystyle{
\frac{1}{J-\kappa - J_0+1}
\sum_{j=J_0}^{J-\kappa}
\left(
y^j - 
\langle y \rangle_{\rm g}
\right)
\left(
y^{j+\kappa} - 
\langle y \rangle_{\rm g}
\right)
}
}{
\displaystyle{
  \langle \left(  y - \langle y \rangle_{\rm g} \right)^2 \rangle_{\rm g}  }}.
$$
We observe that the autocorrelation factor decreases as $\tau$ increases, and is roughly $0.8$ for $\tau=1$,  for all probes considered.

\section{\emph{A posteriori} assessment of the POD accuracy}
\label{sec:cv_POD}
As explained in the main body of the paper, 
POD relies on a snapshot set $\{ \mathring{u}^k \}_{k=1}^K$
to generate a $N$-dimensional approximation  space for the  full trajectory $\{ \mathring{u}^j \}_{j=J_0}^J$ in the limit $J \to \infty$. The snapshot set is associated with the sampling times $\{t_{\rm s}^k := T_0 + \Delta t_{\rm s} k \}_{k=1}^K$, where $T_0 = t_{\rm g}^{J_0}$, $K$  is the cardinality of the snapshot set, and $\Delta t_{\rm s}$ is the sampling period.

The choices of $\Delta t_{\rm s}$ and $K$ 
are a trade-off between (i) information content of the snapshot set, and 
(ii) computational resources. The snapshot set should  be rich enough to accurately estimate the first $N$ POD modes associated with the full trajectory $\{ \mathring{u}^j \}_{j=J_0}^J$ in the limit $J \to \infty$.
On the other hand, it is well-known that POD suffers from  (i) the quadratic growth in $K$ in computational complexity for computing the Gramian, and  for computing the symmetric eigen-decomposition; and (ii) the memory requirements related to the storage of the snapshots, which scale linearly with $K$. Furthermore, by increasing $\Delta t_{\rm s}$ and $K$, we ultimately increase the number of time steps performed by the spectral element solver  --- which is given by $J = \left(T_0 + \Delta t_{\rm s} K  \right)/ \Delta t$.

In this Appendix, we propose a cross-validation 
(CV, see, e.g., \cite{kohavi1995study} and \cite[Chapter 7.10]{hastie2009elements})
strategy to estimate the $\ell^2$-averaged projection error  associated with the POD reduced space over the full-trajectory, 
$$
\mathcal{E}(\{ \mathring{u}^j \}_{j=J_0}^J, \mathcal{Z}^{\rm u})
=
\frac{1}{J + 1 -J_0}
\,
\sum_{j=J_0}^J 
\,
\|  \mathring{u}^j  -  \Pi_{  \mathcal{Z}^{\rm u}  }^V\mathring{u}^j     \|_V^2.
$$
Estimates of this quantity might be employed to decide whether or not to  acquire new snapshots and/or to increase the dimension $N$ of the reduced space.
On the other hand, evaluations of the autocorrelation factor introduced in  Appendix \ref{sec:lid_driven_appendix} can be used to assess \emph{a posteriori} the amount of redundancy in the  snapshot set.

 Since in our setting the snapshots are correlated in time 
(cf. Appendix \ref{sec:lid_driven_appendix}), 
we here rely on the $h$-block variant proposed in \cite{burman1994cross} (see also \cite{racine2000consistent}). The approach relies on the assumption that the snapshot set is associated with a stationary process: under this assumption, the covariance matrix between $\mathring{u}^j$ and $\mathring{u}^{j+\kappa}$ is only a function of $\kappa$, and approaches $0$ as $\kappa \to \infty$. 
The key idea of $h$-block CV is to reduce the training set by removing the $h$ observations preceding and following the observation in the test set.
In  section \ref{sec:method_cv}, we adapt the computational procedure discussed in \cite{burman1994cross} to the particular learning task of interest; then, in  section  \ref{sec:numerics_cv}, we apply the procedure to the case ${\rm Re}= 15000$ to support our choice $K=500$.

Before proceeding with the presentation of the methodology,
we remark that, in the statistics literature, several authors have proposed validation techniques to assess the accuracy of POD (or, equivalently, PCA and 
Karhunen-Lo{\'e}ve) spaces.
We refer to \cite[Chapter 6]{jolliffe2002principal} and to the references therein for a number of different proposals. We further recall the work by Chowdhary and Najm \cite{chowdhary2016bayesian} that relies on a Bayesian framework to account for inaccuracies due to limited sample size. The approaches  presented in 
 \cite{jolliffe2002principal,chowdhary2016bayesian} aim at generating confidence (credible) regions for the estimate of the POD modes; on the other hand,  we are here primarily interested in assessing the 
\emph{out-of-sample} accuracy of the $N$-dimensional POD reduced space $\mathcal{Z}^{\rm u}$ for the full trajectory 
$\{ \mathring{u}^j \}_{j=J_0}^J$. 
For completeness, 
we also  recall that several authors (\cite{paul2015adaptive,himpe2016hapod})
have proposed and analyzed hierarchical POD approaches to reduce the size $K$ of the snapshot set, without significantly compromising the accuracy of the POD space.

\subsection{$h$-block Cross-Validation}
\label{sec:method_cv}

Algorithm \ref{cvPOD} summarizes  the computational procedure for the estimation of the $\ell^2$-averaged projection error 
$\mathcal{E}( \{ \mathring{u}^j \}_{j=J_0}^J, \mathcal{Z}^{\rm u} )$.
We   observe that for $h=0$ the procedure reduces to Leave-One-Out-Cross-Validation (LOOCV).
We further observe that the approach requires the assembling of the Gramian matrix $\mathbb{U}$, and then the solution  to $K$ dense eigenvalue problems of size $K-2h-1$: for the particular problem at hand, the computational cost associated with the procedure is negligible compared to 
the computational cost associated with
the solution to the FOM.
Finally, we emphasize that the procedure relies on the input parameter $h$.
We here propose to choose $h$ based on the analysis of the autocorrelation factor:
recalling the results presented in Appendix    \ref{sec:lid_driven_appendix}, we consider $h=\frac{4}{\Delta t_s} = 4$.

\begin{algorithm}[H]
 \caption{$h$-block Cross-Validation}
  \label{cvPOD}
 
$[  \widehat{\mathcal{E}} ] = $ \texttt{hblock}-CV $(\{  \mathring{u}^k  \}_{k=1}^K, h, N )$
 
 \small
 \emph{Inputs:}  $\{  \mathring{u}^k  \}_{k=1}^K$ = snapshot set,  $h$ = correlation parameter, $N$ = size of the POD space.
 
\emph{Output:} $\widehat{\mathcal{E}}=$  CV estimate of  $\mathcal{E} ( \{ \mathring{u} \}_{j=J_0}^J, \mathcal{Z}^{\rm u} )$.

\normalsize 
 
  \begin{algorithmic}[1]
\For{$k=1,\ldots,K$ }

\State
$[ \mathcal{Z}^{\rm u,(k)} := {\rm span}  \{  \zeta_n^{(k)} \}_{n=1}^{N} ]
=
\texttt{POD}_V
\left(
\{  \mathring{u}^1, \ldots, \mathring{u}^{k-h-1}, \mathring{u}^{k +h+1}, \ldots, \mathring{u}^K   \}_{k=1}^K,
N
\right)$

\EndFor

\State
Compute the CV estimate $\widehat{\mathcal{E}}$ as 
$
\widehat{\mathcal{E}} = \frac{1}{K}
\sum_{k=1}^K
\,
\| \mathring{u}^k -  \Pi_{ \mathcal{Z}^{\rm u,(k)}      }^V  \mathring{u}^k   \|_V^2.$
  \end{algorithmic}
 \end{algorithm}

\subsection{Results for ${\rm Re}=15000$}
\label{sec:numerics_cv}

\begin{figure}[h!]
\centering
{  \includegraphics[width=0.33\textwidth]
 {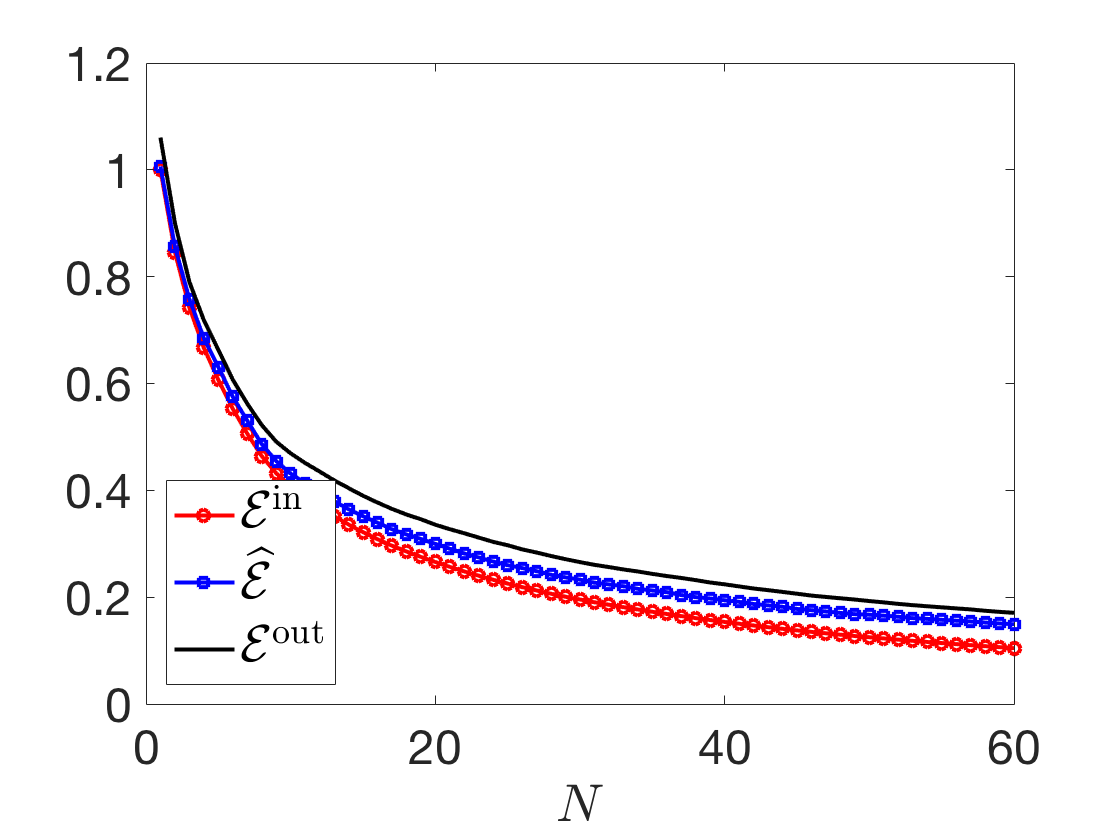}}
  
 \caption{A Cross-Validation procedure for   the \emph{a posteriori} assessment of the POD accuracy.
 Behavior of 
 $\mathcal{E}^{\rm in}$, $\widehat{\mathcal{E}}$, and $\mathcal{E}^{\rm out}$
 with $N$.
 All quantities are  normalized by  $\mathcal{E}^{\rm in}(N=1)$
 ($K=500$, $h=2$, ${\rm Re}=15000$). 
  }
 \label{fig:hblock_cv}
  \end{figure}

Figure \ref{fig:hblock_cv} shows the behavior of  $\widehat{\mathcal{E}}$ for different values of $N$ for the snapshot set $\{ \mathring{u}^k \}_{k=1}^K$ associated with $\{t_{\rm s}^k = 500 + k  \}_{k=1}^{K=500}$, and ${\rm Re}=15000$. We compare results with the in-sample estimate
$$
\mathcal{E}^{\rm in}
=
\frac{1}{K} \sum_{k=1}^K \, \|  \mathring{u}^k  - \Pi_{\mathcal{Z}^{\rm u}}^V \mathring{u}^k     \|_V^2,
$$
and the out-of-sample estimate
$$
\mathcal{E}^{\rm out}
=
\frac{1}{K} \sum_{k=K+1}^{2K} \, \|  \mathring{u}^k  - \Pi_{\mathcal{Z}^{\rm u}}^V \mathring{u}^k     \|_V^2,
$$
where $\{  \mathring{u}^k\}_{k=K+1}^{2K}$ 
are associated with the sampling times $\{t_{\rm s}^k = 500 + k  \}_{k=K+1}^{2K}$.
For visualization purposes, we normalize  $\widehat{\mathcal{E}}$,  $\mathcal{E}^{\rm in}$, and
$\mathcal{E}^{\rm out}$
 by $\mathcal{E}^{\rm in}(N=1)$:  for $N=60$, $\widehat{\mathcal{E}} \approx  15\%  \times \mathcal{E}^{\rm in}(N=1)$,
 $\mathcal{E}^{\rm in} \approx  10\%  \times \mathcal{E}^{\rm in}(N=1)$, and
  $\mathcal{E}^{\rm out} \approx  17 \%  \times \mathcal{E}^{\rm in}(N=1)$. 
  We  observe that   $\widehat{\mathcal{E}} $ is a more accurate estimate 
of   $\mathcal{E}^{\rm out} $ compared to the in-sample estimate
$\mathcal{E}^{\rm in}$. 

\section{On the definition of ROM stability}
\label{sec:ROM_stability}
Based on the results of section \ref{sec:reproduction_problem},
 we could take a pragmatic view of the long-time ROM stability.
Given the reduced space $\mathcal{Z}^{\rm u} \subset V_{\rm div}$, we define the best-fit errors associated with mean flow and TKE:
$$
e_1^{\rm opt}
:=
\frac{ \|   \langle \mathring{u} \rangle_{\rm g}  - \Pi_{\mathcal{Z}^{\rm u}}^{V}  \langle \mathring{u} \rangle_{\rm g}  \|_V }{
\| \langle u  \rangle_{\rm g}        \|_V},
\quad
e_2^{\rm opt}
:=
\frac{   |     \langle  {\rm TKE}  \rangle_{\rm s}     -   \langle  {\rm TKE}^{\rm opt}( \cdot; \mathcal{Z}^{\rm u}  )  \rangle_{\rm s}    |    }{ \langle 
 {\rm TKE}
  \rangle_{\rm s} },
$$
where $\Pi_{\mathcal{Z}^{\rm u}}^{V} $ denotes the projection operator associated with the  $V$ inner product on the subspace $\mathcal{Z}^{\rm u}$, and 
$ {\rm TKE}^{\rm opt}( t; \mathcal{Z}^{\rm u}  ) = \frac{1}{2} \int_{\Omega} \| 
\Pi_{\mathcal{Z}^{\rm u}}^{L^2} ( \mathring{u}(t)  -   \langle \mathring{u} \rangle_{\rm g}    )  \|_2^2 \, dx $. Then, we  define the effective stability constants as the ratios between the optimal mean error  and the actual error:
\begin{equation}
\label{eq:stability_constants}
{\rm m}( \mathcal{Z}^{\rm u} )
:=
\frac{  \|  \langle \mathring{u} \rangle_{\rm g}  -  \langle \hat{u} \rangle_{\rm g}   \|_V  }{    
\|  \langle u \rangle_{\rm g}     \|_V
e_1^{\rm opt}},
\qquad
{\rm \sigma}( \mathcal{Z}^{\rm u} )
:=
\frac{ |  \langle  {\rm TKE} \rangle_{\rm s}  -    \langle \widehat{ {\rm TKE}} \rangle_{\rm s}  |   }{   \langle  {\rm TKE} \rangle_{\rm s}   e_2^{\rm opt}}.
\end{equation}
The stability constants ${\rm m}$ and ${\rm \sigma}$ can be used to quantitatively measure the stability of the ROM.
We observe that, 
in the limit $T \to \infty$,
our definition of long-time stability for ROMs is independent of transient dynamics.
We further observe that 
a ROM of dimension $N$ is stable  if and only if mean and variance of the time coefficients
$\{ a_n^j \}_j$ are correctly estimated for $n=1,\ldots,N$.
We finally remark that our definition of stability
 is close to the one proposed in 
\cite{balajewicz2012stabilization}; however, while the definition in \cite{balajewicz2012stabilization} is tailored to $L^2$ POD spaces and $R_g= \langle u \rangle_{\rm g}$, our definition applies to any reduced space and to any choice of the lift.

\begin{figure}[h!]
\centering
\subfloat[]  
{  \includegraphics[width=0.30\textwidth]
 {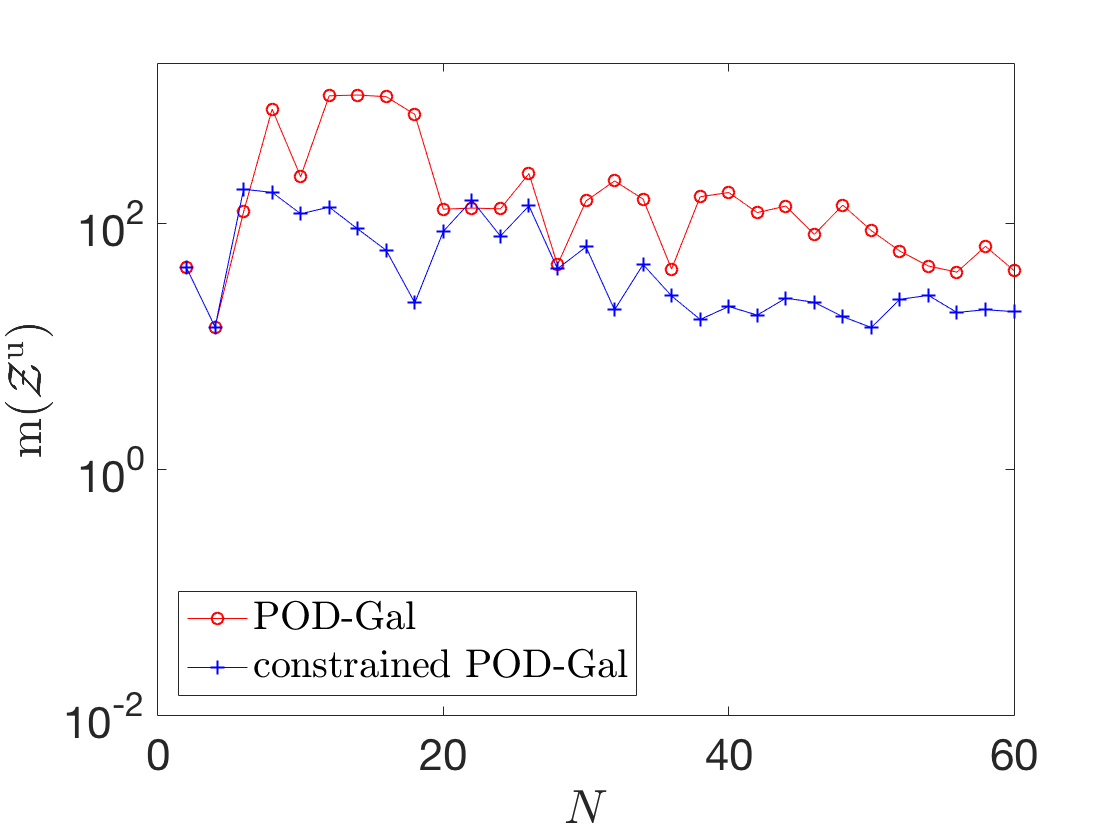}}
 ~~
  \subfloat[ ]  
{  \includegraphics[width=0.30\textwidth]
 {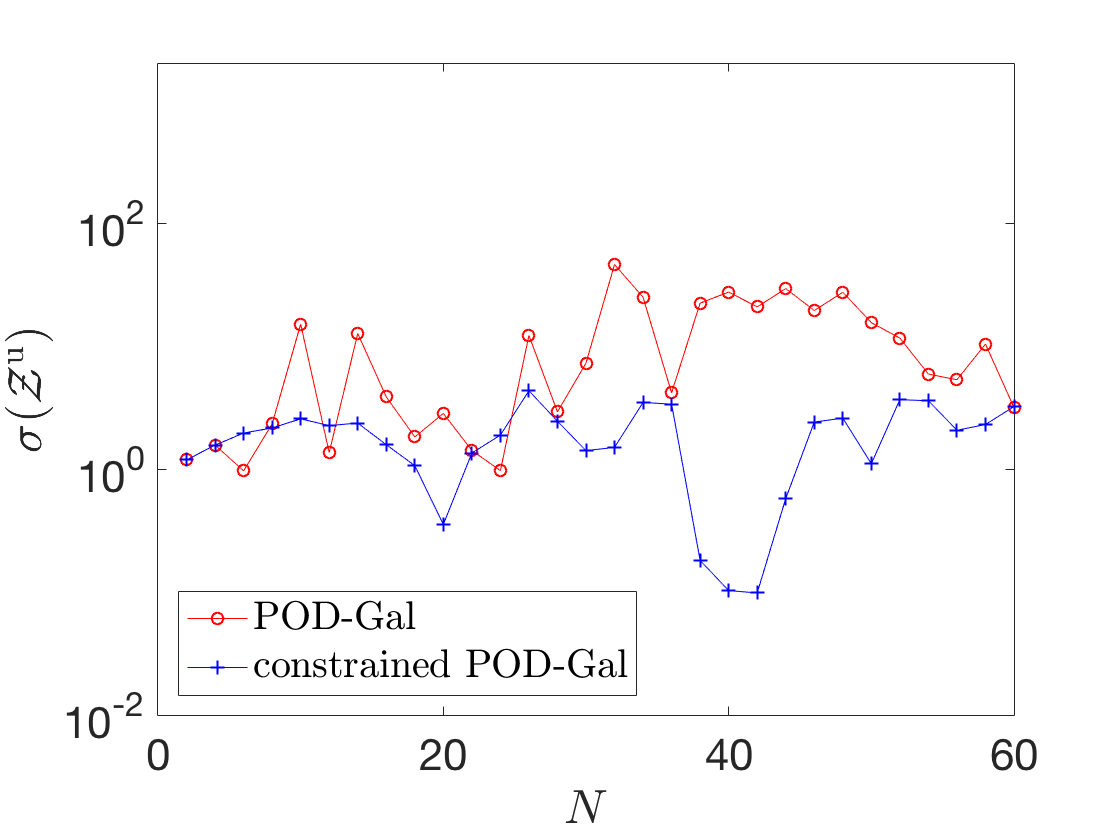}} 
  
 \caption{The solution reproduction problem; 
behavior of  ${\rm m}( \mathcal{Z}^{\rm u} )$ and 
 ${\rm \sigma}( \mathcal{Z}^{\rm u} )$
\eqref{eq:stability_constants}  
  for 
POD-Galerkin and   
constrained POD-Galerkin
 (${\rm Re} = 15000$, $\epsilon=0.01$).
  }
 \label{fig:stab_constants_constrained}
  \end{figure}

Figure \ref{fig:stab_constants_constrained} shows the behavior of  ${\rm m}( \mathcal{Z}^{\rm u} )$ and 
 ${\rm \sigma}( \mathcal{Z}^{\rm u} )$ defined in \eqref{eq:stability_constants}    
  for POD-Galerkin and   constrained POD-Galerkin for ${\rm Re } = 15000$:
 our constrained POD-Galerkin ROM is more stable --- according to the definition given in this Appendix ---
 than the standard Galerkin ROM.
 
\section{Robustness of the constrained formulation}
\label{sec:constrained_stability}
In this Appendix, we present a number of numerical results that provide further insights about the constrained formulation proposed in this paper.
In greater detail, we study the activation rate of the box constraints,   the dependence of the solution to the choice of $\epsilon$,  and   the behavior of  $m_n^{\rm FOM,u}$ and $M_n^{\rm FOM,u}$  defined in \eqref{eq:golden_choice_parameters} 
with respect to the Reynolds number ${\rm Re}$.

\begin{figure}[h!]
\centering
 \subfloat[$N=40$]
{  \includegraphics[width=0.30\textwidth]
 {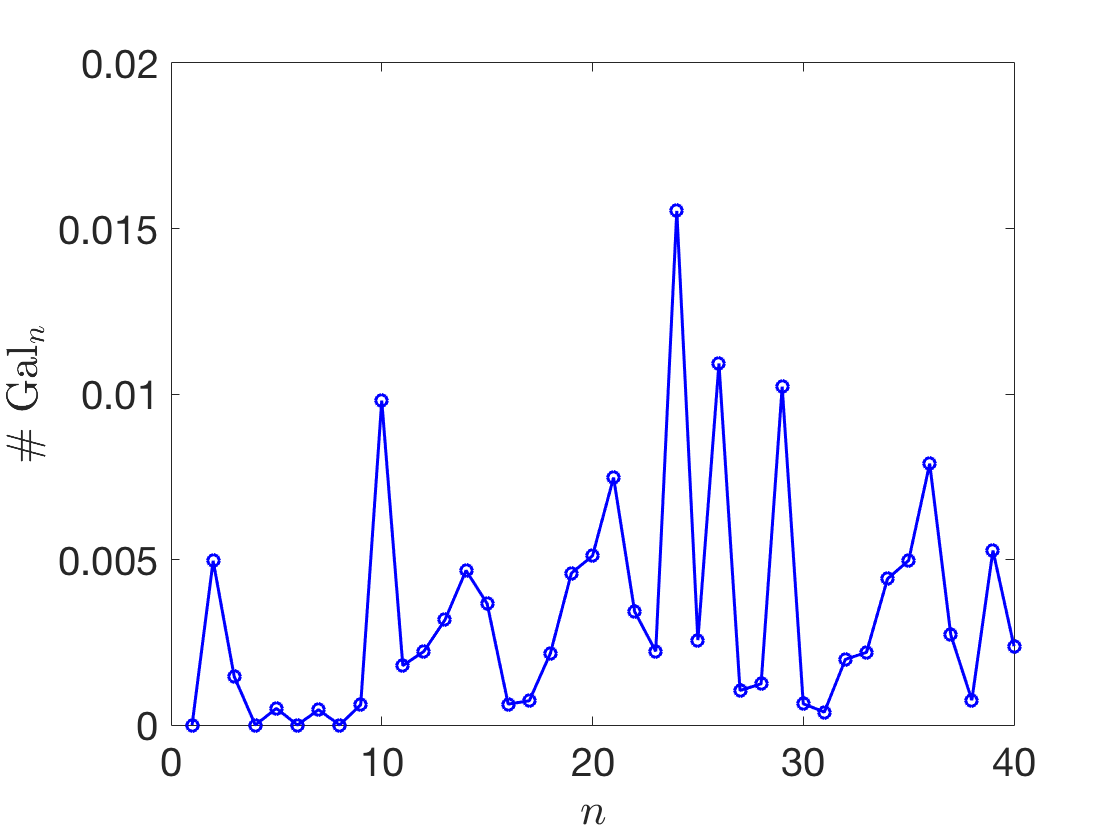}}
   ~~
 \subfloat[$N=60$]
{  \includegraphics[width=0.30\textwidth]
 {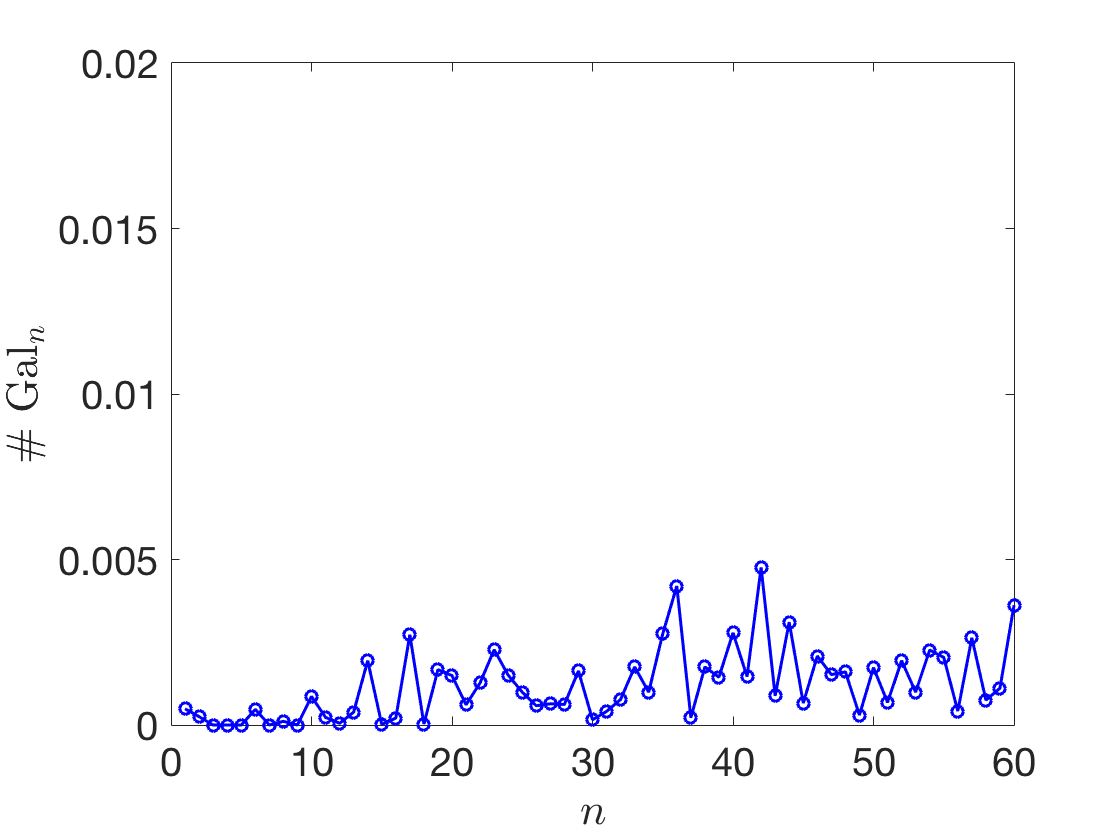}}
 
 \caption{The solution reproduction problem; activity of the box constraints for two values of $N$.
 (${\rm Re} = 15000$, $\epsilon=0.01$).
  }
 \label{fig:cgal_box_constraints}
  \end{figure}

\begin{figure}[h!]
\centering
\subfloat[ $N=40$]  
{  \includegraphics[width=0.30\textwidth]
 {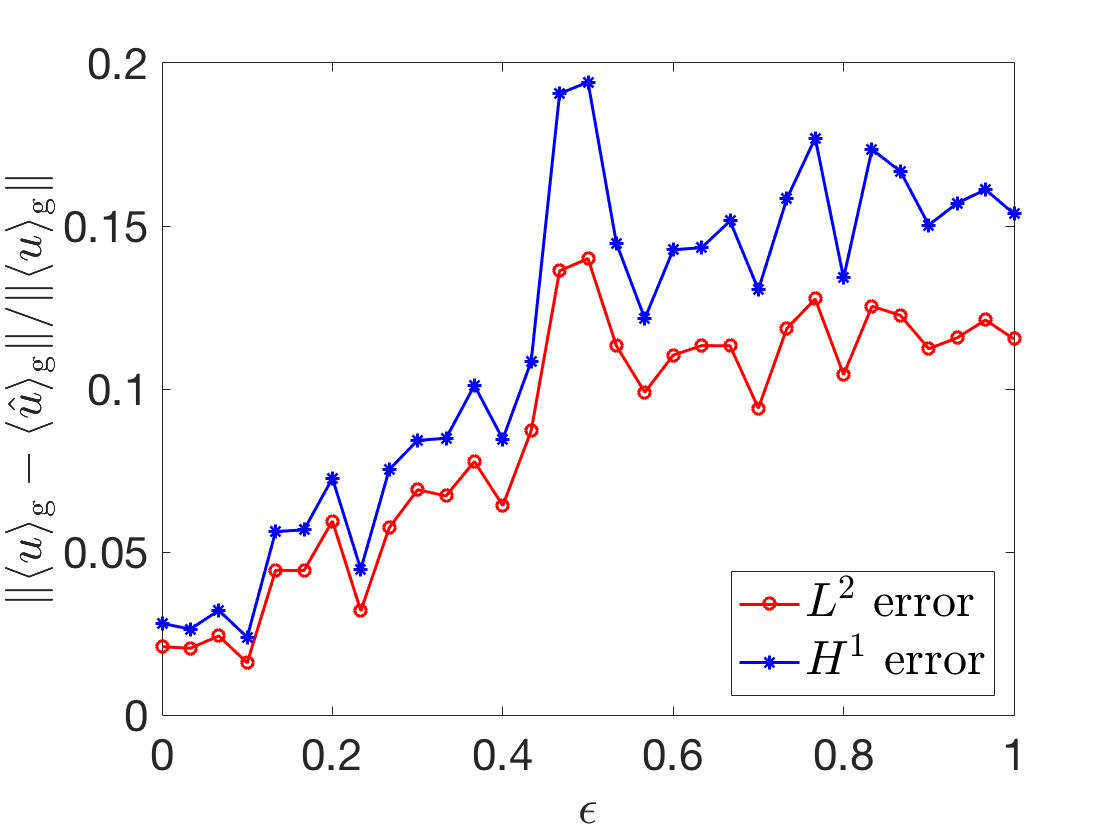}}
 ~~
  \subfloat[ $N=40$]  
{  \includegraphics[width=0.30\textwidth]
 {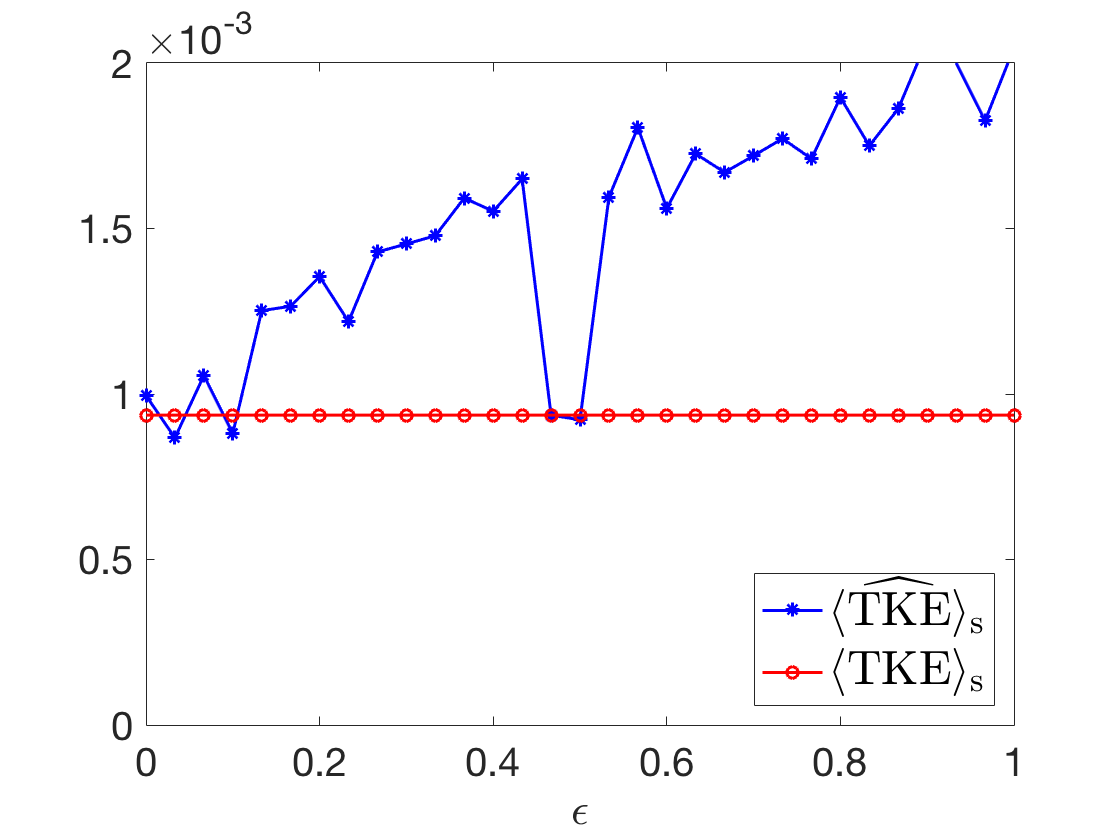}} 
  ~~
  \subfloat[ $N=40$]  
{  \includegraphics[width=0.30\textwidth]
 {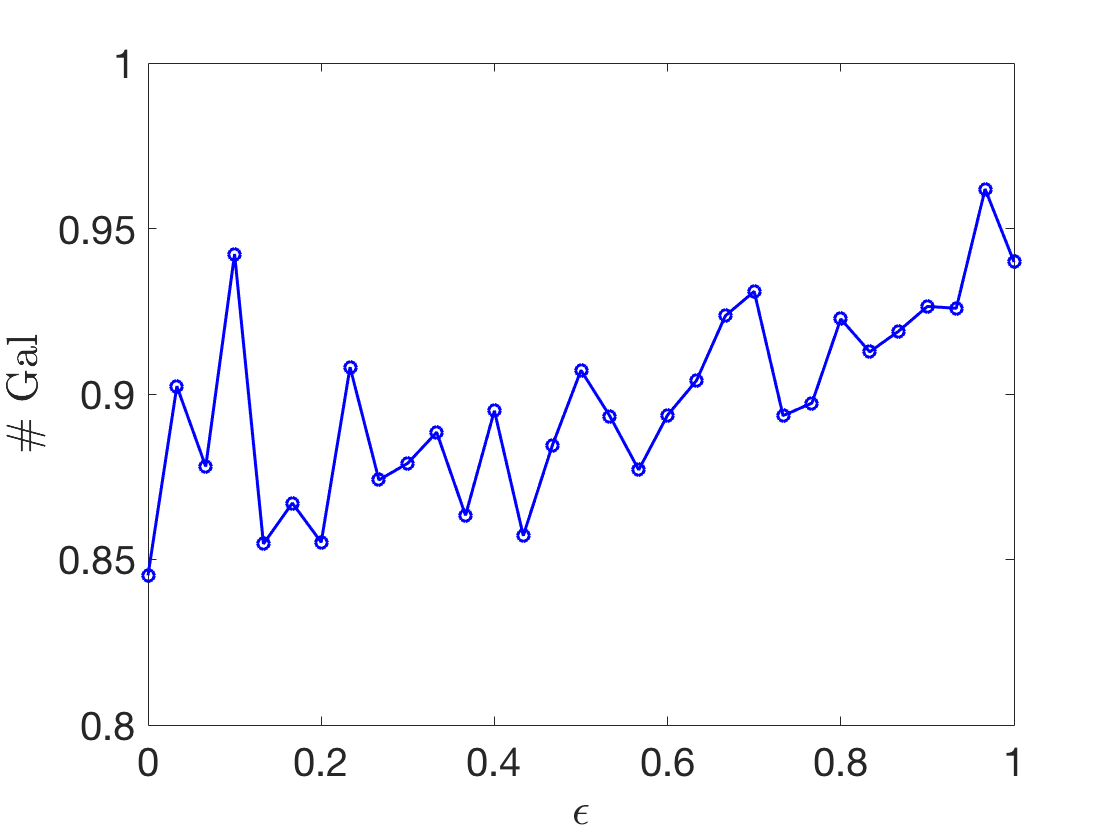}} 
 
\subfloat[ $N=60$]  
{  \includegraphics[width=0.30\textwidth]
 {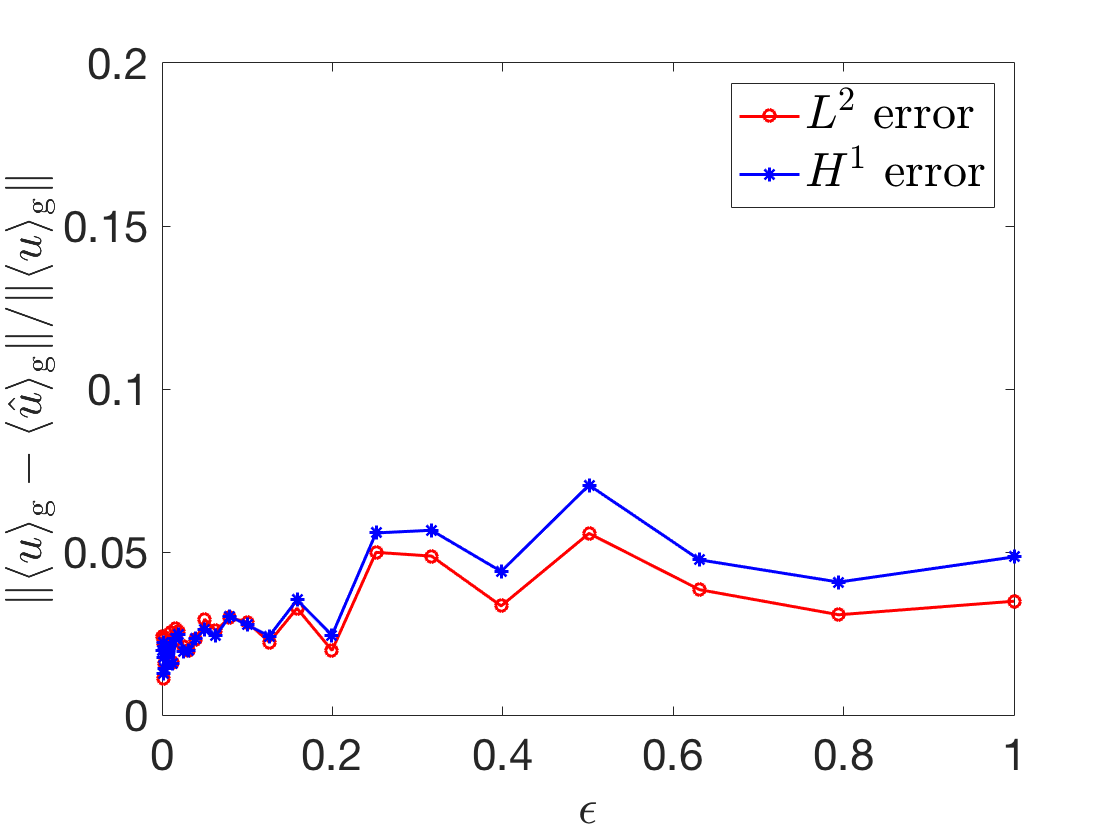}}
 ~~
  \subfloat[ $N=60$]  
{  \includegraphics[width=0.30\textwidth]
 {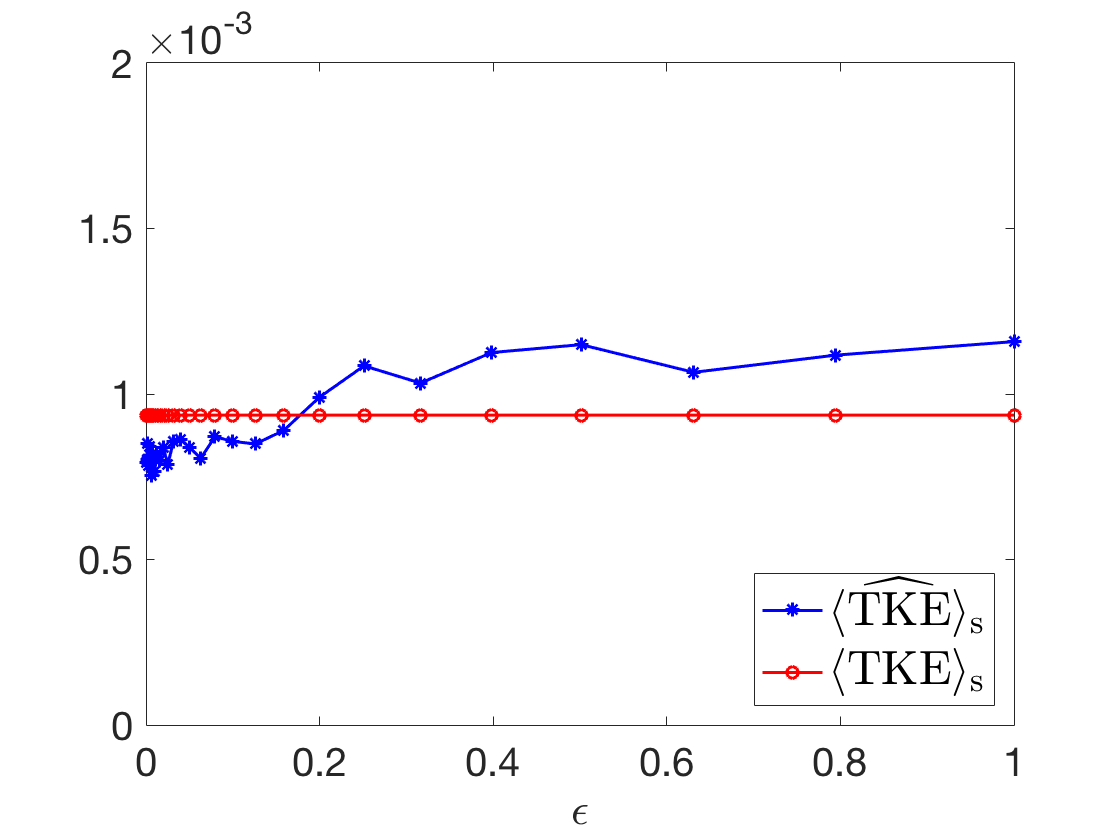}} 
  ~~
  \subfloat[ $N=60$]  
{  \includegraphics[width=0.30\textwidth]
 {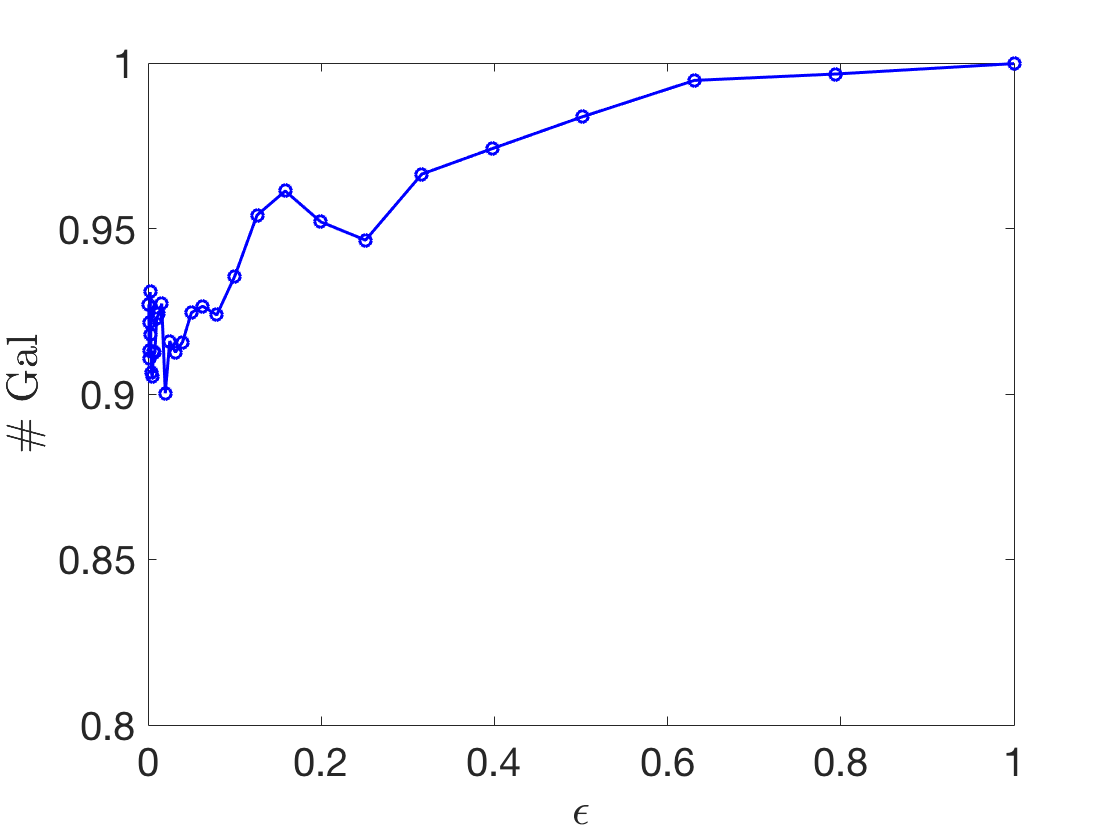}} 
 
 \caption{The solution reproduction problem; sensitivity with respect to  $\epsilon$ for  constrained POD-Galerkin for two values of $N$.
 Figures (a) and (d): behavior of the relative $L^2$ and $H^1$ errors.
 Figures (b) and (e): behavior of the mean TKE.
 Figures (c) and (f): percentage of pure Galerkin solves.
 (${\rm Re} = 15000$).
  }
 \label{fig:cgal_numerics_eps}
  \end{figure}

In Figure \ref{fig:cgal_box_constraints}, we study the behavior of the activation rate of each box constraint for two values of $N$ for ${\rm Re} = 15000$, $\{ t_{\rm s}^k = 500 + k  \}_{k=1}^{K=500}$.
In greater detail, we count how many times the $n$-th component of the solution to Galerkin satisfies the prescribed constraints:
$$
\mbox{\#}
{\rm Gal}_n
:=
\frac{1}{J-J_0}
\sum_{j=J_0+1}^J
\, 
\mathbbm{1} \left(
(\mathbf{a}_{\rm Gal}^j)_n 
\in [\alpha_n, \beta_n]
\right),
\quad
n=1,\ldots,N.
$$
We observe that the behavior with $n$ of $\mbox{\#} {\rm Gal}_n$ is irregular, and strongly depends on $N$.
This suggests that selecting  \emph{a priori} the active constraints might be impractical.
In Figure \ref{fig:cgal_numerics_eps}, we study the behavior of the relative error in the mean flow prediction, the behavior of the mean TKE, and the behavior of the activation rate of the box constraints
$$
\mbox{\#}
{\rm Gal} 
:=
\frac{1}{J-J_0}
\sum_{j=J_0+1}^J
\, 
\mathbbm{1} \left(
(\mathbf{a}_{\rm Gal}^j)_n 
\in [\alpha_n, \beta_n],
\,
n=1,\ldots,N
\right),
$$
with respect to $\epsilon$, for two values of $N$.
We observe that for $\epsilon \lesssim \bar{\epsilon}=0.1$ results do not seem to depend on the value of $\epsilon$. This provides evidence that the current approach is robust with respect to the choice of $\epsilon$. We further observe that for all values of $\epsilon$ considered
$\mbox{\#}
{\rm Gal} (N=40) \gtrsim 0.85$ and 
$\mbox{\#}
{\rm Gal} (N=60) \gtrsim 0.90$. Therefore, our constrained formulation corrects the original formulation only for $10-15\%$ time steps.
For this reason, we envision that efficient implementations of the constrained ROM might be nearly as inexpensive as the Galerkin ROM. We further observe that  
$\mbox{\#} {\rm Gal}$ increases as $N$ increases: this can be explained by observing that the POD-Galerkin ROM becomes more and more accurate as $N$ increases, and thus requires less corrections.

\begin{figure}[h!]

\hspace{-80pt}
{  \includegraphics[width=1.30\textwidth]
 {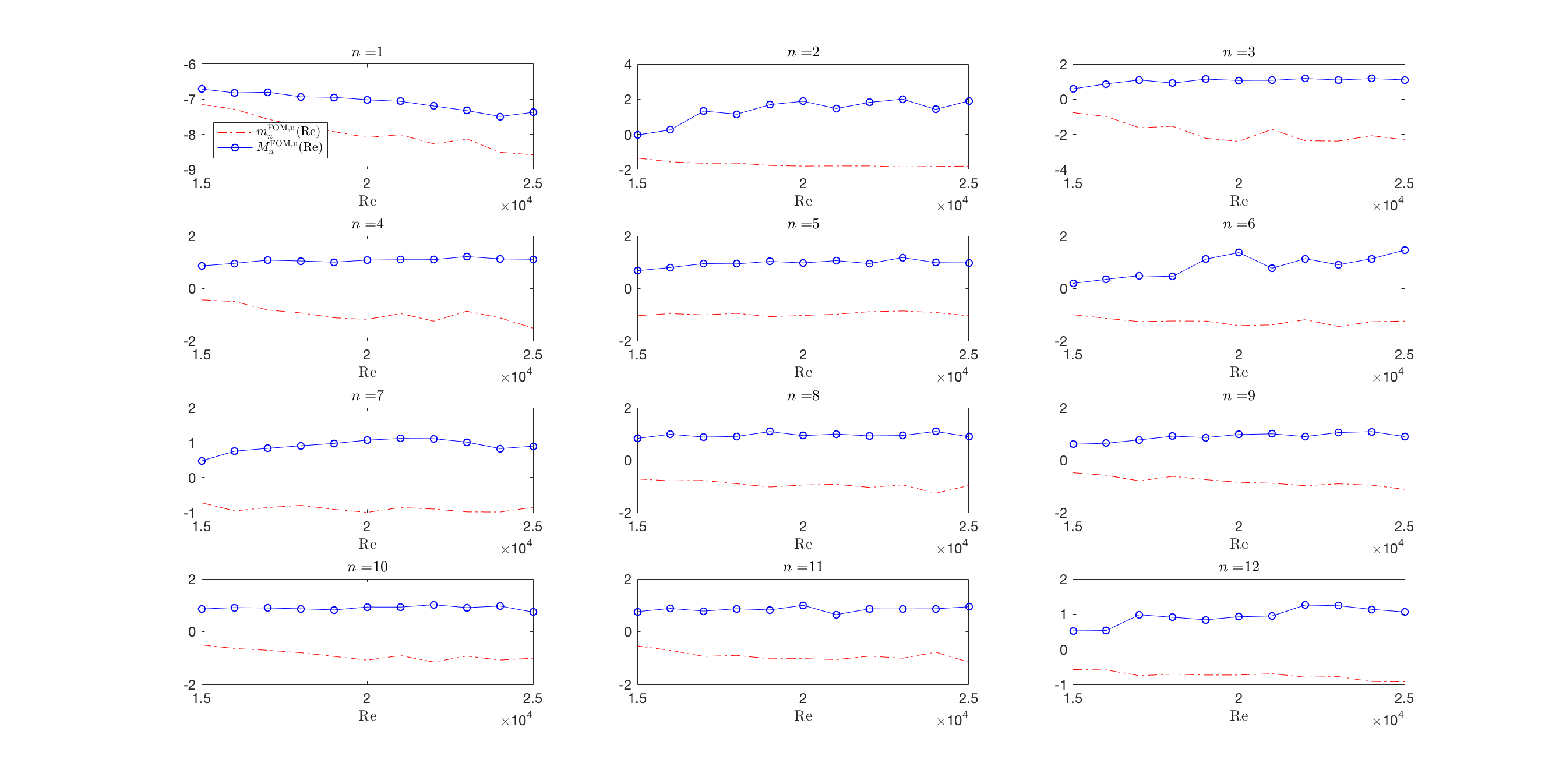}}

 \caption{The parametric problem; behavior of 
$m_n^{\rm FOM,u}({\rm Re})$ and $M_n^{\rm FOM,u}({\rm Re})$ with ${\rm Re}$.
The POD space is generated from the DNS data  for ${\rm Re}=20000$.
  }
 \label{fig:behavior_alpha_beta}
  \end{figure}

Figure \ref{fig:behavior_alpha_beta} investigates the behavior of 
$m_n^{\rm FOM,u}$ and $M_n^{\rm FOM,u}$  defined in \eqref{eq:golden_choice_parameters}
with respect to the Reynolds number ${\rm Re}$.
For this test, we consider the POD space associated with 
${\rm Re}=20000$ and the sampling times
$\{ t_{\rm s}^k = 500 + k  \}_{k=1}^{K=1000}$, and we show results for $n=1,\ldots,12$. Results suggest that  the sensitivity of 
$m_n^{\rm FOM,u}$ and $M_n^{\rm FOM,u}$ with ${\rm Re}$ 
are relatively modest if compared to  $M_n^{\rm FOM,u} - m_n^{\rm FOM,u}$. 

\section{On the problem of $p$-refinement }
\label{sec:p_refinement}
We here illustrate the  major issue associated with the combination of POD modes associated with different values of the parameter.
We here simulate the application of the first two iterations of the POD-$p$Greedy algorithm as proposed in \cite{haasdonk2008reduced}. 
In more detail, we consider the following test.
\begin{enumerate}
\item
Generate DNS data for ${\rm Re}=20000$, and use them to build the 
$N_1=60$-dimensional  POD space.
\item
Generate DNS data for ${\rm Re}=15000$, and build the 
$N_2=60$-dimensional  POD space for the set of snapshots
$\{  (\mathbb{I} - \Pi_{\mathcal{Z}^{\rm u}}^V)
\mathring{u}^k({\rm Re}=15000) 
\}_{k=1}^K$.
\item
Define $\mathcal{Z}^{\rm u}= {\rm span}\{ \zeta_n \}_{n=1}^{N_1+N_2}$ where
$\zeta_1,\ldots,\zeta_{N_1}$ are associated with Step 1 and 
$\zeta_{N_1+1},\ldots,\zeta_{N_2}$ are associated with Step 2.
\item
Perform a convergence study in $N$ for ${\rm Re}=15000$ and 
${\rm Re}=20000$ for both pure Galerkin and constrained Galerkin ($\epsilon=0.01$).
\end{enumerate}
We consider here $T_0=500$ and $T=1500$, $\{ t_{\rm s}^k=500+k \}_{k=1}^{K=1000}$.
We recall that for ${\rm Re}=15000$ (cf. section \ref{sec:reproduction_problem}) we were able to obtain  accurate ROMs 
  for $N \gtrsim 40$ both in terms of mean flow prediction and TKE.

\begin{figure}[h!]
\centering
 \subfloat[${\rm Re}=15000$]
{  \includegraphics[width=0.30\textwidth]
 {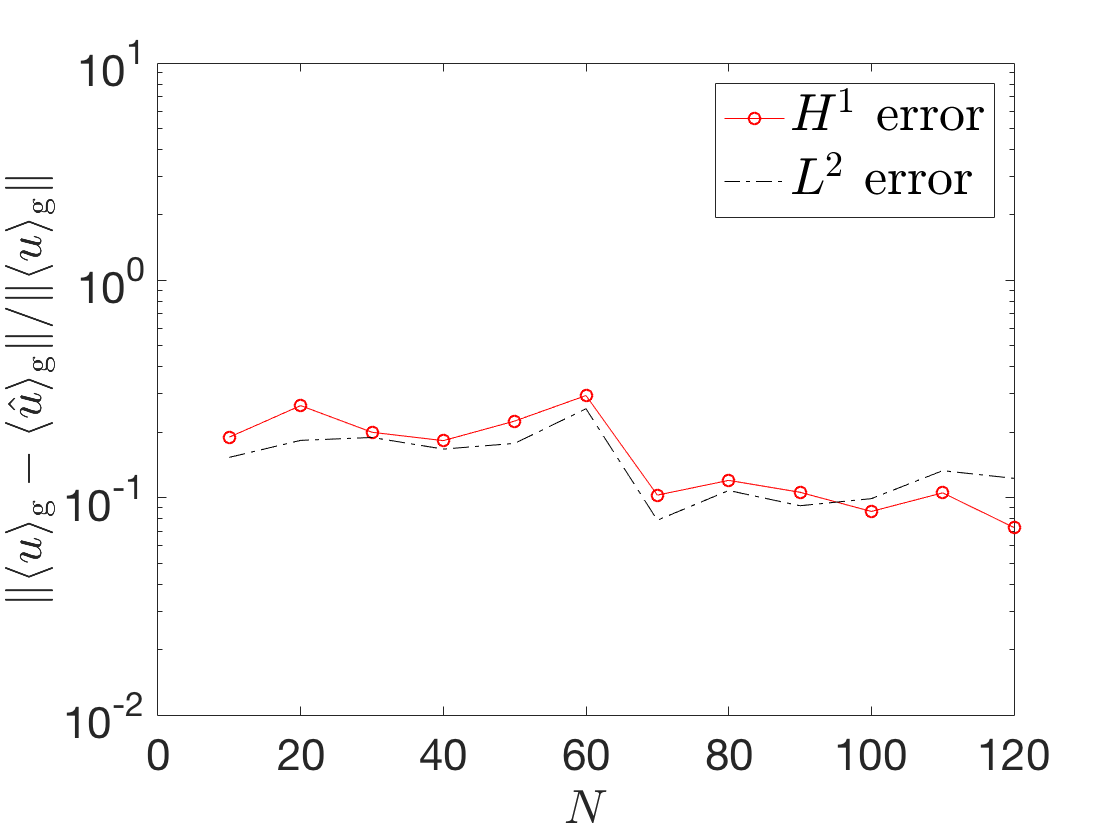}}
   ~~
 \subfloat[${\rm Re}=15000$]
{  \includegraphics[width=0.30\textwidth]
 {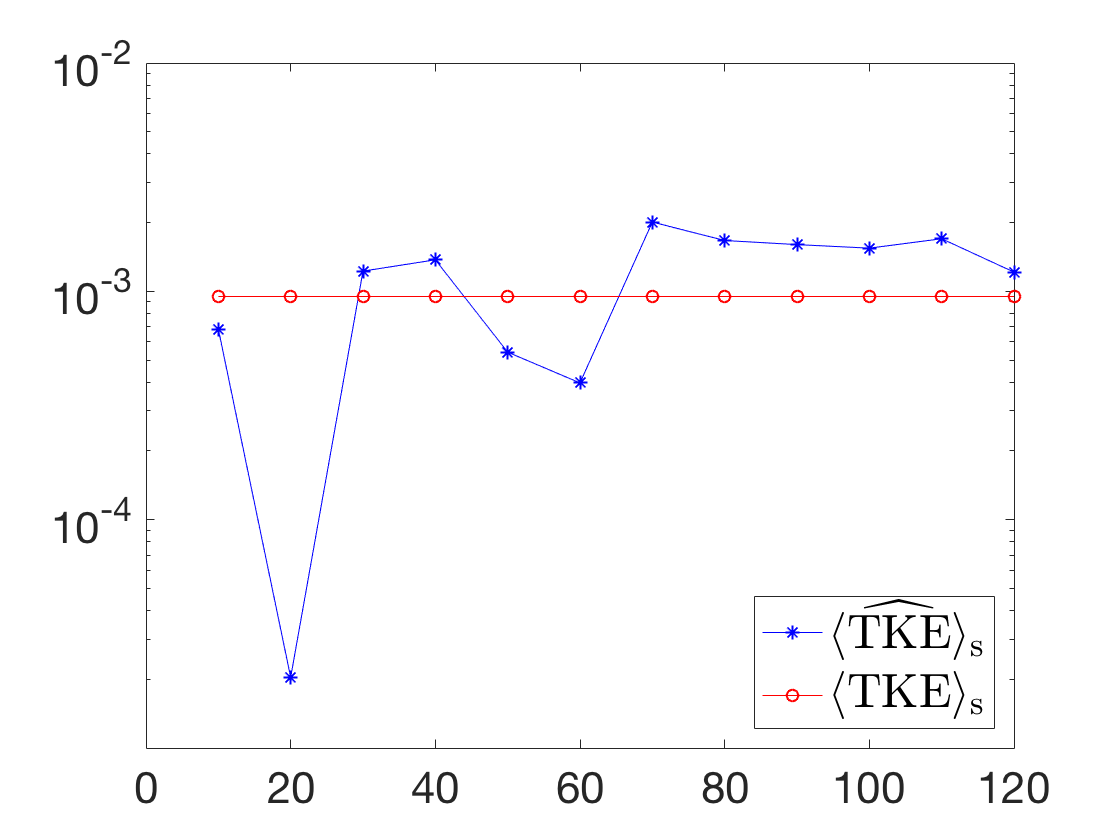}}

 \subfloat[${\rm Re}=20000$]
{  \includegraphics[width=0.30\textwidth]
 {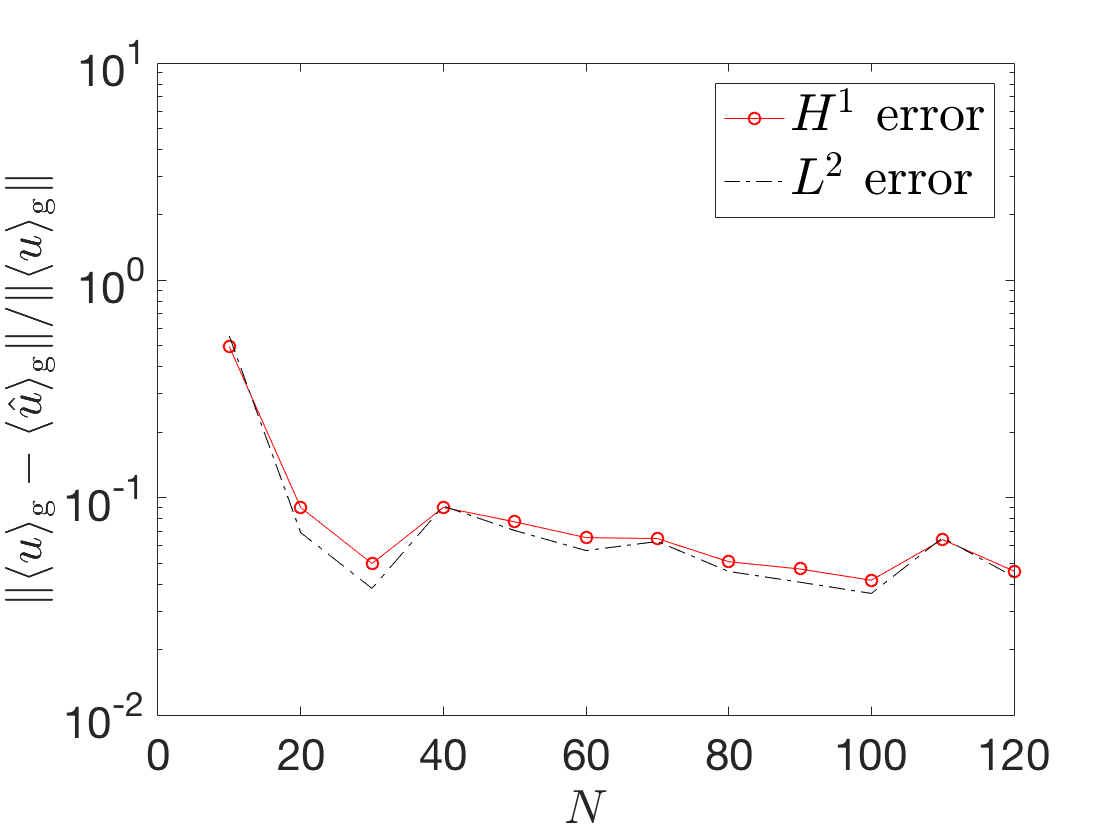}}
   ~~
 \subfloat[${\rm Re}=20000$]
{  \includegraphics[width=0.30\textwidth]
 {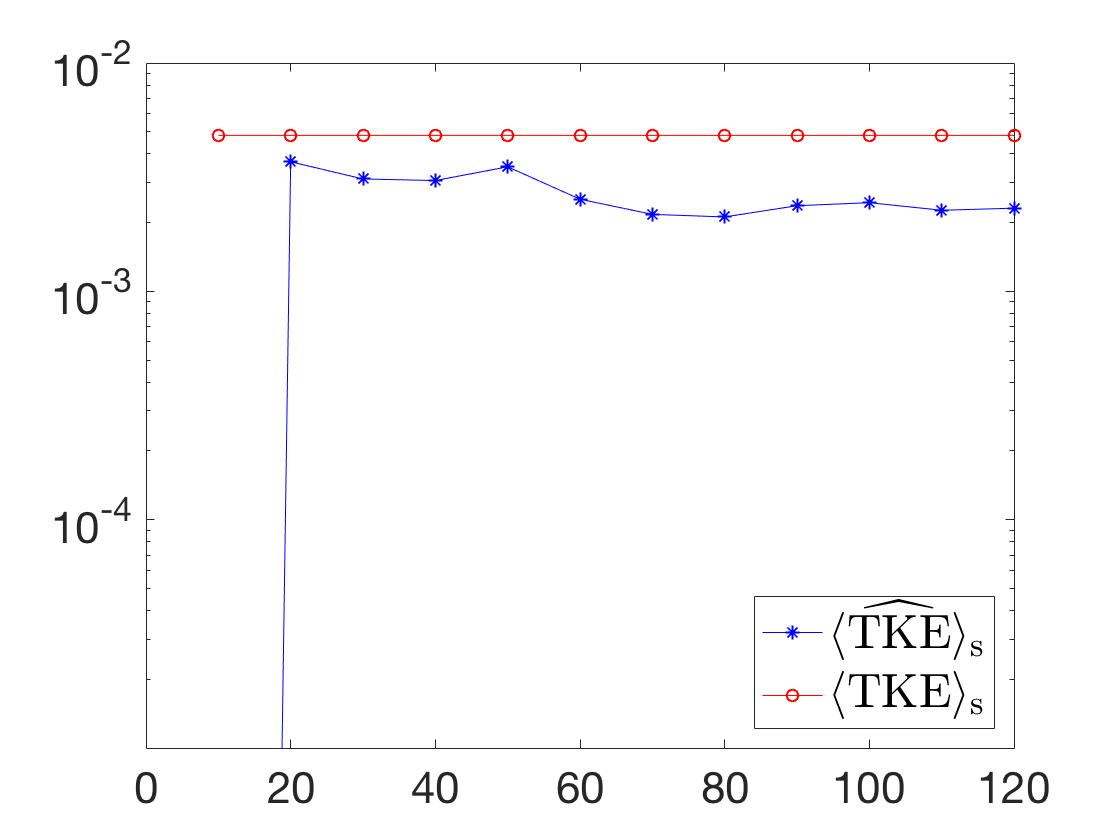}}
   
 \caption{The parametric problem;
on the problems of $p$-refinement.
Performance of constrained Galerkin
 ($\epsilon=0.01$, $N_1=60$, $N_2=60$).
  }
 \label{fig:pGreedy_failure}
  \end{figure}

\begin{figure}[h!]
\centering
 \subfloat[${\rm Re}=15000$, $N=120$]
{  \includegraphics[width=0.30\textwidth]
 {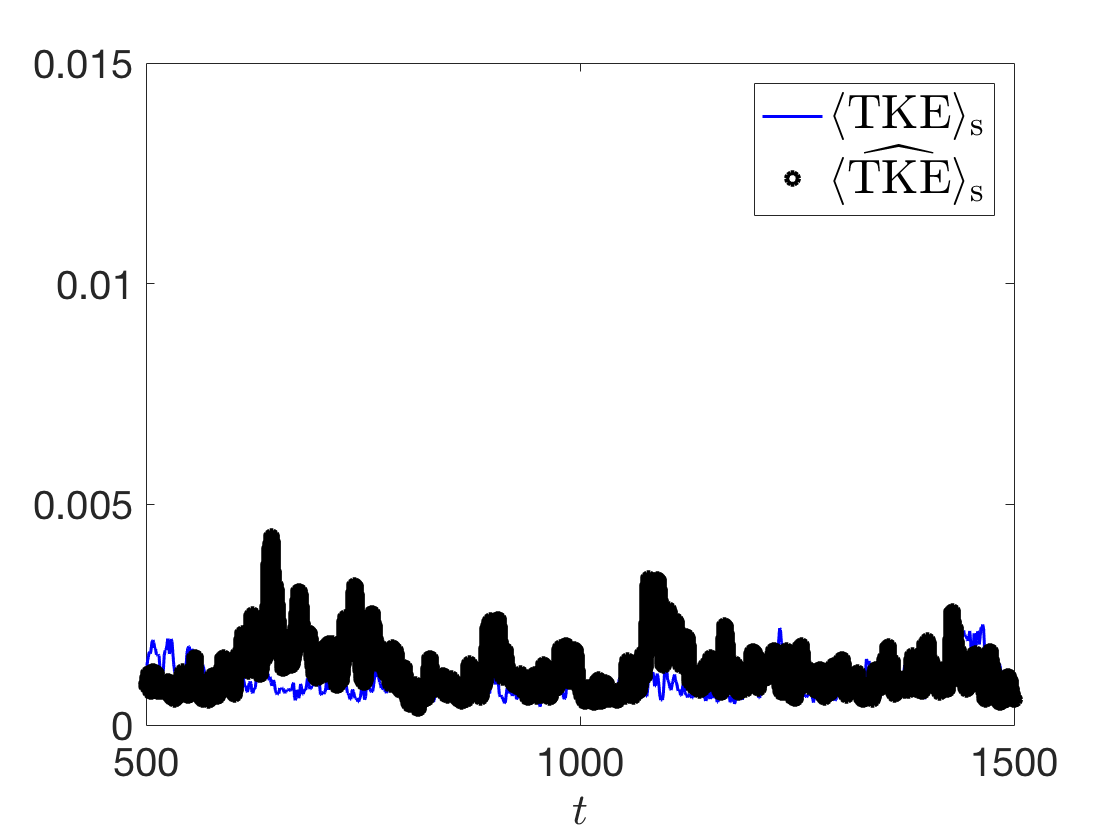}}
   ~~
 \subfloat[${\rm Re}=20000$, $N=120$]
{  \includegraphics[width=0.30\textwidth]
 {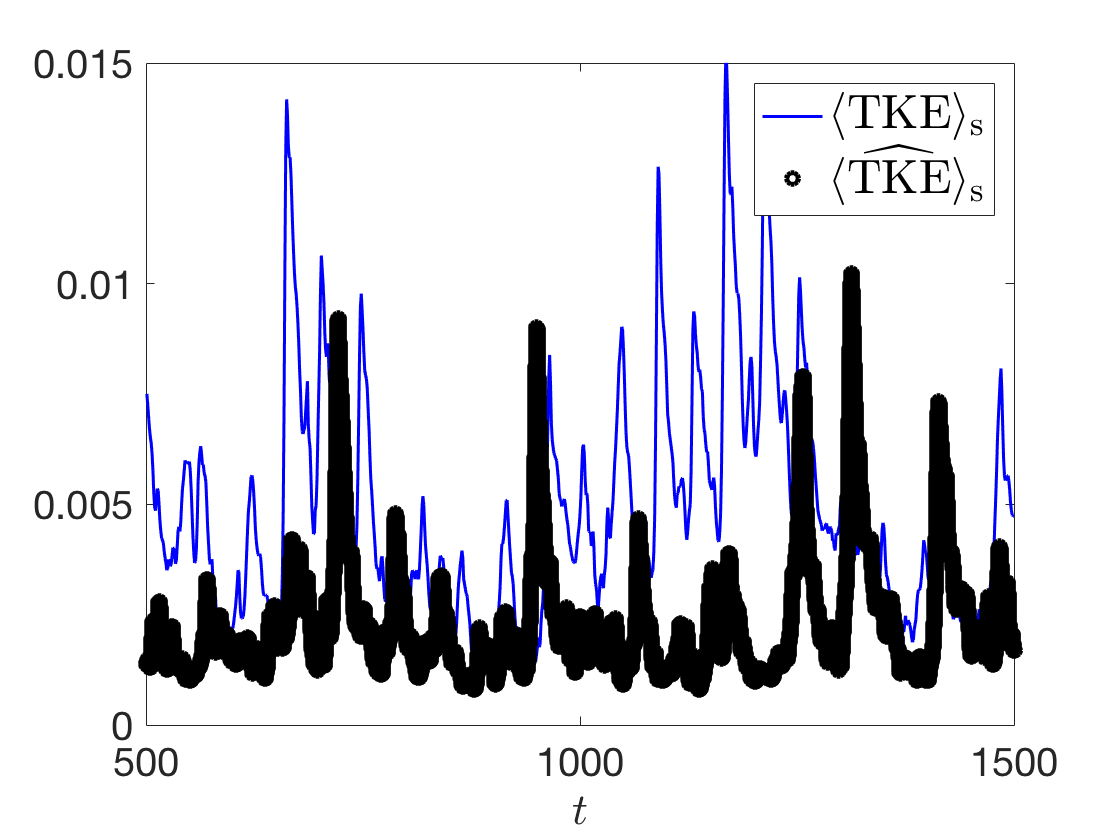}}

 \caption{
The parametric problem;
on the problems of $p$-refinement. Behavior of the TKE with time for constrained Galerkin 
 ($\epsilon=0.01$, $N_1=60$, $N_2=60$).
  }
 \label{fig:pGreedy_failure_time}
  \end{figure}

\begin{figure}[h!]
\centering
 \subfloat[${\rm Re}=15000$, $N=120$]
{  \includegraphics[width=0.30\textwidth]
 {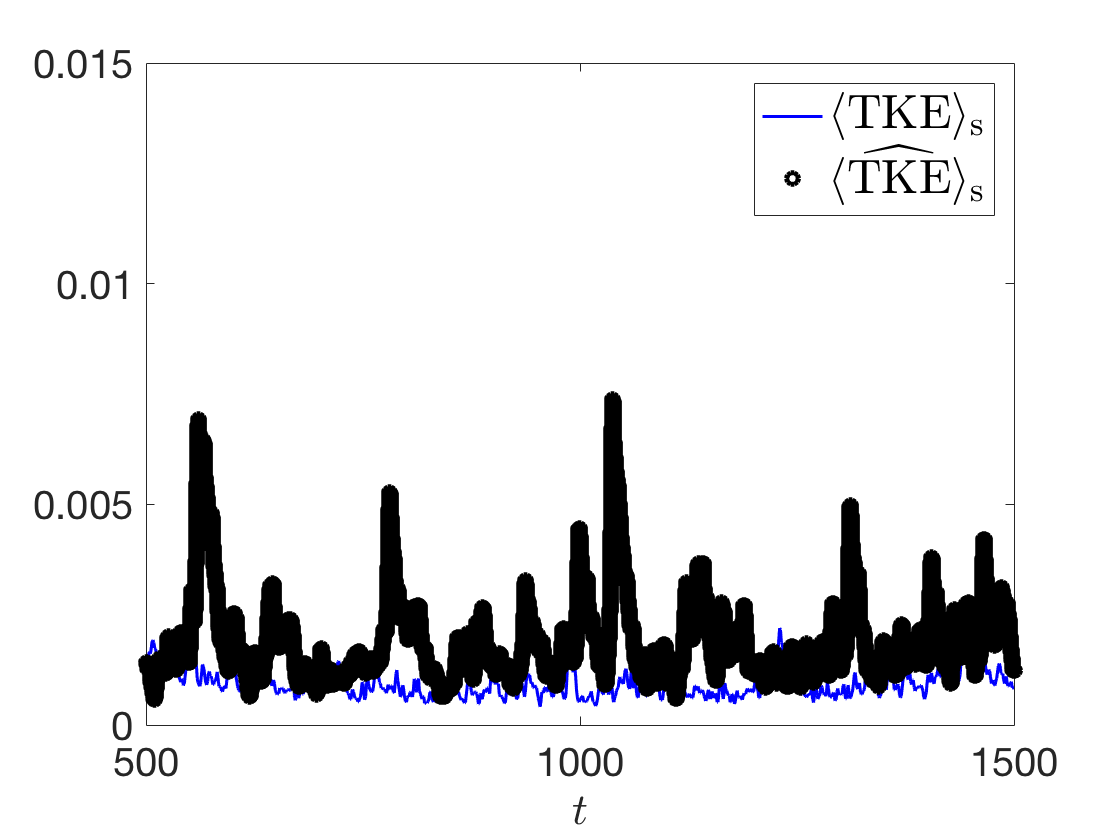}}
   ~~
 \subfloat[${\rm Re}=20000$, $N=120$]
{  \includegraphics[width=0.30\textwidth]
 {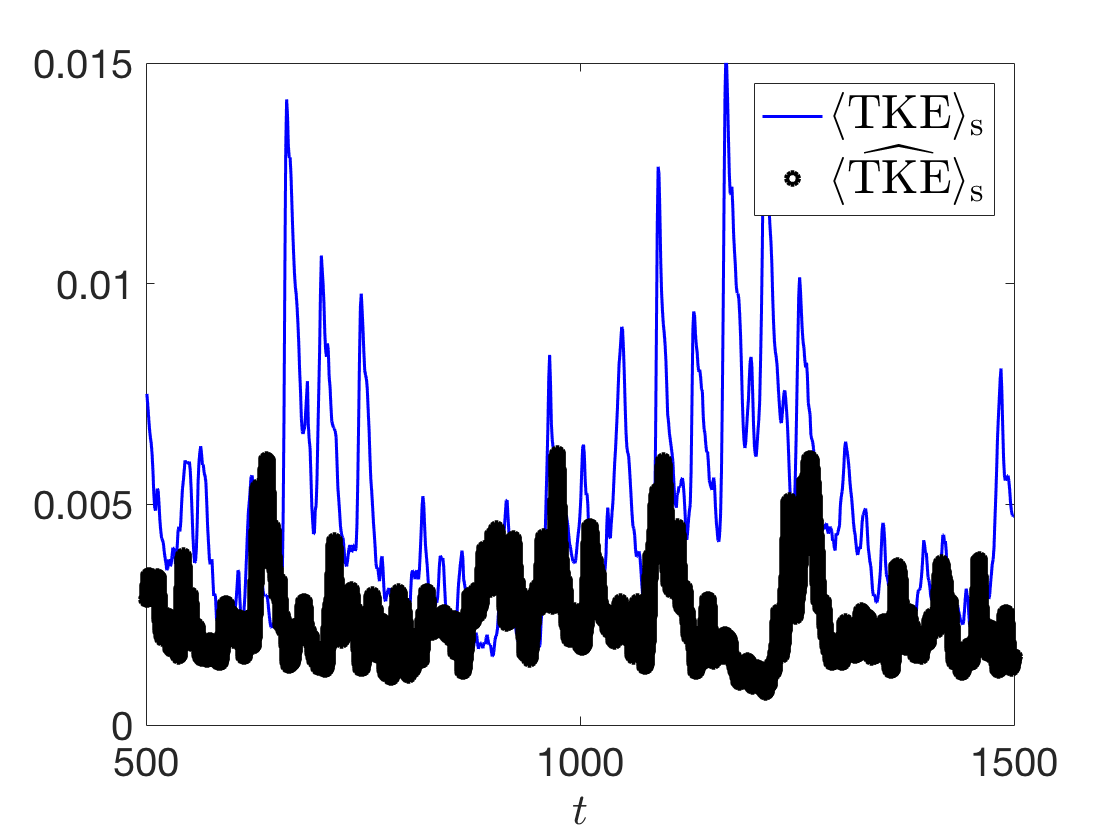}}

 \caption{
The parametric problem;
on the problems of $p$-refinement. Behavior of the TKE with time for unconstrained Galerkin 
 ($N_1=60$, $N_2=60$). 
 }
 \label{fig:pGreedy_failure_time_plain}
  \end{figure}

Figures \ref{fig:pGreedy_failure} and \ref{fig:pGreedy_failure_time} show the results of this test for the constrained-Galerkin ROM.
Figures \ref{fig:pGreedy_failure}(a) 
and (c) show  the behavior of the relative error in mean flow prediction 
for the constrained formulation, 
for ${\rm Re}=15000$ and ${\rm Re}=20000$, respectively.
We here compute lower and upper bounds $\{ \alpha_n \}_n$ and $\{ \beta_n \}_n$ using
\eqref{eq:choice_coefficients} with $\epsilon=0.01$.
Figures \ref{fig:pGreedy_failure}(b) 
and (d) show  the behavior of the mean TKE for 
the same values of the Reynolds number.
Similarly, 
Figures \ref{fig:pGreedy_failure_time}(a) and (b) show the behavior of the TKE in time for $N=120$.
Finally, Figures \ref{fig:pGreedy_failure_time_plain} (a) and (b) show the behavior of the TKE in time for $N=120$ for the unconstrained formulation.
Results --- especially for ${\rm Re}=15000$ ---
show the key issue of combining modes associated with different parameters.
For $N_2=60$ (and $N_1+N_2=120$), the error in mean flow prediction is roughly $10\%$, and we also significantly overestimate the mean  and the peaks of the TKE.
As expected, these issues are  even more severe for the unconstrained formulation: the behavior  with time of the  TKE predicted by the unconstrained ROM  is roughly the same for the two values of the Reynolds number considered.

We offer a physical explanation for the poor performance of POD-$p$Greedy. 
As observed in  Appendix  \ref{sec:lid_driven_appendix}, for sufficiently large values of ${\rm Re}$
 eddies are ejected into the core region of the cavity.   This instability is observed 
for  ${\rm Re}=20000$, but is not observed for  ${\rm Re}=15000$.
As a result, the ejection of the eddies into the core region of the cavity is well-represented by the POD space
associated with ${\rm Re}=20000$, and then, by construction, by  the final reduced space $\mathcal{Z}^{\rm u}$. The presence of modes associated with the core eddies makes the ROM more prone to show this instability even for values of the Reynolds number at which the full-order solution does not show it.

\section{Offline/online computational decomposition for the residual indicator}
\label{sec:error_indicator_offline_online}
We here describe the offline/online computational decomposition for the computation of the residual error indicator introduced in this paper. 
We  omit the subscript $\ell$ associated with the partition of the parameter domain to simplify notation.
We first introduce the Riesz representers:
\begin{equation}
\label{eq:riesz_representers}
\begin{array}{l}
(\xi_n^{\rm m},v )_V
=
(\zeta_n, v)_{L^2(\Omega)},
\quad
(\xi_n^{\rm a},v )_V
=
(\zeta_n, v)_{V},
\quad
(\xi_n^{\rm cg},v )_V
=
c(R_g, \zeta_n, v),
\\[3mm]
(\xi_{m,n}^{\rm c},v )_V
=
c(\zeta_n, \zeta_m, v),
\quad
(\xi_n^{\rm mg},v )_V
=
c(\zeta_n,R_g,  v),
\quad
(\xi_1^{\rm f},v )_V
=
 (R_g,  v)_V,
\\[3mm]
(\xi_2^{\rm f},v )_V
=
 c(R_g,R_g,  v),
 \\
\end{array}
\end{equation}
for $n=1,\ldots,N$ and for all $v \in V_{\rm div}$. Then, it is easy to verify that, if 
$w^j = \sum_{n=1}^N a_n^j \zeta_n$ for $j=0,\ldots,J$, we can rewrite $\langle R \rangle$ as follows:

\clearpage
\newpage
\begin{subequations}
\label{eq:time_avg_res_offline_online}
\begin{equation}
\begin{array}{ll}
\displaystyle{
\langle R  \rangle \left(
\{ w^j \}_{j=0}^J, v; \, {\rm Re}
\right)
=
}
&
\displaystyle{
\Big(
\sum_{n=1}^N
\Big(
\xi_n^{\rm m} \left(
\frac{ a_n^{J} - a_n^{J_0}   }{T-T_0}
\right)
+
\xi_n^{\rm a}
\left( \frac{1}{\rm Re} \bar{a}_n^+ \right)
+
\xi_n^{\rm cg}
 \bar{a}_n^+ 
 }
 \\[4mm]
 &
\displaystyle{
 +
\sum_{m=1}^N
\xi_{m,n}^{\rm c}
\bar{c}_{m,n}
+
\xi_n^{\rm mg}
 \bar{a}_n^-
\Big),
\, v
\Big)_V
}
\\[4mm]
&
\displaystyle{
+
\frac{1}{\rm Re}
(\xi_1^{\rm f}, v)_V
+
(\xi_2^{\rm f}, v)_V;
}
\\
\end{array}
\end{equation}
where
\begin{equation}
\bar{a}_n^+
=
\frac{\Delta t}{T-T_0}
\sum_{j=J_0+1}^J
\, a_n^j,
\quad
\bar{a}_n^-
=
\frac{\Delta t}{T-T_0}
\sum_{j=J_0}^{J-1}
\, a_n^j, \quad
\bar{c}_{m,n}
=
\frac{\Delta t}{T-T_0}
\, \sum_{j=J_0}^{J-1}
a_m^{j+1} a_n^j.
\end{equation}
\end{subequations}
Equation \eqref{eq:time_avg_res_offline_online} can be rewritten as 
$\langle R  \rangle \left(
\{ w^j \}_{j=0}^J, v; \, {\rm Re}
\right) = \sum_{m=1}^M \, \Theta_i( \{ \mathbf{a}^j  \}_j; {\rm Re} )  \tilde{\xi}_m$ where 
$M=N^2 + 3N+2$ and
$$
\begin{array}{l}
\displaystyle{
[\tilde{\xi}_1,\ldots,\tilde{\xi}_M]
=[\xi_1^{\rm m}, \ldots,\xi_N^{\rm m}, \xi_1^{\rm a}, \ldots,  \xi_N^{\rm a},\xi_1^{\rm cg}, \ldots,\xi_N^{\rm cg}, 
}
\\[2mm]
\hfill
\displaystyle{
\xi_{1,1}^{\rm c}, \ldots,  \xi_{N,N}^{\rm c}, 
\xi_1^{\rm mg}, \ldots,\xi_N^{\rm mg},\xi_1^{\rm f},\xi_2^{\rm f}] }
\\[4mm]
\displaystyle{
[ \Theta_1,\ldots, \Theta_M]
=
\Big[
\frac{a_1^{J} -a_1^{J_0}}{T-T_0},
\ldots,
\frac{a_N^{J} -a_N^{J_0}}{T-T_0},
\frac{ \bar{a}_1^+   }{\rm Re},
\ldots,
\frac{ \bar{a}_N^+   }{\rm Re},
\bar{a}_1^+,
\ldots,
\bar{a}_N^+,
}
\\[3mm]
\hfill
\displaystyle{
\bar{c}_{1,1},
\ldots,
\bar{c}_{N,N},
\bar{a}_1^-,
\ldots,
\bar{a}_N^-,
\frac{1}{\rm Re},
1
\Big].
}
\end{array}
$$
Therefore, recalling the Riesz representation theorem, we find 
\begin{equation}
\label{eq:final_error_estimator_offline_online}
\Delta^{\rm u}( \{ w^j \}_j; {\rm Re} )
=
\sqrt{   
\boldsymbol{\Theta}^T
\Sigma
\boldsymbol{\Theta} 
},
\qquad
\boldsymbol{\Theta}
=
\boldsymbol{\Theta}( \{ \mathbf{a}^j \}_j; {\rm Re}   ),
\end{equation}
where $\Sigma_{i,i'}=( \tilde{\xi}_i, \tilde{\xi}_{i'}  )_V$.
Equation \eqref{eq:final_error_estimator_offline_online} clarifies the offline/online decomposition:
during the offline stage, we compute the Riesz representers \eqref{eq:riesz_representers} ---
this corresponds to the solution to $M$ Stokes problems ---
and we assemble the matrix $\Sigma$; 
during the online stage, we compute the vector 
$\boldsymbol{\Theta}$ and we exploit \eqref{eq:final_error_estimator_offline_online} to compute the error estimator $\Delta^{\rm u}$.

\section*{Acknowledgements}
The authors thank Prof. Paul Fischer (UIUC), and Dr. Elia Merzari (Argonne National Lab) for their support with the software \texttt{Nek5000}.


\bibliographystyle{unsrt}

\bibliography{bib_folder/bib_MOR_RB,bib_folder/bib_MOR,bib_folder/bib_RBF,bib_folder/bib_data_assimilation,bib_folder/bib_statistics,bib_folder/bib_general_theory,bib_folder/bib_MOR_scrbe,bib_folder/validation,bib_folder/bib_my_publications,bib_folder/bib_inverse,bib_folder/bib_SHM,bib_folder/bib_optimization,bib_folder/bib_ROM_NS,bib_folder/bib_DMD}


\end{document}